\documentclass[10pt,a4paper]{article}
\usepackage{graphics}
\usepackage{amsmath,amssymb,amsthm,graphics,graphicx,epic,eepic,multicol,ascmac}
\usepackage[mathscr]{euscript}
\usepackage[all]{xy}
\usepackage{fancybox}
\usepackage{url}
\usepackage{ytableau}
\usepackage{here}
\usepackage[toc,page]{appendix}

\numberwithin{equation}{section}
\theoremstyle{theorem}
\newtheorem{thm}{Theorem}[section]
\newtheorem{prop}[thm]{Proposition}
\newtheorem{lem}[thm]{Lemma}

\newtheorem{conj}[thm]{Conjecture}

\theoremstyle{definition}
\newtheorem{defn}[thm]{Definition}
\newtheorem{ex}[thm]{Example}

\def\al{\alpha}

\def\wht(#1){\widehat{\ #1\ }}

\newcommand{\ch}{\mathrm{ch}}

\newcommand{\lbr}{\begin{bmatrix}}
\newcommand{\rbr}{\end{bmatrix}}

\newcommand{\cd}{commutative diagram }

\def\al{\alpha}

\def\beneme{\begin{enumerate}}
\def\beq{\begin{equation}}
\def\beqn{\begin{eqnarray}}
\def\beqnn{\begin{eqnarray*}}

\def\bfii0{{\bf i_0}}

\def\bbra#1,#2,#3{\left\{\begin{array}{c}\hspace{-5pt}
#1;#2\\ \hspace{-5pt}#3\end{array}\hspace{-5pt}\right\}}
\def\cd{\cdots}
\def\ci(#1,#2){c_{#1}^{(#2)}}
\def\Ci(#1,#2){C_{#1}^{(#2)}}
\def\mpp(#1,#2,#3){#1^{(#2)}_{#3}}
\def\bCi(#1,#2){\ovl C_{#1}^{(#2)}}
\def\ch(#1,#2){c_{#2,#1}^{-h_{#1}}}
\def\cc(#1,#2){c_{#2,#1}}

\def\di(#1,#2){D_{#1}^{(#2)}}
\def\dbi(#1,#2){\ovl D_{#1}^{(#2)}}

\def\eneme{\end{enumerate}}

\def\eeq{\end{equation}}
\def\eeqn{\end{eqnarray}}
\def\eeqnn{\end{eqnarray*}}

\def\gau#1,#2{\left[\begin{array}{c}\hspace{-5pt}#1\\
\hspace{-5pt}#2\end{array}\hspace{-5pt}\right]}

\def\ify{\infty}
\def\io{\iota}
\def\ji(#1,#2){j_{#1}^{(#2)}}

\def\km{k^{(-)}}

\def\lan{\langle}

\def\Lm{\Lambda}

\def\nd{\noindent}

\def\ovl{\overline}

\def\qed{\hfill\framebox[2mm]{}}
\def\QQ{\mathbb Q}

\def\ran{\rangle}

\def\TY(#1,#2,#3){#1^{(#2)}_{#3}}

\def\vp{\varphi}

\def\xxi(#1,#2,#3){\displaystyle {}^{#1}\Xi^{(#2)}_{#3}}
\def\xsi(#1,#2,#3){\displaystyle {}^{#1}\Sigma^{(#2)}_{#3}}
\def\xE(#1,#2,#3){\displaystyle {}^{#1}E_{#2}[#3]}
\def\xF(#1,#2){\displaystyle {}^{#1}F_{#2}}
\def\xx(#1,#2){\displaystyle {}^{#1}\Xi_{#2}}
\def\W1{W(\varpi_1)}

\def\ZZ{\mathbb Z}

\def\m@th{\mathsurround=0pt}
\def\fsquare(#1,#2){
\hbox{\vrule$\hskip-0.4pt\vcenter to #1{\normalbaselines\m@th
\hrule\vfil\hbox to #1{\hfill$\scriptstyle #2$\hfill}\vfil\hrule}$\hskip-0.4pt
\vrule}}

\newcommand{\ba}{\begin{array}}
\newcommand{\ea}{\end{array}}

\newcommand{\eq}{\begin{eqnarray}}
\newcommand{\eneq}{\end{eqnarray}}

\title{\textbf{\large{A conjecture on monomial realizations and polyhedral realizations for crystal bases
}}}
\author{\normalsize{YUKI KANAKUBO\thanks{Faculty of Basic Natural Science, Ibaraki University : {j\_chi\_sen\_you\_ky@eagle.sophia.ac.jp}}}
}
\date{}

\begin{document}

\maketitle
\vspace{-10pt}

\begin{abstract}
Crystal bases are powerful combinatorial tools in the representation theory of quantum groups $U_q(\mathfrak{g})$ for a symmetrizable Kac-Moody algebras $\mathfrak{g}$.
The polyhedral realizations are combinatorial descriptions of the crystal base $B(\infty)$ for Verma modules
in terms of the set of integer points of a polyhedral cone, which equals the string cone when $\mathfrak{g}$ is finite dimensional simple.
It is a fundamental and natural problem to find explicit forms of the polyhedral cone.
The monomial realization expresses crystal bases $B(\lambda)$ of integrable highest weight representations as Laurent monomials with double indexed variables.
In this paper, we give a conjecture between explicit forms of the polyhedral cones
and monomial realizations. We prove the conjecture is true when 
$\mathfrak{g}$ is a classical Lie algebra, a rank $2$ Kac-Moody algebra or a classical affine Lie algebra. 
\end{abstract}

\section{Introduction}

Crystal bases are introduced for combinatorial study of representations of quantum groups $U_q(\mathfrak{g})$ for symmetrizable Kac-Moody algebras $\mathfrak{g}$ over $\mathbb{C}$
and we can express them by using a bunch of combinatorial objects such as Young tableaux, Young walls, path models and so on \cite{JMMO,KN,Kang,KL,Lit94,Lit95}.
In this paper, we focus on two expressions, {\it monomial realizations} and {\it polyhedral realizations}.
The monomial realization is introduced in \cite{Kas03,Nak03}, which expresses crystal bases $B(\lambda)$ of integrable highest weight representations as Laurent monomials with double indexed variables.
The following is an example for $\mathfrak{g}=\mathfrak{s}\mathfrak{l}_3(\mathbb{C})$. The crystal base $B(\Lambda_1)$ is expressed as
\begin{equation}\label{intro-ex1}
X_{s,1}\rightarrow\frac{X_{s,2}}{X_{s+1,1}}\rightarrow \frac{1}{X_{s+1,2}}
\end{equation}
and $B(\Lambda_2)$ is expressed as
\begin{equation}\label{intro-ex2}
X_{s,2}\rightarrow\frac{X_{s+1,1}}{X_{s+1,2}}\rightarrow \frac{1}{X_{s+2,1}}.
\end{equation}

The polyhedral realization is invented in \cite{NZ}, which expresses elements in crystal bases $B(\infty)$
of Verma modules (or the negative part $U_q^-(\mathfrak{g})\subset U_q(\mathfrak{g})$) in terms of the set of integer points of a polyhedral cone. It is defined as the image ${\rm Im}(\Psi_{\iota})$ of an embedding $\Psi_{\iota}:B(\infty)\rightarrow \mathbb{Z}^{\infty}$
of crystals associated with an infinite sequence $\iota=(\cdots,i_3,i_2,i_1)$ of indices from $I$. Here,
$I=\{1,2,\cdots,n\}$ is an index set for simple roots of $\mathfrak{g}$.
The following is an example of ${\rm Im}(\Psi_{\iota})$ ($\cong B(\infty)$)
 for $\mathfrak{g}=\mathfrak{s}\mathfrak{l}_3(\mathbb{C})$ and $\iota=(\cdots,2,1,2,1)$:
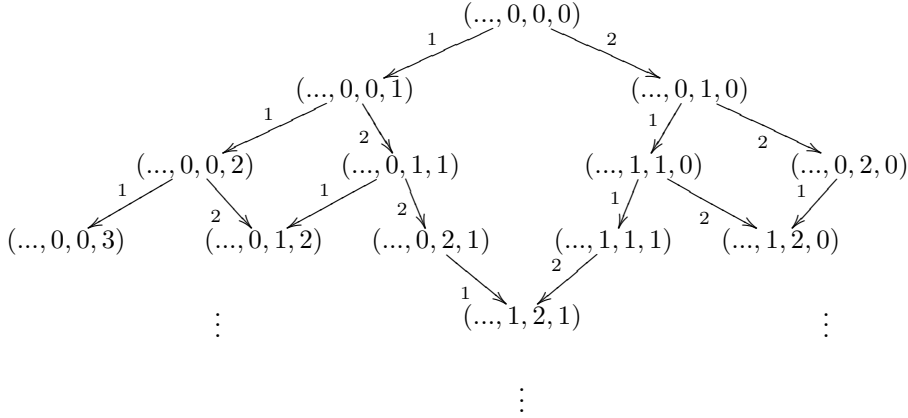
\[
\begin{xy}
(0,0)*{(...,0,0,0)}="0",
(22,-10) *{(...,0,1,0)}="2",
(-22,-10) *{(...,0,0,1)}="1",
(-43,-20) *{(...,0,0,2)}="11",
(-16,-20) *{(...,0,1,1)}="21",
(16,-20) *{(...,1,1,0)}="12",
(43,-20) *{(...,0,2,0)}="22",
(-60,-30) *{(...,0,0,3)}="111",
(-34,-30) *{(...,0,1,2)}="211",
(-12,-30) *{(...,0,2,1)}="221",
(12,-30) *{(...,1,1,1)}="112",
(34,-30) *{(...,1,2,0)}="122",
(0,-40) *{(...,1,2,1)}="1221",
(-40,-40) *{\vdots}="dot1",
(40,-40) *{\vdots}="dot2",
(0,-50) *{\vdots}="dot3",
\ar@{->} "0";"2"^{2}
\ar@{->} "0";"1"_{1}
\ar@{->} "1";"11"_{1}
\ar@{->} "1";"21"_{2}
\ar@{->} "2";"12"_{1}
\ar@{->} "2";"22"_{2}
\ar@{->} "11";"211"_{2}
\ar@{->} "21";"211"_{1}
\ar@{->} "11";"111"_{1}
\ar@{->} "21";"221"_{2}
\ar@{->} "12";"112"_{1}
\ar@{->} "22";"122"_{1}
\ar@{->} "12";"122"_{2}
\ar@{->} "221";"1221"_{1}
\ar@{->} "112";"1221"_{2}
\end{xy}
\]
Note that the set of elements appearing in the above graph coincides with the set of integer points in
a polyhedral cone:
\begin{equation}\label{intr-po}
{\rm Im}(\Psi_{\iota})=
\{
(\cdots,a_3,a_2,a_1)\in\mathbb{Z}^{\infty} |
a_1\geq0,\ a_2\geq a_3\geq0,\ a_k=0\ (k>3)
\}.
\end{equation}
It is natural problem to find an explicit form of the polyhedral cone. A lot of researchers
are working on this problem. It is known that if $\mathfrak{g}$ is a finite dimensional simple Lie algebra
then the cone coincides with a {\it string cone} in \cite{Lit98}. 
When $\mathfrak{g}$ is a finite dimensional simple Lie algebra and a specific sequence $\iota=(\cdots,n,\cdots,2,1,n,\cdots,2,1)$,
explicit forms of cones are given in \cite{H05,KS,NZ}.
For classical affine Lie algebras $\mathfrak{g}$ and same $\iota$, explicit forms are provided
in \cite{H13,NZ}. When $\mathfrak{g}$ is of type $A_n$ and $\iota$ is a specific one, it is shown that
the inequalities can be obtained from monomial realizations for a system of Demazure crystals \cite{Na}.
In our previous paper \cite{KaN20}, we combinatorially express inequalities defining the polyhedral cone in terms of column tableaux when $\mathfrak{g}$
is of type $A_n$, $B_n$, $C_n$ or $D_n$ and $\iota$ is adapted (see Definition \ref{adapt}).
For a classical affine type
$X=A^{(1)}_{n-1}$, $B^{(1)}_{n-1}$, $C^{(1)}_{n-1}$, $D^{(1)}_{n-1}$, $A^{(2)}_{2n-2}$, $A^{(2)}_{2n-3}$ or $D^{(2)}_{n}$
and an adapted sequence $\iota$, we give a combinatorial description of the inequalities in terms of extended Young diagrams
and Young walls of type $X^L$ \cite{Ka23a, Ka23b, Ka24a}. The notation $X^L$ implies the Langlands dual type for $X$ (see (\ref{Langtype})).
Note that column tableaux, extended Young diagrams
and Young walls are introduced for the combinatorial study of fundamental representations of quantum groups $U_q(\mathfrak{g})$
\cite{JMMO,KN,Kang}.
From these results, one can expect that
the inequalities defining the polyhedral cone are expressed by some combinatorial objects deeply related to 
fundamental representations of quantum groups for general symmetrizable Kac-Moody algebras $\mathfrak{g}$.
Based on this philosophy, we focus on the monomial realizations as the combinatorial objects.

The main purpose of this paper is to give a conjecture between explicit forms of polyhedral realizations 
and monomial realizations (Conjecture \ref{mainconj}).
The detail is as follows: We fix an adapted sequence $\iota$ and
consider the union of monomial realizations $\mathcal{M}_{s,k,\iota}$ with $s\in\mathbb{Z}_{\geq1}$, $k\in I$ for the Langlands dual algebra $\mathfrak{g}^L$.
Using the tropicalization map (subsection \ref{tropsub}), one considers the subset of $\mathbb{Z}^{\infty}$:
\begin{equation}\label{intro-eq1}
\{
\mathbf{a}=(a_{m,j})_{m\in\mathbb{Z}_{\geq1},j\in I}\in\mathbb{Z}^{\infty} | Trop(M)(\mathbf{a})\geq0\quad\text{for all }M\in\bigcup_{s\in\mathbb{Z}_{\geq1},k\in I}\mathcal{M}_{s,k,\iota}
\}.
\end{equation}
Here, we identify $(a_j)_{j\in\mathbb{Z}_{\geq1}}\in\mathbb{Z}^{\infty}$ with $(a_{m,j})_{m\in\mathbb{Z}_{\geq1},j\in I}$ by a rule of subsection \ref{s-d}.
The conjecture claims the above set coincides with ${\rm Im}(\Psi_{\iota})$. For instance, by tropicalizing monomials in (\ref{intro-ex1}), (\ref{intro-ex2}), one obtains
a system of inequalities
\[
a_{s,1}\geq0,\ a_{s,2}-a_{s+1,1}\geq0,\ -a_{s+1,2}\geq0,\quad
a_{s,2}\geq0,\ a_{s+1,1}-a_{s+1,2}\geq0,\ -a_{s+2,1}\geq0\quad (s\in\mathbb{Z}_{\geq1}).
\]
By $a_{s,1}\geq0$, $-a_{s+2,1}\geq0$ and $a_{s,2}\geq0$, $-a_{s+1,2}\geq0$, we see that $a_{m+2,1}=a_{m+1,2}=0$ for all $m\geq1$.
Simplifying other inequalities, the set (\ref{intro-eq1}) is equal to
\[
\{
\mathbf{a}=(a_{m,j})_{m\in\mathbb{Z}_{\geq1},j\in I}\in\mathbb{Z}^{\infty} |
a_{1,2}\geq a_{2,1}\geq0,\ a_{1,1}\geq0,\ 
a_{m+2,1}=a_{m+1,2}=0\ (m\in\mathbb{Z}_{\geq1})
\},
\]
which coincides with ${\rm Im}(\Psi_{\iota})$ in (\ref{intr-po}) as set.

We will prove the conjecture is true when 
\begin{itemize}
\item
$\mathfrak{g}$ is a finite dimensional simple Lie algebra of type $A_n,B_n,C_n$ or $D_n$,
\item $\mathfrak{g}$ is a rank $2$ Kac-Moody algebra.
\item $\mathfrak{g}$ is an affine Lie algebra of type $A^{(1)}_{n-1},B^{(1)}_{n-1},C^{(1)}_{n-1},D^{(1)}_{n-1},A^{(2)}_{2n-2},A^{(2)}_{2n-3}$ or $D^{(2)}_n$. 
\end{itemize}
For the proof, we use our previous expression of inequalities given in 
\cite{KaN20,Ka23a, Ka24a,Ka24b}. 

The organization
of paper is as follows:
In Sect.2, we recall the definition of crystals and polyhedral realizations. In particular, we recall
procedures to compute the explicit forms of ${\rm Im}(\Psi_{\iota})$.
In Sect.3, we review monomial realizations of crystal bases.
Sect.4 is devoted to give our conjecture and main theorems.
In Sect.5, we prove one of theorems, which gives a sufficient condition for the conjecture.
We show Theorem \ref{mainthm2} in Sect.6, 7 and 8
that claims the conjecture is true for specific types.

%We assume (\ref{mainassump}) and
%let us prove
%\[
%\{
%S_{j_1}'\cdots S_{j_m}'x_{s,k} |m\geq0,\ j_1,\cdots,j_m\in\mathbb{Z}_{\geq1} 
%\}=Trop(\mathcal{M}_{s,k,\iota}).
%\]
%Here, $\mathcal{M}_{s,k,\iota}$ is a monomial realization for $B(\Lambda_k)$ with highest weight vector $X_{s,k}$
%associated with an adapted sequence $\iota$.

\vspace{2mm}

\nd
{\bf Acknowledgements}
This work was supported by JSPS KAKENHI Grant Number JP24K22825.

\section{Crystals and polyhedral realizations for $B(\infty)$}

\subsection{Notation}

We set $I:=\{1,2,\cdots,n\}$ with $n\in \mathbb{Z}_{\geq1}$.
Let $\mathfrak{g}$ be a symmetrizable Kac-Moody algebra over $\mathbb{C}$ with a generalized Cartan matrix $(a_{i,j})_{i,j\in I}$, 
Cartan subalgebra $\mathfrak{h}$,
weight lattice $P\subset \mathfrak{h}^*$, set of
simple roots $\{\alpha_i\}_{i\in I}$ and set of simple coroots $\{h_i\}_{i\in I}$, Weyl group $W$. 
Let $P^+:=\{\lambda\in P |\lan h_i,\lambda \ran\mathbb{Z}_{\geq0},\ \text{for all }i\in I\}$ be the set of dominant integral weights,
$\lan ,\ran:\mathfrak{h}\times \mathfrak{h}^*\rightarrow \mathbb{C}$ denote the canonical pairing
and $\Lambda_i\in P^+$ be the $i$-th fundamental weight for $i\in I$.
We obtain $\lan h_i,\alpha_j\ran=a_{i,j}$ and $\lan h_i,\Lambda_j\ran=\delta_{i,j}$. 
Let $U_q(\mathfrak{g})$
be the quantized universal enveloping algebra of $\mathfrak{g}$ with indeterminant $q$, which is an associative $\mathbb{C}(q)$-algebra
generated by $e_i$, $f_i$ ($i\in I$) and $q^h$ ($h\in P^*=\{h\in\mathfrak{h}|\lan h,P\ran\subset\mathbb{Z}\}$).
It has a subalgebra $U_q^-(\mathfrak{g})$ generated by $f_i$ ($i\in I$).

It is known that for $\lambda\in P^+$, every integrable highest weight representation $V(\lambda)$ has a crystal base $(L(\lambda),B(\lambda))$
and $U_q^-(\mathfrak{g})$ also has a crystal base $(L(\infty),B(\infty))$ (\cite{Kas90,Kas91,Lus}).
For $X=A_n,B_n,C_n$ or $D_n$, we define the Langlands dual type $X^L$ as
$X^L=X$ ($X=A_n$, $D_n$) and $(B_n)^L=C_n$, $(C_n)^L=B_n$.
We define
\begin{equation}\label{Langtype}
X^L=X\ (X=A^{(1)}_{n-1},D^{(1)}_{n-1}),\ (C^{(1)}_{n-1})^L=D^{(2)}_{n},\ (D^{(2)}_{n})^L=C^{(1)}_{n-1},\ 
(A^{(2)}_{2n-3})^L=B^{(1)}_{n-1},\ (B^{(1)}_{n-1})^L=A^{(2)}_{2n-3}.
\end{equation}
We also define $(A^{(2)}_{2n-2})^L=A^{(2)\dagger}_{2n-2}$ and $(A^{(2)\dagger}_{2n-2})^L=A^{(2)}_{2n-2}$. Note that $(X^L)^L=X$.
We often omit the rank and simply write $A_n$, $B_n$, $C_n$, $D_n$ as $A,B,C,D$ for simplicity.
Let $\mathfrak{g}^L$ denote the Kac-Moody algebra whose generalized Cartan matrix
is the transposed matrix of $(a_{i,j})_{i,j\in I}$. For two integers $m,l$ with $m\leq l$, we set $[m,l]:=\{m,m+1,\cdots,l-1,l\}$.

The numbering of Dynkin diagrams and affine Dynkin diagrams are as follows:
\[
\begin{xy}
(-8,0) *{A_{n} : }="A1",
(0,0) *{\bullet}="1",
(0,-3) *{1}="1a",
(10,0)*{\bullet}="2",
(10,-3)*{2}="2a",
(20,0)*{\ \cdots\ }="3",
(30,0)*{\bullet}="4",
(30,-3)*{n-1}="4a",
(40,0)*{\bullet}="5",
(40,-3)*{n}="5a",
(52,0) *{B_{n} : }="C1",
(60,0) *{\bullet}="11",
(60,-3) *{1}="11a",
%(85,5) *{>}=">",
(70,0)*{\bullet}="22",
(70,-3)*{2}="22a",
(80,0)*{\ \cdots\ }="33",
(90,0)*{\bullet}="44",
(90,-3)*{n-1}="44a",
%(115,5)*{<}="<",
(100,0)*{\bullet}="55",
(100,-3)*{n}="55a",
\ar@{-} "1";"2"^{}
\ar@{-} "2";"3"^{}
\ar@{-} "3";"4"^{}
\ar@{-} "4";"5"^{}
\ar@{-} "11";"22"^{}
\ar@{-} "22";"33"^{}
\ar@{-} "33";"44"^{}
\ar@{=>} "44";"55"^{}
\end{xy}
\]
\[
\begin{xy}
(-8,5) *{C_{n} : }="A2",
(0,5) *{\bullet}="1",
(0,2) *{1}="1a",
%(5,5) *{>}=">",
(10,5)*{\bullet}="2",
(10,2)*{2}="2a",
(20,5)*{\ \cdots\ }="3",
(30,5)*{\bullet}="4",
(30,2)*{n-1}="4a",
%(35,5)*{>}=">",
(40,5)*{\bullet}="5",
(40,2)*{n}="5a",
(52,5) *{D_{n} : }="D2",
(60,5) *{\bullet}="11",
(60,2) *{1}="11a",
%(85,5) *{<}="<",
(70,5)*{\bullet}="22",
(70,2)*{2}="22a",
(80,5)*{\ \cdots\ }="33",
(90,5)*{\bullet}="44",
(88,2)*{n-2}="44a",
%(115,5)*{>}=">",
(100,0)*{\bullet}="55a",
(100,10)*{\bullet}="55b",
(100,14)*{n-1}="55bb",
(100,-4)*{n}="55aa",
\ar@{-} "1";"2"^{}
\ar@{-} "2";"3"^{}
\ar@{-} "3";"4"^{}
\ar@{<=} "4";"5"^{}
\ar@{-} "11";"22"^{}
\ar@{-} "22";"33"^{}
\ar@{-} "33";"44"^{}
\ar@{-} "44";"55a"^{}
\ar@{-} "44";"55b"^{}
\end{xy}
\]
\[
\begin{xy}
(52,5) *{G_{2} : }="D2",
(60,5) *{\bullet}="11",
(60,5.8) *{}="11u",
(60,5.1) *{}="11c",
(60,4.5) *{}="11d",
(60,2) *{1}="11a",
%(85,5) *{<}="<",
(70,5)*{\bullet}="22",
(67.5,5.8)*{}="22u",
(70,5.1)*{}="22c",
(67.5,4.5)*{}="22d",
(70,2)*{2}="22a",
(68,5.1)*{>}="ya",
\ar@{-} "11d";"22d"^{}
\ar@{-} "11c";"22c"^{}
\ar@{-} "11u";"22u"^{}
\end{xy}
\]
\[
\begin{xy}
(-8,0) *{A_{1}^{(1)} : }="A1",
(0,0) *{\bullet}="1",
(0,-3) *{1}="1a",
(10,0)*{\bullet}="2",
(10,-3) *{2}="2a",
\ar@{<=>} "1";"2"^{}
\end{xy}
\]
\[
\begin{xy}
(-15,5) *{A_{n-1}^{(1)} \ (n\geq 3) : }="A1",
(20,10) *{\bullet}="n",
(20,13) *{n}="na",
(0,0) *{\bullet}="1",
(0,-3) *{1}="1a",
(10,0)*{\bullet}="2",
(10,-3)*{2}="2a",
(20,0)*{\ \cdots\ }="3",
(30,0)*{\bullet}="4",
(30,-3)*{n-2}="4a",
(40,0)*{\bullet}="5",
(40,-3)*{n-1}="5a",
(65,5) *{B_{n-1}^{(1)} \ (n\geq 4) : }="B1",
(80,10) *{\bullet}="11",
(80,0) *{\bullet}="111",
(80,13) *{1}="111a",
(80,-3) *{2}="11a",
%(85,5) *{>}=">",
(90,5)*{\bullet}="22",
(90,2)*{3}="22a",
(100,5)*{\ \cdots\ }="33",
(110,5)*{\bullet}="44",
(110,2)*{n-1}="44a",
%(115,5)*{<}="<",
(120,5)*{\bullet}="55",
(120,2)*{n}="55a",
\ar@{-} "1";"2"^{}
\ar@{-} "2";"3"^{}
\ar@{-} "3";"4"^{}
\ar@{-} "4";"5"^{}
\ar@{-} "1";"n"^{}
\ar@{-} "n";"5"^{}
\ar@{-} "11";"22"^{}
\ar@{-} "111";"22"^{}
\ar@{-} "22";"33"^{}
\ar@{-} "33";"44"^{}
\ar@{=>} "44";"55"^{}
\end{xy}
\]
\[
\begin{xy}
(-15,5) *{C_{n-1}^{(1)} \ (n\geq 3) : }="C1",
(0,5) *{\bullet}="11",
(0,2) *{1}="11a",
%(85,5) *{>}=">",
(10,5)*{\bullet}="22",
(10,2)*{2}="22a",
(20,5)*{\ \cdots\ }="33",
(30,5)*{\bullet}="44",
(30,2)*{n-1}="44a",
%(115,5)*{<}="<",
(40,5)*{\bullet}="55",
(40,2)*{n}="55a",
(65,5) *{D_{n-1}^{(1)} \ (n\geq 5) : }="D1",
(80,10) *{\bullet}="111",
(80,0) *{\bullet}="1111",
(80,13) *{1}="111a1",
(80,-3) *{2}="11a1",
%(85,5) *{>}=">",
(90,5)*{\bullet}="221",
(90,2)*{3}="22a1",
(100,5)*{\ \cdots\ }="331",
(110,5)*{\bullet}="441",
(110,2)*{n-2}="44a1",
%(115,5)*{<}="<",
(120,10)*{\bullet}="551",
(120,0)*{\bullet}="5511",
(120,12)*{n-1}="55a1",
(120,-2)*{n}="55a11",
\ar@{=>} "11";"22"^{}
\ar@{-} "22";"33"^{}
\ar@{-} "33";"44"^{}
\ar@{<=} "44";"55"^{}
\ar@{-} "111";"221"^{}
\ar@{-} "1111";"221"^{}
\ar@{-} "221";"331"^{}
\ar@{-} "331";"441"^{}
\ar@{-} "441";"551"^{}
\ar@{-} "441";"5511"^{}
\end{xy}
\]
\[
\begin{xy}
(-15,5) *{A_{2n-2}^{(2)}\ (n\geq 3) : }="A2",
(0,5) *{\bullet}="1",
(0,2) *{1}="1a",
%(5,5) *{>}=">",
(10,5)*{\bullet}="2",
(10,2)*{2}="2a",
(20,5)*{\ \cdots\ }="3",
(30,5)*{\bullet}="4",
(30,2)*{n-1}="4a",
%(35,5)*{>}=">",
(40,5)*{\bullet}="5",
(40,2)*{n}="5a",
(65,5) *{D_{n}^{(2)} \ (n\geq 3) : }="D2",
(80,5) *{\bullet}="11",
(80,2) *{1}="11a",
%(85,5) *{<}="<",
(90,5)*{\bullet}="22",
(90,2)*{2}="22a",
(100,5)*{\ \cdots\ }="33",
(110,5)*{\bullet}="44",
(110,2)*{n-1}="44a",
%(115,5)*{>}=">",
(120,5)*{\bullet}="55",
(120,2)*{n}="55a",
\ar@{=>} "1";"2"^{}
\ar@{-} "2";"3"^{}
\ar@{-} "3";"4"^{}
\ar@{=>} "4";"5"^{}
\ar@{<=} "11";"22"^{}
\ar@{-} "22";"33"^{}
\ar@{-} "33";"44"^{}
\ar@{=>} "44";"55"^{}
\end{xy}
\]
\[
\begin{xy}
(65,5) *{A_{2n-3}^{(2)} \ (n\geq 4) : }="A2",
(80,10) *{\bullet}="11",
(80,0) *{\bullet}="111",
(80,13) *{1}="111a",
(80,-3) *{2}="11a",
%(85,5) *{>}=">",
(90,5)*{\bullet}="22",
(90,2)*{3}="22a",
(100,5)*{\ \cdots\ }="33",
(110,5)*{\bullet}="44",
(110,2)*{n-1}="44a",
%(115,5)*{<}="<",
(120,5)*{\bullet}="55",
(120,2)*{n}="55a",
\ar@{-} "11";"22"^{}
\ar@{-} "111";"22"^{}
\ar@{-} "22";"33"^{}
\ar@{-} "33";"44"^{}
\ar@{<=} "44";"55"^{}
\end{xy}
\]
Replacing our numbering $1,2,\cdots,n-1,n,n+1$ of $A_{2n}^{(2)}$ with $n,n-1,\cdots,1,0$, one obtains the numbering in \cite{Kang}.
We define the type $A_{2n-2}^{(2)\dagger}$ by the following diagram whose numbering is same as \cite{Kang}:
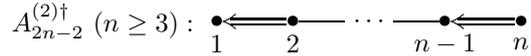
\[
\begin{xy}
(-15,5) *{A_{2n-2}^{(2)\dagger}\ (n\geq 3) : }="A2",
(0,5) *{\bullet}="1",
(0,2) *{1}="1a",
%(5,5) *{>}=">",
(10,5)*{\bullet}="2",
(10,2)*{2}="2a",
(20,5)*{\ \cdots\ }="3",
(30,5)*{\bullet}="4",
(30,2)*{n-1}="4a",
%(35,5)*{>}=">",
(40,5)*{\bullet}="5",
(40,2)*{n}="5a",
\ar@{<=} "1";"2"^{}
\ar@{-} "2";"3"^{}
\ar@{-} "3";"4"^{}
\ar@{<=} "4";"5"^{}
\end{xy}
\]

\subsection{Crystals}

We briefly review the crystals.

\begin{defn}
Let $B$ be a set and we suppose that there are maps
${\rm wt}:B\rightarrow P$, $\varepsilon_i:B\rightarrow\mathbb{Z}\sqcup\{-\infty\}$, $\varphi_i:B\rightarrow\mathbb{Z}\sqcup\{-\infty\}$,
$\tilde{f}_i:B\rightarrow B\sqcup\{0\}$ and $\tilde{e}_i:B\rightarrow B\sqcup\{0\}$ for $i\in I$.
Here, $0$ and $-\infty$ are additional elements.
When the following conditions hold, the set $B$ together with these maps is called a {\it crystal}: For $b,b'\in B$,
\begin{enumerate}
\item[(1)]${\rm wt}(\tilde{e}_kb)={\rm wt}(b)+\alpha_k$ if $\tilde{e}_k(b)\in B$,
\qquad${\rm wt}(\tilde{f}_kb)={\rm wt}(b)-\alpha_k$ if $\tilde{f}_k(b)\in B$,
\item[(2)]$\varphi_k(b)=\varepsilon_k(b)+\langle h_k,{\rm wt}(b)\rangle$,
\item[(3)] $\varepsilon_k(\tilde{e}_k(b))=\varepsilon_k(b)-1,\quad
\varphi_k(\tilde{e}_k(b))=\varphi_k(b)+1$\ if $\tilde{e}_k(b)\in B$, 
\item[(4)] $\varepsilon_k(\tilde{f}_k(b))=\varepsilon_k(b)+1,\ \ 
\varphi_k(\tilde{f}_k(b))=\varphi_k(b)-1$\ if $\tilde{f}_k(b)\in B$, 
\item[(5)] $\tilde{f}_k(b)=b'$ if and only if $b=\tilde{e}_k(b')$,
\item[(6)] if $\varphi_k(b)=-\infty$ then $\tilde{e}_k(b)=\tilde{f}_k(b)=0$.
\end{enumerate}
\end{defn}
It is well known that the sets $B(\infty)$, $B(\lambda)$ have crystal structures.

\begin{defn}
\begin{enumerate}
\item[$(1)$] For two crystals $B_1$, $B_2$,
a map $\psi : B_1\sqcup\{0\}\rightarrow B_2\sqcup\{0\}$ is said to be a {\it strict morphism} 
and denoted by $\psi : B_1\rightarrow B_2$
if $\psi(0)=0$ and the following conditions hold: For $k\in I$ and $b\in B_1$,
\begin{itemize}
\item if $\psi(b)\in B_2$ then for $i\in I$, it holds
\[
{\rm wt}(\psi(b))={\rm wt}(b),\quad
\varepsilon_i(\psi(b))=\varepsilon_i(b),\quad
\varphi_i(\psi(b))=\varphi_i(b),
\]
\item it holds $\tilde{e}_i(\psi(b))=\psi(\tilde{e}_i(b))$ and $\tilde{f}_i(\psi(b))=\psi(\tilde{f}_i(b))$ for $i\in I$, where
we understand $\tilde{e}_i(0)=\tilde{f}_i(0)=0$.
\end{itemize}
\item[$(2)$] If a strict morphism 
$\psi : B_1\sqcup\{0\}\rightarrow B_2\sqcup\{0\}$ is injective then
$\psi$ is said to be a {\it strict embedding} and denoted by $\psi:B_1 \hookrightarrow B_2$.
If $\psi$ is bijective then $\psi$ is said to be an {\it isomorphism}.
\end{enumerate}
\end{defn}

\subsection{An embedding}

Let
$\io=(\cd,i_r,\cd,i_2,i_1)$ be a sequence of indices from $I$ such that
\begin{equation}
{\hbox{
$i_r\ne i_{r+1}$ for $r\in\mathbb{Z}_{\geq1}$ and $\sharp\{r\in\mathbb{Z}_{\geq1}| i_r=k\}=\ify$ for all $k\in I$.}}
\label{seq-con}
\end{equation}
We can define a crystal structure on the set
\[
\ZZ^{\ify}
:=\{(\cd,a_r,\cd,a_2,a_1)| a_r\in\ZZ\text{ for }r\in\mathbb{Z}_{\geq1}
\text{ and it holds }a_r=0\,\,{\rm for}\ r\gg 0\},
\]
associated with $\iota$ (see subsection 2.4 of \cite{NZ}) and denote it by $\ZZ^{\ify}_{\iota}$.

\begin{thm}\cite{Kas93,NZ}\label{emb}
There is the unique strict embedding of crystals
\begin{equation}
\Psi_{\io}:B(\ify)\hookrightarrow \ZZ^{\ify}_{\io}
\label{psi}
\end{equation}
such that the highest weight vector $u_{\ify}\in B(\ify)$ is mapped to
$\textbf{0}:=(\cd,0,\cd,0,0)\in \ZZ^{\ify}_{\io}$: $\Psi_{\io} (u_{\ify}) = \textbf{0}$.
\end{thm}

\subsection{Nakashima-Zelevinsky's procedure}

Following \cite{NZ}, let us recall a procedure to compute an explicit form of ${\rm Im}(\Psi_{\io})$.
We set
\[
\QQ^{\ify}:=\{\textbf{a}=
(\cd,a_r,\cd,a_2,a_1)| a_r \in \QQ\text{ for }r\in\mathbb{Z}_{\geq1}
\text{ and it holds }
a_r = 0\,\,{\rm for\ } r \gg 0\}.
\]
For $r\in \mathbb{Z}_{\geq1}$, 
one defines $x_r\in(\QQ^{\ify})^*$ as
$x_r(\cd,a_r,\cd,a_2,a_1)=a_r$. We understand $x_r:=0$ when $r\in \mathbb{Z}_{<1}$.
For $r\in \mathbb{Z}_{\geq1}$, we set
\[
r^{(+)}:={\rm min}\{l\in\mathbb{Z}_{\geq1}\ |\ l>r\,\,{\rm and }\,\,i_r=i_{l}\},\ \ 
r^{(-)}:={\rm max} \{l\in\mathbb{Z}_{\geq1}\ |\ l<r\,\,{\rm and }\,\,i_r=i_{l}\}\cup\{0\},
\]
and
\begin{equation}
\beta_r:=
x_r+\sum_{r<j<r^{(+)}}\lan h_{i_r},\al_{i_j}\ran x_j+x_{r^{(+)}}\in (\QQ^{\ify})^*,\qquad
\beta_0:=0\in (\QQ^{\ify})^*.
\label{betak}
\end{equation}
We define a piecewise-linear operator 
$S_r=S_{r,\io}$ on $(\QQ^{\ify})^*$ by
\begin{equation}
S_r(\vp):=
\begin{cases}
\vp-\vp_r\beta_r & {\mbox{ if }}\vp_r>0,\\
 \vp-\vp_r\beta_{r^{(-)}} & {\mbox{ if }}\vp_r\leq 0.
\end{cases}
\label{Sk}
\end{equation}
We often simply write the map $S_r(\vp)$ as $S_r\vp$. 
Let us define
\begin{eqnarray}
\Xi_{\io} &:=  &\{S_{j_l}\cd S_{j_2}S_{j_1}x_{j_0}\,|\,
l\in\mathbb{Z}_{\geq0},j_0,j_1,\cd,j_l\geq1\},
\label{Xi_io}\\
\Sigma_{\io} & := &
\{\textbf{x}\in \ZZ^{\ify}\subset \QQ^{\ify}\,|\,\vp(\textbf{x})\geq0\,\,{\rm for}\,\,
{\rm any}\,\,\vp\in \Xi_{\io}\}.
\end{eqnarray}
We consider the following {\it positivity condition} on $\iota$:
\begin{equation}
\text{for any }
\vp=\sum_{k\in\mathbb{Z}_{\geq1}}\vp_kx_k\in \Xi_{\io},\ \text{if }\km=0\text{ then }\vp_k\geq0.
\label{posi}
\end{equation}

\begin{thm}\cite{NZ}\label{polyhthm}
For a sequence $\iota$ of indices satisfying $(\ref{seq-con})$ 
and $(\ref{posi})$, we have 
${\rm Im}(\Psi_{\io})=\Sigma_{\io}$.
\end{thm}

\subsection{Modified procedure}

Modifying the procedure in the previous subsection,
we get another procedure to compute ${\rm Im}(\Psi_{\iota})$.
Let
$\io=(\cd,i_r,\cd,i_2,i_1)$ be a sequence of indices from $I$ satisfying (\ref{seq-con}).
For each $r\in\mathbb{Z}_{\geq1}$,
a map
$S_r'=S_{r,\io}':(\QQ^{\ify})^*\rightarrow (\QQ^{\ify})^*$ is defined as follows:
For $\vp=\sum_{l\in\mathbb{Z}_{\geq1}}c_{l}x_{l}\in(\QQ^{\ify})^*$,
\begin{equation}
S_r'(\vp):=
\begin{cases}
\vp-\beta_r & {\mbox{ if }}\ c_r>0,\\
 \vp+\beta_{r^{(-)}} & {\mbox{ if }}\ c_r< 0,\\
 \vp &  {\mbox{ if }}\ c_r= 0.
\end{cases}
\label{SkSk}
\end{equation}
Here, $\beta_r$ is defined in (\ref{betak}). When $c_r>0$, we call $S'_r$ a {\it positive action}.
Just as in the previous subsection, one defines
\begin{eqnarray}
\Xi_{\io}' &:=  &\{S_{j_{l}}'\cd S_{j_2}'S_{j_1}'x_{j_0}\,|\,
l\in\mathbb{Z}_{\geq0},j_0,j_1,\cd,j_{l}\in\mathbb{Z}_{\geq1}\}, \label{xiprime}\\
\Sigma_{\io}' & := &
\{\textbf{a}\in \ZZ^{\ify}\subset \QQ^{\ify}\,|\,\vp(\textbf{a})\geq0\,\,{\rm for}\,\,
{\rm any}\,\,\vp\in \Xi_{\io}'\}.\nonumber
\end{eqnarray}
We impose the following $\Xi'$-{\it positivity condition} on $\io$:
\begin{equation}
\text{for any }
\vp=\sum_{l\in\mathbb{Z}_{\geq1}} c_{l}x_{l}\in \Xi_{\io}',\text{ if }l^{(-)}=0\text{ then }c_{l}\geq0.
\label{posiposi}
\end{equation}

\begin{thm}\label{inf-thm}\cite{Ka23a}
Let $\io=(\cdots,i_2,i_1)$ be a sequence satisfying $(\ref{seq-con})$ 
and $(\ref{posiposi})$. Then we get
${\rm Im}(\Psi_{\io})=\Sigma_{\io}'$.
\end{thm}

\subsection{Adapted sequences}

In this paper, we will consider only {\it adapted} sequences.

\begin{defn}\label{adapt}\cite{KaN20}
Let $A=(a_{i,j})_{i,j\in I}$ be the symmetrizable generalized Cartan matrix of $\mathfrak{g}$ and $\io=(\cdots,i_3,i_2,i_1)$
be a sequence that satisfies $(\ref{seq-con})$.
If the following condition holds then the sequence $\iota$ is said to be {\it adapted} to $A$:
For each $i,j\in I$ such that $a_{i,j}<0$, the subsequence of $\iota$ consisting of all $i$, $j$ is either
\[
(\cdots,i,j,i,j,i,j,i,j)\quad {\rm or}\quad (\cdots,j,i,j,i,j,i,j,i).
\]
When the matrix $A$ is fixed, we shortly say $\iota$ is {\it adapted}.
\end{defn}

For an adapted sequence $\iota=(\cdots,i_3,i_2,i_1)$,
we define a set of integers $(p_{i,j})_{i,j\in I;a_{i,j}<0}$ by
\begin{equation}\label{pij}
p_{i,j}=\begin{cases}
1 & {\rm if}\ {\rm the\ subsequence\ of\ }\iota{\rm\ consisting\ of}\ i,j\ {\rm is}\ (\cdots,j,i,j,i,j,i),\\
0 & {\rm if}\ {\rm the\ subsequence\ of\ }\iota{\rm\ consisting\ of}\ i,j\ {\rm is}\ (\cdots,i,j,i,j,i,j).
\end{cases}
\end{equation}
We can verify that if $a_{i,j}<0$ then
\begin{equation}\label{pij2}
p_{i,j}+p_{j,i}=1.
\end{equation}

\subsection{An identification of single indices with double indices}\label{s-d}

For a fixed sequence $\iota=(\cdots,i_2,i_1)$ satisfying (\ref{seq-con}),
we identify the set of single indices $\mathbb{Z}_{\geq1}$ with the set of double indices $\mathbb{Z}_{\geq1}\times I$ as follows:
We identify
each single index $r\in\mathbb{Z}_{\geq1}$ with a double index $(s,k)\in \mathbb{Z}_{\geq1}\times I$
when $i_r=k$ and $k$ appears $s$ times in $i_r$, $i_{r-1}$, $\cdots,i_1$.
For example, when $\iota=(\cdots,2,1,3,2,1,3,2,1,3)$, single indices
$\cdots,6,5,4,3,2,1$ are identified with double indices
\[
\cdots,(2,2),(2,1),(2,3),(1,2),(1,1),(1,3).
\]
The notation $x_{r}$, $\beta_{r}$, $S_r$ and $S'_{r}$ in Sect.2 are also written as
\[
x_r=x_{s,k},\quad \beta_r=\beta_{s,k},\quad S_r=S_{s,k}, \text{ and } S'_r=S'_{s,k}.
\]
When $(s,k)\notin \mathbb{Z}_{\geq1}\times I$, we understand $x_{s,k}:=0$.
By this identification and the ordinary order on $\mathbb{Z}_{\geq1}$ ($1<2<3<4<5<6<\cdots$),
we can naturally define an order on $\mathbb{Z}_{\geq1}\times I$. For $\iota=(\cdots,2,1,3,2,1,3,2,1,3)$, the order is
$\cdots>(2,2)>(2,1)>(2,3)>(1,2)>(1,1)>(1,3)$.
Using the notation in (\ref{pij}), if $\iota$ is adapted then $\beta_{s,k}$ is in the form
\begin{equation}\label{betask}
\beta_{s,k}=x_{s,k}+x_{s+1,k}+\sum_{j\in I; a_{k,j}<0} a_{k,j}x_{s+p_{j,k},j}.
\end{equation}
One can verify that
 $\{\beta_{t,i}\}_{t\in\mathbb{Z}_{\geq1},i\in I}$ is $\mathbb{Z}$-linearly independent.

\section{Monomial realizations}

Let us review the monomial realizations for crystal bases of highest weight representations.
We consider the set of Laurent monomials as follows:
\begin{gather}\label{ydef}
{\mathcal Y}:=\left\{X=\prod\limits_{s \in \mathbb{Z},\ i \in I}
X_{s,i}^{\zeta_{s,i}}\, \Bigg| \,\zeta_{s,i} \in \mathbb{Z},\ 
\zeta_{s,i}= 0\text{ except for finitely many }(s,i) \right\}.
\end{gather}
For $X=\prod\limits_{s \in \mathbb{Z},\; i \in I} X_{s,i}^{\zeta_{s,i}}\in {\mathcal Y}$, one sets
$\text{wt}(X):= \sum\limits_{s,i}\zeta_{s,i}\Lambda_i$ and

\vspace{-3mm}

\begin{gather*}
\varphi_i(X):=\max\left\{\! \sum\limits_{k\leq s}\zeta_{k,i}  \,|\, s\in \mathbb{Z} \!\right\},\!
\quad
\varepsilon_i(X):=\varphi_i(X)-\text{wt}(X)(h_i)=
\max\left\{\! -\sum\limits_{k> s}\zeta_{k,i}  \,|\, s\in \mathbb{Z} \!\right\}.
\end{gather*}
We fix an adapted sequence $\iota=(\cdots,i_3,i_2,i_1)$ and take $p_{i,j}$ as in (\ref{pij})
and put
\begin{equation}\label{ask}
A_{s,k}:=X_{s,k}X_{s+1,k}\prod\limits_{j\in I; a_{j,k}<0}X_{s+p_{j,k},j}^{a_{j,k}}\ \ \ (s\in\mathbb{Z},\ k\in I).
\end{equation}
For $i\in I$, let us define actions of Kashiwara operators as follows:
\begin{gather*}
\tilde{f}_iX:=
\begin{cases}
A_{n_{f_i},i}^{-1}X & \text{if} \quad \varphi_i(X)>0,
\\
0 & \text{if} \quad \varphi_i(X)=0,
\end{cases}
\quad
\tilde{e}_iX:=
\begin{cases}
A_{n_{e_i},i}X & \text{if} \quad \varepsilon_i(X)>0,
\\
0 & \text{if} \quad \varepsilon_i(X)=0,
\end{cases}
\end{gather*}
where we set
\begin{gather*}
n_{f_i}:=\min \left\{r\in\mathbb{Z} \,\Bigg|\, \varphi_i(X)= \sum\limits_{k\leq r}\zeta_{k,i}\right\},
\qquad
n_{e_i}:=\max \left\{r\in\mathbb{Z} \,\Bigg|\, \varphi_i(X)= \sum\limits_{k\leq r}\zeta_{k,i}\right\}.
\end{gather*}

\begin{thm}\label{mono-real}\cite{Kas03,Nak03}
\begin{enumerate}
\item[(i)] The set ${\mathcal Y}$ together with the above maps
${\rm wt}$, $\varepsilon_i$, $\varphi_i$ and $\tilde{e}_i$, $\tilde{f}_i$ $(i\in I)$ is a crystal. %%%P_{\rm cl}-crystal?
\item[(ii)] Taking $X \in {\mathcal Y}$ as $\tilde{\varepsilon}_i(X)=0$ for all $i \in I$, the set
\[
\{\tilde{f}_{j_m}\cdots\tilde{f}_{j_1}X | m\in\mathbb{Z}_{\geq0},\ j_1,\cdots,j_m \in I \}\setminus\{0\}
\]
is isomorphic to $B({\rm wt}(X))$. We denote 
$\{\tilde{f}_{j_m}\cdots\tilde{f}_{j_1}X_{s,k} | m\in\mathbb{Z}_{\geq0},\ j_1,\cdots,j_m \in I \}\setminus\{0\}$ by $\mathcal{M}_{s,k,\iota}$ for
$s\in\mathbb{Z}$ and $k\in I$, which is isomorphic to $B(\Lambda_k)$.
\end{enumerate}
\end{thm}

The following is straightforward from the definitions as above:
\begin{lem}\label{fundlem}
Let $X=\prod\limits_{s \in \mathbb{Z},\ i \in I}
X_{s,i}^{\zeta_{s,i}}\in {\mathcal Y}$ and $j\in I$.
\begin{enumerate}
\item[(1)] If $\tilde{f}_j X\neq 0$ and $\tilde{f}_j X=A_{s,j}^{-1}X$ with some $s\in\mathbb{Z}$
then $\zeta_{s,j}>0$.
\item[(2)] If $\tilde{e}_j X\neq 0$ and $\tilde{e}_j X=A_{s,j}X$ with some $s\in\mathbb{Z}$
then $\zeta_{s+1,j}<0$.
\item[(3)] If $\zeta_{s,j}<0$ with some $s\in\mathbb{Z}$ and $\zeta_{s',j}=0$ for all $s'>s$ then
$\varepsilon_j(X)>0$ so that
$\tilde{e}_j X\neq0$.
\end{enumerate}
\end{lem}

\begin{ex}\label{exc2-1}
Let $\mathfrak{g}$ be of type $C_2$ and $\iota=(\cdots,2,1,2,1,2,1)$. Then $\mathcal{M}_{s,1,\iota}$ is
\[
X_{s,1}\overset{1}{\rightarrow} \frac{X_{s,2}}{X_{s+1,1}}\overset{2}{\rightarrow} \frac{X_{s+1,1}}{X_{s+1,2}}\overset{1}{\rightarrow} \frac{1}{X_{s+2,1}}
\]
and $\mathcal{M}_{s,2,\iota}$ is
\[
X_{s,2}\overset{2}{\rightarrow}\frac{X_{s+1,1}^2}{X_{s+1,2}}\overset{1}{\rightarrow}\frac{X_{s+1,1}}{X_{s+2,1}}\overset{1}{\rightarrow}\frac{X_{s+1,2}}{X_{s+2,1}^2}\overset{2}{\rightarrow}\frac{1}{X_{s+2,2}}.
\]
\end{ex}

\begin{ex}\label{exA11-1}
Next, we set $\mathfrak{g}$ is of type $A^{(1)}_1$ and $\iota=(\cdots,2,1,2,1)$.
Then the partial crystal graph of $\mathcal{M}_{s,1,\iota}$ is as follows:
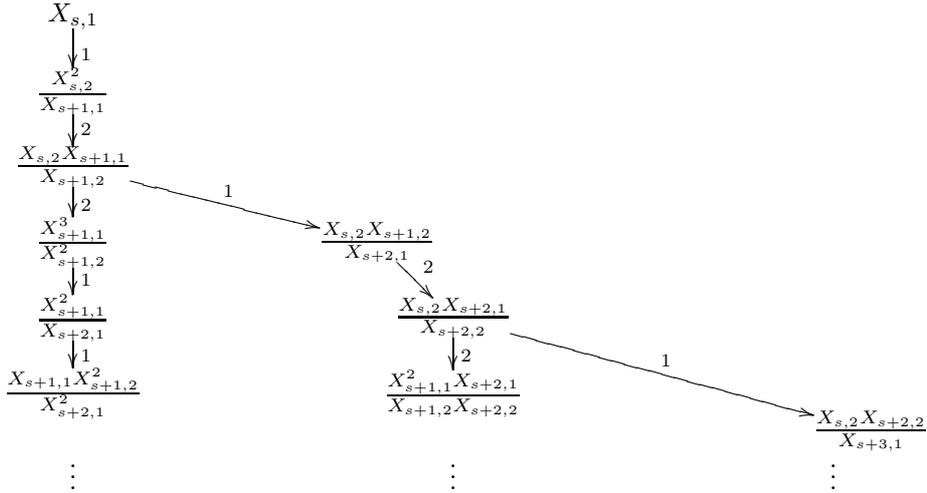
\[
\begin{xy}
(0,0) *{X_{s,1}}="1",
(0,-10) *{\frac{X_{s,2}^2}{X_{s+1,1}}}="2",
(0,-20) *{\frac{X_{s,2}X_{s+1,1}}{X_{s+1,2}}}="3",
(0,-30) *{\frac{X_{s+1,1}^3}{X_{s+1,2}^2}}="4-1",
(40,-30) *{\frac{X_{s,2}X_{s+1,2}}{X_{s+2,1}}}="4-2",
(0,-40) *{\frac{X_{s+1,1}^2}{X_{s+2,1}}}="5-1",
(50,-40) *{\frac{X_{s,2}X_{s+2,1}}{X_{s+2,2}}}="5-2",
(0,-50) *{\frac{X_{s+1,1}X_{s+1,2}^2}{X_{s+2,1}^2}}="6-1",
(50,-50) *{\frac{X_{s+1,1}^2X_{s+2,1}}{X_{s+1,2}X_{s+2,2}}}="6-2",
(105,-55) *{\frac{X_{s,2}X_{s+2,2}}{X_{s+3,1}}}="6-3",
(0,-60) *{\vdots}="dot1",
(50,-60) *{\vdots}="dot2",
(100,-60) *{\vdots}="dot3",
\ar@{->} "1";"2"^{1}
\ar@{->} "2";"3"^{2}
\ar@{->} "3";"4-1"^{2}
\ar@{->} "3";"4-2"^{1}
\ar@{->} "4-1";"5-1"^{1}
\ar@{->} "4-2";"5-2"^{2}
\ar@{->} "5-1";"6-1"^{1}
\ar@{->} "5-2";"6-2"^{2}
\ar@{->} "5-2";"6-3"^{1}
\end{xy}
\]
Similarly, 
the partial crystal graph of $\mathcal{M}_{s,2,\iota}$ is as follows:
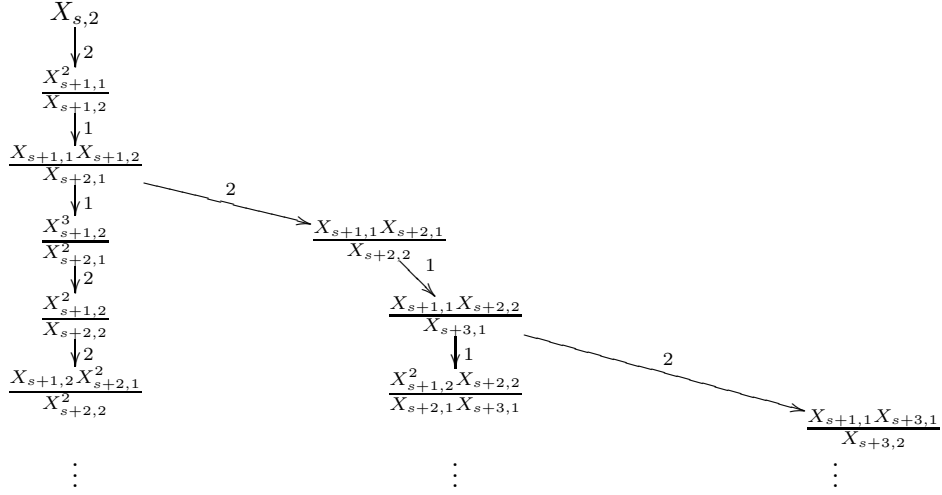
\[
\begin{xy}
(0,0) *{X_{s,2}}="1",
(0,-10) *{\frac{X_{s+1,1}^2}{X_{s+1,2}}}="2",
(0,-20) *{\frac{X_{s+1,1}X_{s+1,2}}{X_{s+2,1}}}="3",
(0,-30) *{\frac{X_{s+1,2}^3}{X_{s+2,1}^2}}="4-1",
(40,-30) *{\frac{X_{s+1,1}X_{s+2,1}}{X_{s+2,2}}}="4-2",
(0,-40) *{\frac{X_{s+1,2}^2}{X_{s+2,2}}}="5-1",
(50,-40) *{\frac{X_{s+1,1}X_{s+2,2}}{X_{s+3,1}}}="5-2",
(0,-50) *{\frac{X_{s+1,2}X_{s+2,1}^2}{X_{s+2,2}^2}}="6-1",
(50,-50) *{\frac{X_{s+1,2}^2X_{s+2,2}}{X_{s+2,1}X_{s+3,1}}}="6-2",
(105,-55) *{\frac{X_{s+1,1}X_{s+3,1}}{X_{s+3,2}}}="6-3",
(0,-60) *{\vdots}="dot1",
(50,-60) *{\vdots}="dot2",
(100,-60) *{\vdots}="dot3",
\ar@{->} "1";"2"^{2}
\ar@{->} "2";"3"^{1}
\ar@{->} "3";"4-1"^{1}
\ar@{->} "3";"4-2"^{2}
\ar@{->} "4-1";"5-1"^{2}
\ar@{->} "4-2";"5-2"^{1}
\ar@{->} "5-1";"6-1"^{2}
\ar@{->} "5-2";"6-2"^{1}
\ar@{->} "5-2";"6-3"^{2}
\end{xy}
\]

\end{ex}

\section{Main results}

We fix an adapted sequence $\iota=(\cdots,i_2,i_1)$.

\subsection{Tropicalizations}\label{tropsub}

We set
\begin{gather*}
{\mathcal Y}':=\left\{X=\prod\limits_{s \in \mathbb{Z}_{\geq1},\ k\in I}
X_{s,k}^{\zeta_{s,k}}\, \Bigg| \,\zeta_{s,k} \in \mathbb{Z},\ 
\zeta_{s,k}= 0\text{ except for finitely many }(s,k) \right\}
\end{gather*}
as a subset of ${\mathcal Y}$ in (\ref{ydef}) and define
\begin{equation}\label{H-def}
{\mathcal H}:=
\left\{
\sum_{s\in\mathbb{Z}_{\geq1},\ k\in I} \zeta_{s,k}x_{s,k} | \zeta_{s,k}\in\mathbb{Z},\ \zeta_{s,k}=0\text{ except for finitely many }(s,k)
\right\}.
\end{equation}
Here, $x_{s,k}\in(\mathbb{Q}^{\infty})$ is the notation as in subsection \ref{s-d}. We define a bijection
\[
Trop:{\mathcal Y}'\rightarrow {\mathcal H}
\]
as
\[
Trop\left(\prod\limits_{s \in \mathbb{Z}_{\geq1},\ k\in I}
X_{s,k}^{\zeta_{s,k}}\right)=\sum_{s\in\mathbb{Z}_{\geq1},\ k\in I} \zeta_{s,k}x_{s,k}.
\]
We define $DeTrop:{\mathcal H}\rightarrow {\mathcal Y}'$
as its inverse map:
\[
DeTrop\left(
\sum_{s\in\mathbb{Z}_{\geq1},\ k\in I} \zeta_{s,k}x_{s,k}
\right)
=\prod\limits_{s \in \mathbb{Z}_{\geq1},\ k\in I}
X_{s,k}^{\zeta_{s,k}}.
\]
Note that $Trop(A_{s,k})$ for $\mathfrak{g}$ coincides with
$\beta_{s,k}$ for $\mathfrak{g}^L$ by (\ref{betask}) and (\ref{ask}) when $s\geq1$.

\subsection{A conjecture and theorems}

\begin{conj}\label{mainconj}
Let $\iota$ be an adapted sequence and $\Psi_{\iota}$ be the map in Theorem \ref{emb} for $\mathfrak{g}$.
Let $\mathcal{M}_{s,k,\iota}$ be the set of monomials in Theorem \ref{mono-real} (ii) for $\mathfrak{g}^L$. Then
\[
{\rm Im}(\Psi_{\iota})=
\left\{
\mathbf{a}\in\mathbb{Z}^{\infty} |
\varphi(\mathbf{a})\geq0\ \text{for all }\varphi\in \bigcup_{s\in\mathbb{Z}_{\geq1},\ k\in I}Trop(\mathcal{M}_{s,k,\iota})
\right\}
\]
\end{conj}

Using maps in $(\ref{SkSk})$, we set
\[
\Xi'_{s,k,\iota}:=
\{
S_{j_1}'\cdots S_{j_m}'x_{s,k} |m\in\mathbb{Z}_{\geq0},\ j_1,\cdots,j_m\in\mathbb{Z}_{\geq1} 
\}
\]
for $s\in\mathbb{Z}_{\geq1}$, $k\in I$. Let $\Xi^{'+}_{s,k,\iota}$ be the set consisting of
$S_{j_1}'\cdots S_{j_m}'x_{s,k}$ such that $S'_{j}$ acts by positive actions ($j=j_1,\cdots,j_m$).
Clearly, $\Xi^{'+}_{s,k,\iota}\subset \Xi'_{s,k,\iota}$ holds.

We consider the condition
\begin{equation}\label{mainassump}
\Xi^{'+}_{s,k,\iota}= \Xi'_{s,k,\iota} \quad \text{for all }s\in\mathbb{Z}_{\geq1},k\in I.
\end{equation}
Note that $\Xi'$-positivity condition follows from the assumption (\ref{mainassump}).

\begin{thm}\label{mainthm1}
If the condition (\ref{mainassump}) holds then
Conjecture \ref{mainconj} is true.
\end{thm}

\begin{thm}\label{mainthm2}
\begin{enumerate}
\item[(1)]
When $\mathfrak{g}$ is a finite dimensional simple Lie algebra of type $A_n$, $B_n$, $C_n$ or $D_n$, the 
Conjecture \ref{mainconj} is true.
\item[(2)] When $\mathfrak{g}$ is a Kac-Moody algebra of rank $2$, the 
Conjecture \ref{mainconj} is true.
\item[(3)] When $\mathfrak{g}$ is a classical affine Lie algebra
of type
$A^{(1)}_{n-1},B^{(1)}_{n-1}$, $C^{(1)}_{n-1}$, $D^{(1)}_{n-1}$, $A^{(2)}_{2n-2}$, $A^{(2)}_{2n-3}$ or
$D^{(2)}_{n}$, the 
Conjecture \ref{mainconj} is true.
\end{enumerate}
\end{thm}

\begin{ex}\label{exc2-2}
We consider the setting as in Example \ref{exc2-1}. By the map Trop, we get the following homomorphisms from $\mathcal{M}_{s,k,\iota}$
($k=1,2$):
\[
x_{s,1},\ x_{s,2}-x_{s+1,1},\ x_{s+1,1}-x_{s+1,2},\ -x_{s+2,1},\quad
x_{s,2},\ 2x_{s+1,1}-x_{s+1,2},\ x_{s+1,1}-x_{s+2,1},\ x_{s+1,2}-2x_{s+2,1},\ -x_{s+2,2}.
\]
According to Theorem \ref{mainthm2} (1),
an explicit form of ${\rm Im}(\Psi_{\iota})$ of type $B_2(=(C_2)^L)$ is obtained from these homomorphisms.
By $a_{s,1}\geq0$, $-a_{s+2,1}\geq0$ for all $s\in\mathbb{Z}_{\geq1}$, we get $a_{m,1}=0$ when $m\geq3$
and by 
$a_{s,2}\geq0$, $-a_{s+2,2}\geq0$ for all $s\in\mathbb{Z}_{\geq1}$, we also get $a_{m,2}=0$ when $m\geq3$.
Simplifying the inequalities, one obtains
\[
{\rm Im}(\Psi_{\iota})
=
\{
\mathbf{a}=(a_{m,j})\in\mathbb{Z}^{\infty} | a_{1,2}\geq a_{2,1}\geq a_{2,2}\geq0,\ a_{1,1}\geq0,\ a_{m,1}=a_{m,2}=0\ (m\geq3)
\}.
\]
\end{ex}

\begin{ex}
We consider the same setting as in Example \ref{exA11-1}, that is,
$\mathfrak{g}$ is of type $A^{(1)}_1$ and $\iota=(\cdots,2,1,2,1)$.
Using Theorem \ref{mainthm2} (2), one obtains several inequalities defining 
${\rm Im}(\Psi_{\iota})$ of type $A^{(1)}_1(=(A^{(1)}_1)^L)$ from $\mathcal{M}_{s,k,\iota}$
($k=1,2$):
\[
{\rm Im}(\Psi_{\iota})
=
\left\{
\mathbf{a}=(a_{m,j})\in\mathbb{Z}^{\infty} | 
\begin{array}{l}
a_{s,1}\geq0,\ 2a_{s,2}-a_{s+1,1}\geq0,\ a_{s,2}+a_{s+1,1}-a_{s+1,2}\geq0,\ 3a_{s+1,1}-2a_{s+1,2}\geq0,\\
a_{s,2}+a_{s+1,2}-a_{s+2,1}\geq0,\ 2a_{s+1,1}-a_{s+2,1}\geq0,\ a_{s,2}+a_{s+2,1}-a_{s+2,2}\geq0,\\
a_{s+1,1}+2a_{s+1,2}-2a_{s+2,1}\geq0,\ 2a_{s+1,1}+a_{s+2,1}-a_{s+1,2}-a_{s+2,2}\geq0,\\
a_{s,2}+a_{s+2,2}-a_{s+3,1}\geq0,\cdots\\
a_{s,2}\geq0,\ 2a_{s+1,1}-a_{s+1,2}\geq0,\ a_{s+1,1}+a_{s+1,2}-a_{s+2,1}\geq0,\ 3a_{s+1,2}-2a_{s+2,1}\geq0,\\
a_{s+1,1}+a_{s+2,1}-a_{s+2,2}\geq0,\ 2a_{s+1,2}-a_{s+2,2}\geq0,\ a_{s+1,1}+a_{s+2,2}-a_{s+3,1}\geq0,\\
a_{s+1,2}+2a_{s+2,1}-2a_{s+2,2}\geq0,\ 2a_{s+1,2}+a_{s+2,2}-a_{s+2,1}-a_{s+3,1}\geq0,\\
a_{s+1,1}+a_{s+3,1}-a_{s+3,2}\geq0
\end{array}
\right\}.
\]

\end{ex}

\section{Proof of Theorem \ref{mainthm1}}

By the assumption (\ref{mainassump}), we get
$\Xi_{s,k,\iota}^{'+}= \Xi'_{s,k,\iota}$. For any $\varphi= S_{j_1}'\cdots S_{j_r}'x_{s,k}\in \Xi_{s,k,\iota}^{'+}$
with positive actions $S'_{j_1},\cdots,S'_{j_r}$, we can define
\begin{equation}\label{mdef}
|\varphi|=r.
\end{equation}
This is well-defined by the linear independence of $\{\beta_{t,i}\}_{t\in\mathbb{Z}_{\geq1},i\in I}$.
If we can prove $\Xi'_{s,k,\iota}=Trop(\mathcal{M}_{s,k,\iota})$ then our claim follows from Theorem \ref{inf-thm}.
The inclusion
$\Xi'_{s,k,\iota} \supset Trop(\mathcal{M}_{s,k,\iota})$ is clear from Lemma \ref{fundlem} (1),(2) so that
we show $\Xi'_{s,k,\iota} \subset Trop(\mathcal{M}_{s,k,\iota})$.
For any $\varphi\in \Xi'_{s,k,\iota}$,
let us show $\varphi\in Trop(\mathcal{M}_{s,k,\iota})$ via induction on $|\varphi|$.
When $|\varphi|=0$ then 
$\varphi=x_{s,k}$ 
so that the claim is clear so we assume $|\varphi|>0$.
Note that by $\varphi\in \Xi'_{s,k,\iota}=\Xi_{s,k,\iota}^{'+}$, one can write $\varphi$ as
\[
\varphi=x_{s,k}-\sum_{(t,i)\in\mathbb{Z}_{\geq s}\times I}c_{t,i}\beta_{t,i}
\]
with non-negative coefficients $c_{t,i}\in\mathbb{Z}_{\geq0}$ and $c_{t,i}=0$ except for finitely many $(t,i)$.
We take
\[
(r,j):={\rm max}\{(t,i)\in\mathbb{Z}_{\geq s}\times I | c_{t,i}>0\} 
\]
in the order defined in the subsection \ref{s-d}.
It follows from the definition of $\beta_{t,i}$ and adaptedness of $\iota$ that the coefficient of $x_{r+1,j}$ in $\varphi$
is negative. By the maximality of $(r,j)$ and adaptedness of $\iota$, we see that
the coefficients of $x_{l,j}$ ($l\geq r+2$) in
$\varphi$ are equal to $0$.
Thus, writing $\varphi=Trop(M)$ with $M=DeTrop(\varphi)$,
we get $\tilde{e}_jM\neq0$ from Lemma \ref{fundlem} (3)
and
\[
Trop(\tilde{e}_jM)=Trop(M)+\beta_{m,j}
\]
with some $m\in \mathbb{Z}_{\geq0}$. Note that since the coefficient of $x_{m+1,j}$ in $\varphi$ is negative 
by Lemma \ref{fundlem} (2)
so 
it holds $Trop(\tilde{e}_jM)=S'_{m+1,j}\varphi \in \Xi'_{s,k,\iota}$. It follows by (\ref{mainassump}) that
$|Trop(\tilde{e}_jM)|<|\varphi|$. % \beta no linear independence yori 
Thus, the induction assumption implies
$Trop(\tilde{e}_jM)\in Trop(\mathcal{M}_{s,k,\iota})$. By $\tilde{e}_jM\in \mathcal{M}_{s,k,\iota}$,
we have
$M=\tilde{f}_j\tilde{e}_jM\in\mathcal{M}_{s,k,\iota}$. Thus, we get $\varphi=Trop(M)\in Trop(\mathcal{M}_{s,k,\iota})$.\qed

\section{Finite case}

Let $\mathfrak{g}$ be of type $X=A_n$, $B_n$, $C_n$ or $D_n$ and $\iota$ be adapted.

\subsection{Boxes and tableaux}

For $k$ $(2\leq k\leq n)$, we put
\[
P^X(k):=\begin{cases}
p_{2,1}+p_{3,2}+\cdots+p_{n-2,n-3}+p_{n,n-2} & {\rm if}\ k=n\ {\rm and}\ X=D_n, \\
p_{2,1}+p_{3,2}+p_{4,3}+\cdots+p_{k,k-1} & {\rm if}\ {\rm otherwise}
\end{cases}
\]
and $P^X(0)=P^X(1)=P^X(n+1)=0$. 
We often write $P^X(k)=P(k)$ when $X$ is fixed.
By $p_{i,j}\in\{0,1\}$,
it holds for $k$, $l\in I$ such that $k\geq l$,
\begin{equation}\label{pineq}
P(k)\geq P(l),
\end{equation}
\begin{equation}\label{pineq2}
(k-l)+P(l)\geq P(k),
\end{equation}
except for the case $X=D_n$, $k=n$ and $l=n-1$. 

We define the following (partial) ordered sets $J_{A}$, $J_{B}$, $J_{C}$ and $J_{D}$:
\begin{itemize}
\item[(i)] $J_{A}:=\{1,2,\cdots,n,n+1\}$ with the order $1<2<\cdots<n<n+1$.
\item[(ii)]
$J_{B}=J_{C}:=\{1,2,\cdots,n,\overline{n},\cdots,\overline{2},\overline{1}\}$ with the order
\[
1<2<\cdots<n<\overline{n}<\cdots<\overline{2}<\overline{1}.\]
\item[(iii)]
$J_{D}:=\{1,2,\cdots,n,\overline{n},\cdots,\overline{2},\overline{1}\}$ with the partial order 
\[
 1< 2<\cdots< n-1<\ ^{n}_{\overline{n}}\ < \overline{n-1}< \cdots< \overline
 2< \overline 1.\]
\end{itemize}
For $j\in\{1,2,\cdots,n\}$, we set $|j|=|\overline{j}|=j$.

Next, we define the boxes and tableaux as homomorphisms. Note
that the notation is slightly different from \cite{KaN20}.
\begin{defn}\cite{KaN20}\label{box-def}
\begin{enumerate}
\item For $1\leq j\leq n+1$ and $s\in\mathbb{Z}$, we define
\[
\fbox{$j$}^{A}_{s}:=x_{s+P^A(j),j}-x_{s+P^A(j-1)+1,j-1}\in (\QQ^{\ify})^*.\]
\item
For $1\leq j\leq n$ and $s\in\mathbb{Z}$, we define
\[
\fbox{$j$}^{C}_{s}
:=x_{s+P^C(j),j}-x_{s+P^C(j-1)+1,j-1}\in (\QQ^{\ify})^*,\]
\[
\fbox{$\overline{j}$}^{C}_{s}
:=x_{s+P^C(j-1)+n-j+1,j-1}-x_{s+P^C(j)+n-j+1,j}\in (\QQ^{\ify})^*.
\]
\item
For $1\leq j\leq n-1$ and $s\in\mathbb{Z}$, we set
\[
\fbox{$j$}^{B}_{s}
:=x_{s+P^B(j),j}-x_{s+P^B(j-1)+1,j-1},
\quad
\fbox{$n$}^{B}_{s}:=2x_{s+P^B(n),n}-x_{s+P^B(n-1)+1,n-1}\in (\QQ^{\ify})^*,\]
\[
\fbox{$0$}^{B}_{s}:=x_{s+P^B(n),n}-x_{s+P^B(n)+1,n}\in (\QQ^{\ify})^*,
\]
\[
\fbox{$\overline{n}$}^{B}_{s}
:=x_{s+P^B(n-1)+1,n-1}-2x_{s+P^B(n)+1,n},
\quad
\fbox{$\overline{j}$}^{B}_{s}
:=x_{s+P^B(j-1)+n-j+1,j-1}-x_{s+P^B(j)+n-j+1,j}\in (\QQ^{\ify})^*,
\]
\[
\fbox{$\overline{n+1}$}^{B}_{s}
:=x_{s+P^B(n),n}
\in (\QQ^{\ify})^*.
\]
\item 
For $s\in\mathbb{Z}$, we set
\[
\fbox{$j$}^{D}_{s}
:=x_{s+P^D(j),j}-x_{s+P^D(j-1)+1,j-1}\in (\QQ^{\ify})^*,\ \ (1\leq j\leq n-2,\ j=n),
\]
\[
\fbox{$n-1$}^{D}_{s}
:=x_{s+P^D(n-1),n-1}+x_{s+P^D(n),n}-x_{s+P^D(n-2)+1,n-2}\in (\QQ^{\ify})^*,
\]
\[
 \fbox{$\overline{n}$}^{D}_{s}
:=x_{s+P^D(n-1),n-1}-x_{s+P^D(n)+1,n}\in (\QQ^{\ify})^*,
\]
\[
\fbox{$\overline{n-1}$}^{D}_{s}
:=x_{s+P^D(n-2)+1,n-2}-x_{s+P^D(n-1)+1,n-1}-x_{s+P^D(n)+1,n}\in (\QQ^{\ify})^*,
\]
\[
\fbox{$\overline{j}$}^{D}_{s}
:=x_{s+P^D(j-1)+n-j,j-1}-x_{s+P^D(j)+n-j,j}\in (\QQ^{\ify})^*,\ \ (1\leq j\leq n-2),\]
\[
\fbox{$\overline{n+1}$}^{D}_{s}
:=x_{s+P^D(n),n}\in (\QQ^{\ify})^*.
\]
\item
For $X=A$, $B$, $C$ or $D$,
we set
\[
\begin{ytableau}
j_1 \\
j_2 \\
\vdots \\
j_k
\end{ytableau}^{X}_{s}
:=
\fbox{$
j_k$}_s^{X}
+\fbox{$j_{k-1}$}_{s+1}^{X}+\cdots
+
\fbox{$
j_2 
$}^{X}_{s+k-2}
+
\fbox{$
j_1
$}^{X}_{s+k-1} \in (\QQ^{\ify})^*.
\]
\end{enumerate}
\end{defn}

%\begin{defn}
%Let $j\in I$ and $j_1,j_2\in J_X\cup\{\overline{n+1}\}$.
%For any $s\in\mathbb{Z}_{\geq1}$, if
%$\fbox{$
%j_1
%$}^{X}_s$ has a term $x_{m,j}$ ($m\in\mathbb{Z}_{\geq1}$) with a positive coefficient and
%$\fbox{$
%j_1
%$}^{X}_s$+
%$\fbox{$
%j_2
%$}^{X}_{s+1}$ does not have a term $x_{m,j}$ with a positive coefficient
%then we say the pair $(j_1,j_2)$ is a $j$-positive cancelling pair.
%If $\fbox{$
%j_2
%$}^{X}_{s+1}$ has a term $x_{m,j}$ ($m\in\mathbb{Z}_{\geq1}$) with a negative coefficient and
%$\fbox{$
%j_1
%$}^{X}_s$+
%$\fbox{$
%j_2
%$}^{X}_{s+1}$ does not have a term $x_{m,j}$ with a negative coefficient
%then we say the pair $(j_1,j_2)$ is a $j$-negative cancelling pair.
%\end{defn}

\begin{lem}\cite{KaN20}\label{box-lem}
\begin{enumerate}
\item When $\mathfrak{g}$ is of type $A_n$, it holds
\begin{equation}\label{A-box}
\fbox{$j+1$}^{A}_{s}=\fbox{$j$}^{A}_{s} -\beta_{s+P^A(j),j} \qquad (1\leq j\leq n,\ s\geq 1-P^A(j)).
\end{equation}
\item When $\mathfrak{g}$ is of type $B_n$, it holds
\begin{eqnarray}
\fbox{$j+1$}^{C}_{s}&=&\fbox{$j$}^{C}_{s} -\beta_{s+P^C(j),j} \qquad (1\leq j\leq n-1,\ s\geq 1-P^C(j)), \label{B-box1}\\
\fbox{$\overline{n}$}^{C}_{s}&=&\fbox{$n$}^{C}_{s} -\beta_{s+P^C(n),n} \qquad (s\geq 1-P^C(n)), \label{B-box2}\\
\fbox{$\overline{j-1}$}^{C}_{s}&=&\fbox{$\overline{j}$}^{C}_{s} -\beta_{s+P^C(j-1)+n-j+1,j-1}
 \qquad (2\leq j\leq n,\ s\geq j-P^C(j-1)-n). \label{B-box3}
\end{eqnarray}
\item When $\mathfrak{g}$ is of type $C_n$, it holds
\begin{eqnarray}
\fbox{$j+1$}^{B}_{s}&=&\fbox{$j$}^{B}_{s} -\beta_{s+P^B(j),j} \qquad (1\leq j\leq n-1,\ s\geq 1-P^B(j)), \label{C-box1}\\
\fbox{$0$}^{B}_{s}&=&\fbox{$n$}^{B}_{s} -\beta_{s+P^B(n),n} \qquad (s\geq 1-P^B(n)), \label{C-box2}\\
\fbox{$\overline{n}$}^{B}_{s}&=&\fbox{$0$}^{B}_{s} -\beta_{s+P^B(n),n} \qquad (s\geq 1-P^B(n)), \label{C-box22}\\
\fbox{$\overline{j-1}$}^{B}_{s}&=&\fbox{$\overline{j}$}^{B}_{s} -\beta_{s+P^B(j-1)+n-j+1,j-1} \ 
 (2\leq j\leq n,\ s\geq j-P^B(j-1)-n),\ \ \ \ \ \label{C-box3}\\
\fbox{$\ovl{n+1}$}^{B}_{l+1}+ \fbox{$\ovl{n}$}^{B}_{l}
&=&\fbox{$\ovl{n+1}$}^{B}_{l}-\beta_{l+P^B(n),n}
 \qquad (l\geq 1-P^B(n)). \label{BC-pr3}
\end{eqnarray}

\item When $\mathfrak{g}$ is of type $D_n$, it holds
\begin{eqnarray}
\fbox{$j+1$}^{D}_{s}&=&\fbox{$j$}^{D}_{s} -\beta_{s+P^D(j),j} \qquad (1\leq j\leq n-1,\ s\geq 1-P^D(j)), \label{D-box1}\\
\fbox{$\overline{n}$}^{D}_{s}&=&\fbox{$n-1$}^{D}_{s} -\beta_{s+P^D(n),n} \qquad (s\geq 1-P^D(n)), \label{D-box2}\\
\fbox{$\overline{n-1}$}^{D}_{s}&=&\fbox{$n$}^{D}_{s} -\beta_{s+P^D(n),n} \qquad (s\geq 1-P^D(n)), \label{D-box3}\\
\fbox{$\overline{j-1}$}^{D}_{s}&=&\fbox{$\overline{j}$}^{D}_{s} -\beta_{s+P^D(j-1)+n-j,j-1} \ (2\leq j\leq n,\ s\geq 1+j-P^D(j-1)-n),
\qquad \quad \label{D-box4} \\
\fbox{$\overline{n+1}$}^{D}_{l+2}
+\fbox{$\overline{n}$}^{D}_{l+1}
+\fbox{$\overline{n-1}$}^{D}_{l}
&=&\fbox{$\overline{n+1}$}^{D}_{l} -\beta_{l+P^D(n),n} \qquad (l\geq 1-P^D(n)). \label{D-box5}
\end{eqnarray}
\end{enumerate}
\end{lem}

\begin{defn}
\begin{enumerate}
\item[(1)]
For $X=A_n$ or $C_n$ and $k\in I$, we define
\[
{\rm Tab}^X_{k,\iota} := \{ \begin{ytableau}
j_1 \\
j_2 \\
\vdots \\
j_k
\end{ytableau}^{X}_{s} | j_i\in J_{X},\ s\geq 1-P^X(k),\ j_1<\cdots<j_k \}.
\]
\item[(2)] For $k\in I$ with $k<n$, we define
\[
{\rm Tab}^{B_n}_{k,\iota} := \{ \begin{ytableau}
j_1 \\
j_2 \\
\vdots \\
j_k
\end{ytableau}^{B}_{s} | j_i\in J_{B},\ s\geq 1-P^B(k),\ j_1\leq\cdots\leq j_k,\ \text{if }j_{i}=j_{i+1}\text{ then }j_i=0 \},
\]
\[
{\rm Tab}^{B_n}_{n,\iota}:=
\{ 
\begin{ytableau}
_{\overline{n+1}} \\
j_1 \\
\vdots \\
j_r
\end{ytableau}^{B}_{s} 
| 
\begin{array}{l}
r\in[0,n],\ j_1,\cdots,j_r\in J_{{B}},\\
\overline{n}\leq j_1<\cdots<j_r\leq\overline{1},\ s\geq 1-P^B(n).\\
\end{array}
\}
\]
\item[(3)] For $k\in I$ with $k<n-1$, we degfine
\[
{\rm Tab}^{D_n}_{k,\iota}:=
\{ 
\begin{ytableau}
j_1 \\
j_2 \\
\vdots \\
j_k
\end{ytableau}^{D}_{s} 
| 
\begin{array}{l}
j_1,\cdots,j_k\in J_{D}, \\
j_1\ngeq j_2\ngeq\cdots\ngeq j_k,\ s\geq 1-P^D(k).
\end{array}
\},
\]
\[
{\rm Tab}^{D_n}_{n-1,\iota}:=
\{
\begin{ytableau}
_{\overline{n+1}} \\
j_1 \\
\vdots \\
j_r
\end{ytableau}^{D}_{s} 
| 
\begin{array}{l}
r\in[0,n]\text{ is odd},\ j_1,\cdots,j_r\in J_{D},\\
\overline{n}\leq j_1<\cdots<j_r\leq\overline{1},\ s\geq 1-P^D(n-1).\\
\end{array}
\}
\]
\[
{\rm Tab}^{D_n}_{n,\iota}:=
\{
\begin{ytableau}
_{\overline{n+1}} \\
j_1 \\
\vdots \\
j_r
\end{ytableau}^{D}_{s} 
| 
\begin{array}{l}
r\in[0,n]\text{ is even},\ j_1,\cdots,j_r\in J_{D},\\
\overline{n}\leq j_1<\cdots<j_r\leq\overline{1},\ s\geq 1-P^D(n).\\
\end{array}
\}
\]

%{\rm if}\ j_1=\ovl{n+1}\ {\rm and}\ k\ {\rm is\ even} \ {\rm then}\ k\in[1,n+1],\ \overline{n}\leq j_2<\cdots<j_k\leq\overline{1},
%\ s\geq 1-P(n-1), \\
%{\rm if}\ j_1=\ovl{n+1}\ {\rm and}\ k\ {\rm is\ odd} \ {\rm then}\ k\in[1,n+1],\ \overline{n}\leq j_2<\cdots<j_k\leq\overline{1},
%\ s\geq 1-P(n).

\end{enumerate}
\end{defn}
For simplicity, we write each tableau 
$\begin{ytableau}
j_1 \\
\vdots \\
j_r
\end{ytableau}^{X}_{s} $
as $[j_1,\cdots,j_r]^{X}_s$.
For $s\geq 1$, we define
\[
{\rm Tab}^{X}_{s,k,\iota}:=
\{[j_1,\cdots,j_r]^{X}_{s-P^X(k)}\in {\rm Tab}^{X}_{k,\iota}\}.
\]

\subsection{Closedness of ${\rm Tab}^{X}_{s,k,\iota}$}

\begin{prop}\label{closedness-fin}
Let $\mathfrak{g}$ be of type $X=A_n,B_n,C_n$ or $D_n$. Then
${\rm Tab}^{X^L}_{s,k,\iota}$ is closed under the action of $S'_{m,j}$ for all $m\in \mathbb{Z}_{\geq1},j\in I$. 
\end{prop}

\nd
{\it Proof} For $s\in\mathbb{Z}_{\geq 1-P^{X^L}(k)}$, we need to show ${\rm Tab}^{X^L}_{s+P^{X^L}(k),k,\iota}$ is closed
under the action of $S'_{m,j}$.

\vspace{2mm}

\nd
\underline{Step1}

\vspace{2mm}

First, we show that for any
$T=[j_1,\cdots,j_r]^{X^L}_s\in {\rm Tab}^{X^L}_{s+P^{X^L}(k),k,\iota}$, if the coefficient of $x_{t,j}$ in $T$ is positive (resp. negative)
then $t\geq1$ (resp. $t\geq2$). In this proof, we simply write $P^{X^L}(j)$ as $P(j)$.

Recall that for $X=A, B$ or $C$, we have
$T=[j_1,\cdots,j_r]^{X^L}_s=\sum_{i=1}^{r}\fbox{$j_i$}^{X^L}_{s+r-i}$ and
by Definition \ref{box-def}, one obtains
\begin{equation}
\fbox{$j_i$}^{X^L}_{s+r-i}=
\begin{cases}\label{BC-box}
c(j_i)x_{s+r-i+P(j_i),j_i}-x_{s+r-i+P(j_i-1)+1,j_i-1} & {\rm if}\ j_i\in\{1,2,\cdots,n\},\\
x_{s+r-i+P(n),n}-x_{s+r-i+P(n)+1,n} & {\rm if}\ j_i=0,\\
x_{s+r-i+P(n),n} & {\rm if}\ j_i=\overline{n+1},\\
x_{s+r-i+P(|j_i|-1)+n-|j_i|+1,|j_i|-1}-c(j_i)x_{s+r-i+P(|j_i|)+n-|j_i|+1,|j_i|} & {\rm if}\ j_i\in \{\ovl{n},\cdots,\ovl{1}\},
\end{cases}
\end{equation}
where if $X^L=B_n$ and $j_i\in\{n,\ovl{n}\}$ then $c(j_i)=2$,
otherwise $c(j_i)=1$.
When $j_i\in\{1,2,\cdots,n\}$, it holds $r=k$ and
by $p_{l,l-1}=0$ or $1$ for $2\leq l\leq n$, $j_i\geq i$ and $s\geq 1-P(k)$,
we have
\begin{eqnarray}
s+r-i+P(j_i)&=&s+r-i+p_{2,1}+p_{3,2}+\cdots+p_{j_i,j_i-1}\nonumber \\
&\geq &s+r-i+p_{2,1}+p_{3,2}+\cdots+p_{i,i-1}\nonumber \\
&\geq &s+p_{2,1}+p_{3,2}+\cdots+p_{i,i-1} + p_{i+1,i}+\cdots+p_{r,r-1}\nonumber \\
&=&s+P(r)=s+P(k) \geq 1. \label{A-m1}
\end{eqnarray}
Similarly, if $j_i>i$ then we get
\begin{equation}\label{A-m2}
s+r-i+P(j_i-1)+1\geq2
\end{equation}
and if $j_i=i$ then we have $i=1$ or $j_{i-1}=i-1$ and $-x_{s+r-i+P(j_i-1)+1,j_i-1}$ in
$\fbox{$j_i$}^{X^L}_{s+r-i}$ is cancelled in $\fbox{$j_{i-1}$}^{X^L}_{s+r-i+1}+\fbox{$j_i$}^{X^L}_{s+r-i}$.
When $j_i=0$ or $j_i=\ovl{n+1}$, it is easy to see
\begin{equation}\label{A-m3}
s+r-i+P(n)\geq s+P(k)\geq1,\quad s+r-i+P(n)+1\geq2
\end{equation}
from $r\geq i$, $P(k)\leq P(n)$ ((\ref{pineq})) and $s\geq 1-P(k)$. When 
$j_i\in \{\ovl{n},\cdots,\ovl{1}\}$, it holds 
$P(|j_i|-1)+n-|j_i|+1\geq P(n)$, $P(|j_i|)+n-|j_i|\geq P(n)$ so that
\begin{equation}\label{A-m4}
s+r-i+P(|j_i|-1)+n-|j_i|+1\geq s+P(k)\geq1,\quad s+r-i+P(|j_i|)+n-|j_i|+1\geq2.
\end{equation}
For $X=D$,
we see that 
$T=[j_1,\cdots,j_r]^{D}_s=\sum_{i=1}^{r}\fbox{$j_i$}^{D}_{s+r-i}$, and
by Definition \ref{box-def}, we obtain
\begin{equation}\label{D-box}
\fbox{$j_i$}^{D}_{s+r-i}=
\begin{cases}
x_{s+r-i+P(j_i),j_i}-x_{s+r-i+P(j_i-1)+1,j_i-1} & {\rm if}\ j_i\in[1,n-2]\cup\{n\},\\
x_{s+r-i+P(n-1),n-1}+x_{s+r-i+P(n),n}-x_{s+r-i+P(n-2)+1,n-2} & {\rm if}\ j_i=n-1,\\
x_{s+r-i+P(n),n} & {\rm if}\ j_i= \ovl{n+1}, \\ 
x_{s+r-i+P(n-1),n-1}-x_{s+r-i+P(n)+1,n} & {\rm if}\ j_i= \ovl{n}, \\ 
x_{s+r-i+P(n-2)+1,n-2}-x_{s+r-i+P(n-1)+1,n-1}
-x_{s+r-i+P(n)+1,n} & {\rm if}\ j_i= \ovl{n-1},\\
x_{s+r-i+P(|j_i|-1)+n-|j_i|,|j_i|-1}-x_{s+r-i+P(|j_i|)+n-|j_i|,|j_i|} & {\rm if}\ j_i\geq \ovl{n-2}.
\end{cases}
\end{equation}
When $j_i\in\{1,2,\cdots,n\}$, we have $r=k\leq n-2$ so that $P(k)\leq P(n), P(n-1)$ and just as in (\ref{A-m1}), (\ref{A-m2}), one can show
\begin{equation}\label{D-m1}
s+r-i+P(j_i)\geq s+P(k)\geq1,\quad \text{if }j_i>i\text{ then }
s+r-i+P(j_i-1)+1\geq2.
\end{equation}
When $j_i=\ovl{n+1}$, if $k=n$ then $s+r-i+P(n)\geq1$ is clear by $s\geq 1-P(k)=1-P(n)$. If $k=n-1$ then we have $i=1<r$ so that
\begin{equation}\label{D-m2}
s+r-i+P(n)\geq s+P(k)\geq1.
\end{equation}
When $j_i=\ovl{n}$, if $k\leq n-1$ then $s+r-i+P(n-1)\geq1$ is clear. If $k=n$ then $r\geq3$ and $r-i\geq1$ so that
\begin{equation}\label{D-m3}
s+r-i+P(n-1)\geq s+P(k)\geq1.
\end{equation}
Similarly, the relation
\begin{equation}\label{D-m4}
s+r-i+P(n)+1\geq2
\end{equation}
is clear except for the case
$r=i$ and $k=n-1$.
In this case, we have $T=[\ovl{n+1},\ovl{n}]_s^D=x_{s+P(n-1),n-1}$ and the term $-x_{s+r-i+P(n)+1,n}$ is cancelled.
Similarly, when $j_i=\ovl{n-1}$, we have
\begin{equation}\label{D-m5}
s+r-i+P(n-2)+1\geq s+P(k)\geq1
\end{equation}
and it holds
\begin{equation}\label{D-m6}
s+r-i+P(n-1)+1\geq2,\quad
s+r-i+P(n)+1\geq2
\end{equation}
or terms $-x_{s+r-i+P(n-1)+1,n-1}$ or $-x_{s+r-i+P(n)+1,n}$ are cancelled.

When $j_i\geq \ovl{n-2}$, by $P(|j_i|-1)+n-|j_i|\geq P(n),P(n-1)$ and
$P(|j_i|)+n-|j_i|-1\geq P(n),P(n-1)$, one obtains
\begin{equation}\label{D-m7}
s+r-i+P(|j_i|-1)+n-|j_i|\geq s+P(k)\geq1,\quad 
s+r-i+P(|j_i|)+n-|j_i|\geq 2.
\end{equation}
Therefore, considering (\ref{A-m1})-(\ref{A-m4}) and (\ref{D-m1})-(\ref{D-m7}), if
the coefficient of $x_{t,j}$ in $T$ is positive (resp. negative)
then $t\geq1$ (resp. $t\geq2$).

\vspace{2mm}

\nd
\underline{Step2}

\vspace{2mm}

Next, for any $T=[j_1,\cdots,j_r]^{X^L}_s\in {\rm Tab}^{X^L}_{s+P(k),k,\iota}$, let us show $S'_{m,j}T\in {\rm Tab}^{X^L}_{s+P(k),k,\iota}$.
When the coefficient of
$x_{m,j}$ in $T$ equals $0$, the claim is clear so that we assume the coefficient is non-zero.
First, we consider the case the coefficient of
$x_{m,j}$ in $T$ is positive.

\vspace{2mm}

\nd
\underline{Case 1 : $j<n$}

\vspace{2mm}

There is $i\in[1,r]$ such that
$j_i=j$ or $j_i=\overline{j+1}$
and
$\fbox{$j_i$}^{X^L}_{s+k-i}$ has a term $x_{m,j}$ with a positive coefficient
and it holds either $i=r$ or $(j_i,j_{i+1})\neq (j,j+1), (\overline{j+1},\overline{j}), (\overline{n},n)$.
When $j_i=j$ (resp. $j_{i}=\overline{j+1}$), we see that
\[
S'_{m,j}T=T-\beta_{m,j}=[j_1,\cdots,j_{i-1},j+1,j_{i+1},\cdots,j_r]_s^{X^L}\ (\text{resp. }[j_1,\cdots,j_{i-1},\overline{j},j_{i+1},\cdots,j_r]^{X^L}_s) \in{\rm Tab}^{X^L}_{s+P(k),k,\iota}
\]
by Lemma \ref{box-lem}.

\vspace{2mm}

\nd
\underline{Case 2 : $j=n$}

\vspace{2mm}

We see that there is $i\in[1,r]$ such that 
$\fbox{$j_i$}^{X^L}_{s+k-i}$ has a term $x_{m,n}$ with a positive coefficient and
the triple ($X^L$, $j_i$, $j_{i+1}$) satisfies one of the following: 
\begin{table}[H]
  \begin{tabular}{|c|c|c|c|c|c|c|c|c|} \hline
  $X^L$ & $A_{n}$ & $C_{n}$ & $B_{n}$ & $B_n$ & $B_n$ & $D_n$ & $D_{n}$ & $D_n$ \\ \hline
  $j_i$ & $=n$ & $=n$ & $=n$ & $=0$ & $=\overline{n+1}$ & $=n-1$ & $=n$ & $=\overline{n+1}$  \\ \hline
  $j_{i+1}$ & $\neq n+1$ & $\neq \ovl{n}$ & $\neq \ovl{n}$ & $\neq0,\ovl{n}$ & $\neq \ovl{n}$ & $\neq \ovl{n},\ovl{n-1}$ &
$\neq \ovl{n},\ovl{n-1}$ & $\neq \ovl{n},\ovl{n-1}$ \\ \hline
  \end{tabular}
\end{table}
Here, when $i=r$, we understand $j_{i+1}$ satisfies the above condition.
When $(X^L,j_i)=(A_n,n)$ (resp. $(C_n,n)$, $(B_n,n)$, $(B_n,0)$, $(B_n,\ovl{n+1})$, $(D_n,n-1)$, $(D_n,n)$, $(D_n,\ovl{n+1})$),
let $T'$ be a tableau $[j_1,\cdots,j_{i-1},n+1]$ (resp. $[j_1,\cdots,j_{i-1},\overline{n},j_{i+1},\cdots,j_k]$, 
$[j_1,\cdots,j_{i-1},0,j_{i+1},\cdots,j_k]$, $[j_1,\cdots,j_{i-1},\overline{n},j_{i+1},\cdots,j_k]$, $[\overline{n+1},\overline{n},j_{2},\cdots,j_k]$, 
$[j_1,\cdots,j_{i-1},\overline{n},j_{i+1},\cdots,j_k]$, $[j_1,\cdots,j_{i-1},\overline{n-1},j_{i+1},\cdots,j_k]$, 
$[\overline{n+1},\overline{n},\overline{n-1},j_{2},\cdots,j_k]$) in ${\rm Tab}^{X^L}_{s+P(k),k,\iota}$.
Here, we simply write $[j_1,\cdots,j_r]_s^{X^L}$ as $[j_1,\cdots,j_r]$.
One obtains
\[
S'_{m,n}T=T-\beta_{m,n}=T'
\]
by Lemma \ref{box-lem}.
 
We can similarly show $S'_{m,j}T\in {\rm Tab}^{X^L}_{s+P(k),k,\iota}$ when the coefficient of $x_{m,j}$
is negative.
Thus, we get our claim. \qed

\subsection{Proof of Theorem \ref{mainthm2} (1)}

Let $\mathfrak{g}$ of type $X$ and
we show ${\rm Tab}^{X^L}_{s+P(k),k,\iota}=\Xi_{s+P(k),k,\iota}'=\Xi_{s+P(k),k,\iota}'^+$ for $s\in\mathbb{Z}_{\geq1-P(k)}$,
which yields our claim by Theorem \ref{mainthm1}.
Let $P_{X^L}=\bigoplus_{j\in I}\mathbb{Z}\Lambda_j$ be the weight lattice of type $X^L$
with the following partial order : For $\lambda$, $\mu\in P_{X^L}$,
$\lambda \geq \mu $ if and only if $\lambda - \mu \in \bigoplus_{j\in I}\mathbb{Z}_{\geq0}\alpha_j$.

For the set ${\mathcal H}$ in (\ref{H-def}),
one considers
the $\mathbb{Z}$-linear map
${\rm wt} : {\mathcal H}\rightarrow P_{X^L}$
defined as ${\rm wt}(x_{r,j}):=\Lm_j$ for any $r\in\mathbb{Z}_{\geq1}$ and $j\in I$.
The explicit form (\ref{betask}) means 
\begin{equation}\label{be-al}
{\rm wt}(\beta_{r,j})=\al_{j}\in P_{X^L}.
\end{equation} 
First,
putting $s':=s+P(k)$,
we prove ${\rm Tab}^{X^L}_{s',k,\iota}\supset \Xi_{s',k,\iota}'\supset\Xi_{s',k,\iota}'^+$.
Note that
\[
x_{s',k}=
\begin{cases}
[\overline{n+1}]^{X^L}_s & \text{ if }k=n,\ X^L=B_n\text{ or }D_n,\\
[\overline{n+1},\overline{n}]^{X^L}_s & \text{ if }k=n-1,\ X^L=D_n,\\
[1,2,\cdots,k]^{X^L}_s & \text{otherwise.}
\end{cases}
\]
Hence, $x_{s',k}\in{\rm Tab}^{X^L}_{s',k,\iota}$. By the definition of 
$\Xi_{s',k,\iota}'$ and Proposition \ref{closedness-fin}, we get ${\rm Tab}^{X^L}_{s',k,\iota}\supset \Xi_{s',k,\iota}'$.

Next, let us show ${\rm Tab}^{X^L}_{s',k,\iota}\subset \Xi_{s',k,\iota}'^+$.
For any $T=[j_1,\cdots,j_r]_s\in {\rm Tab}^{X^L}_{s',k,\iota}$, let us show
$T\in \Xi_{s',k,\iota}'^+$ by the induction on the weight ${\rm wt}(T)$ of $T$. Note that
it follows from Lemma \ref{box-lem} and (\ref{be-al}) that ${\rm wt}(T)\leq \Lambda_k={\rm wt}(x_{s',k})$.
When
\[
T=
\begin{cases}
[\overline{n+1}]^{X^L}_s & \text{ if }k=n,\ X^L=B_n\text{ or }D_n,\\
[\overline{n+1},\overline{n}]^{X^L}_s & \text{ if }k=n-1,\ X^L=D_n,\\
[1,2,\cdots,k]^{X^L}_s & \text{otherwise},
\end{cases}
\]
it holds $T=x_{s',k}\in\Xi_{s',k,\iota}'^+$. Hence, we may assume that $T\neq x_{s',k}$.
Using Lemma \ref{box-lem}, (\ref{A-m1}), (\ref{A-m3}), (\ref{A-m4}),
(\ref{D-m1}), (\ref{D-m2}), (\ref{D-m3}), (\ref{D-m5}) and (\ref{D-m7}),  
we can write
\[
T=x_{s',k}-\sum_{(t,j)\in\mathbb{Z}_{\geq s'}\times I}c_{t,j}\beta_{t,j}
\]
with non-negative coefficients $c_{t,j}\in\mathbb{Z}_{\geq0}$.
Except for finitely many $(t,j)$, it holds $c_{t,j}=0$.
We take
\[
(t',j'):={\rm max}\{(t,j)\in\mathbb{Z}_{\geq s'}\times I | c_{t,j}>0\}
\]
in the order defined in the subsection \ref{s-d}.
It follows from the definition of $\beta_{t,j}$ and adaptedness of $\iota$ that the coefficient of $x_{t'+1,j'}$ in $T$
is negative. We can take
\[
(t'',j''):={\rm min}\{(t,j)\in\mathbb{Z}_{\geq s'}\times I | \text{the coefficient of }x_{t,j}\text{ in }T\text{ is negative}\}.
\]
\underline{Case 1 : $j''<n$}

\vspace{2mm}

There exists $i\in[1,r]$ such that 
$\fbox{$j_i$}^{X^L}_{s+r-i}$ has a term $x_{t'',j''}$ with a negative coefficient
and it holds either $j_i=j''+1$ and $j_{i-1}\neq j'',\ovl{n}$ or $j_i=\overline{j''}$ and $j_{i-1}\neq \overline{j''+1}$.
Let $T'=[j_1,\cdots,j_{i-1},j'',j_{i+1},\cdots,j_r]^{X^L}_s\in {\rm Tab}^{X^L}_{s',k,\iota}$ or
$T'=[j_1,\cdots,j_{i-1},\overline{j''+1},j_{i+1},\cdots,j_r]^{X^L}_s\in {\rm Tab}^{X^L}_{s',k,\iota}$.
The minimality of $(t'',j'')$ implies the coefficient of $x_{t''-1,j''}$ in $T'$ is positive.
One obtains
\[
T=T'-\beta_{t''-1,j''}=S'_{t''-1,j''}T'
\]
by Lemma \ref{box-lem}. Note that by Step1 in the proof of Proposition \ref{closedness-fin}, we get $t''>1$.
By the induction assumption, it holds $T'\in \Xi_{s',k,\iota}'^+$, which yields our claim
$T\in \Xi_{s',k,\iota}'^+$.

\vspace{2mm}

\nd
\underline{Case 2 : $j''=n$, $(k,X^L)\neq (n,B_n),(n-1,D_n),(n,D_n)$}

\vspace{2mm}

There exists $i\in[1,r]$ such that
$\fbox{$j_i$}^{X^L}_{s+r-i}$ has a term $x_{t'',n}$ with a negative coefficient and
the triple ($X^L$,$j_i$,$j_{i-1}$) is one of the following:
\begin{table}[H]
  \begin{tabular}{|c|c|c|c|c|c|c|} \hline
  $X^L$ & $A_{n}$ & $C_{n}$ & $B_{n}$ & $B_n$ & $D_n$ & $D_{n}$ \\ \hline
  $j_{i-1}$ & $\neq n$ & $\neq n$ & $\neq n,0$ & $\neq n$ & $\neq n-1,n$ & $\neq n-1,n$ \\ \hline
  $j_{i}$ & $=n+1$ & $=\ovl{n}$ & $=0$ & $=\ovl{n}$ & $=\ovl{n}$ & $=\ovl{n-1}$ \\ \hline
  \end{tabular}
\end{table}
Here, if $i=1$ then we understand $j_{i-1}$ satisfies the above condition.
When $(j_i,X^L)=(n+1,A_n)$ (resp. $(j_i,X^L)=(\overline{n},C_n),(0,B_n),(\overline{n},B_n),(\overline{n},D_n),(\overline{n-1},D_n)$),
putting $T'=[j_1,\cdots,j_{i-1},n]$ (resp.
$[j_1,\cdots,j_{i-1},n,j_{i+1},\cdots,j_r]$, $[j_1,\cdots,j_{i-1},n,j_{i+1},\cdots,j_r]$, $[j_1,\cdots,j_{i-1},0,j_{i+1},\cdots,j_r]$,
$[j_1,\cdots,j_{i-1},n-1,j_{i+1},\cdots,j_r]$, $[j_1,\cdots,j_{i-1},n,j_{i+1},\cdots,j_r]$),
one obtains $T'\in {\rm Tab}^{X^L}_{s',k,\iota}$. Here, we simply write $[j_1,\cdots,j_r]^{X^L}_s$ as $[j_1,\cdots,j_r]$.
It holds 
\[
T=T'-\beta_{t''-1,n}=S'_{t''-1,n}T'
\]
by Lemma \ref{box-lem}. Step1 in the proof of Proposition \ref{closedness-fin} yields $t''>1$.
The induction assumption yields $T'\in \Xi_{s',k,\iota}'^+$ so that
$T\in \Xi_{s',k,\iota}'^+$. 

\vspace{2mm}

\nd
\underline{Case 3 : $j''=n$, $(k,X^L)= (n,B_n),(n-1,D_n)\text{ or }(n,D_n)$}

\vspace{2mm}

It holds either $X^L=B_n$ and $j_1=\overline{n+1}$, $j_2=\overline{n}$ or $X^L=D_n$
and $j_1=\overline{n+1}$, $j_2=\overline{n}$, $j_3=\overline{n-1}$.
In the former case, we set $T':=[\overline{n+1},j_3,\cdots,j_r]_s$
and in the latter case, we set $T':=[\overline{n+1},j_4,\cdots,j_r]_s$.
Then
we have $T'\in {\rm Tab}^{X^L}_{s',k,\iota}$ and
\[
T=T'-\beta_{t''-1,n}=S'_{t''-1,n}T'
\]
by Lemma \ref{box-lem}. Using Step1 in the proof of Proposition \ref{closedness-fin}, we have $t''>1$.
From the induction assumption, it holds $T'\in \Xi_{s',k,\iota}'^+$ so that
$T\in \Xi_{s',k,\iota}'^+$.

\qed

\vspace{2mm}

\section{Rank $2$ case}

We take $\mathfrak{g}$ as a rank $2$ Kac-Moody algebra with a generalized Cartan matrix
$A
=\begin{pmatrix}
2 & -a \\
-b & 2
\end{pmatrix}$ ($a,b\in\mathbb{Z}_{\geq0}$). 
Let 
\[
\iota=(\cdots,2,1,2,1,2,1).
\]
We fix $s\in\mathbb{Z}$, $i\in I$ and define $i'\in I$ as $\{i,i'\}=\{1,2\}=I$.
One defines $N:=|\lan s_2s_1\ran|$ for $s_2s_1\in W$.
Note that $N=\infty$ if $\mathfrak{g}$ is not finite dimensional.
For $m\in\mathbb{Z}_{\geq0}$
such that $2m<N$,
let $P_{m,i}^{(i)}$, $P_{m,i'}^{(i)}$be integers defined by
\[
(s_{i'}s_i)^m\Lambda_i=P^{(i)}_{m,i}\Lambda_i-P^{(i)}_{m,i'}\Lambda_{i'}.
\]
Note that when $2m+1<N$, we have
\[
s_i(s_{i'}s_i)^m\Lambda_i=P^{(i)}_{m+1,i'}\Lambda_{i'}-P^{(i)}_{m,i}\Lambda_{i},
\]
where if $2m+2\geq N$ then we set $P^{(i)}_{m+1,i'}=0$.

Then $\Xi_{s,i,\iota}:=\{S_{j_1}\cdots S_{j_m}x_{s,i}|m\in\mathbb{Z}_{\geq0}\}$ for $\mathfrak{g}^L$ is given by
\[
\Xi_{s,1,\iota}
=\{P^{(1)}_{m,1}x_{s+m,1}-P^{(1)}_{m,2}x_{s+m,2}| m\in\mathbb{Z}_{\geq0}, 2m<N\}\cup
\{P^{(1)}_{m+1,2}x_{s+m,2}-P^{(1)}_{m,1}x_{s+m+1,1}| m\in\mathbb{Z}_{\geq0}, 2m+1<N \}
\]
and
\[
\Xi_{s,2,\iota}
=\{P^{(2)}_{m,2}x_{s+m,2}-P^{(2)}_{m,1}x_{s+m+1,1}|m\in\mathbb{Z}_{\geq0}, 2m<N\}
\cup
\{P^{(2)}_{m+1,1}x_{s+m+1,1}-P^{(2)}_{m,2}x_{s+m+1,2}| m\in\mathbb{Z}_{\geq0}, 2m+1<N \}
\]
by the definition of $S_r$ in (\ref{Sk}). 

\begin{thm}\label{prev}\cite{Ka24b}
Let $\mathcal{M}_{s,i,\iota}$ be the monomial realization for $\mathfrak{g}$ in Theorem \ref{mono-real} (ii).
\begin{enumerate}
\item[(1)]
The set $\mathcal{M}_{s,1,\iota}$ includes Laurent monomials
in the form
\[
\frac{X^{P^{(1)}_{m,1}}_{s+m,1}}{X^{P^{(1)}_{m,2}}_{s+m,2}}\ (1\leq m<\frac{1}{2}N),\quad
\frac{X^{P^{(1)}_{m+1,2}}_{s+m,2}}{X^{P^{(1)}_{m,1}}_{s+m+1,1}}\ (1\leq m<\frac{1}{2}(N-1)).
\]
Ohter Laurent monomials in $\mathcal{M}_{s,1,\iota}$
are expressed by a product of above monomials with exponents of positive rational numbers.
\item[(2)]\ The set $\mathcal{M}_{s,2,\iota}$ includes Laurent monomials
in the form
\[
\frac{X^{P^{(2)}_{m,2}}_{s+m,2}}{X^{P^{(2)}_{m,1}}_{s+m+1,1}}\ (1\leq m<\frac{1}{2}N),\quad
\frac{X^{P^{(2)}_{m+1,1}}_{s+m+1,1}}{X^{P^{(2)}_{m,2}}_{s+m+1,2}}\ (1\leq m<\frac{1}{2}(N-1)).
\]
Ohter Laurent monomials in $\mathcal{M}_{s,2,\iota}$
are expressed by a product of above monomials with exponents of positive rational numbers.
\end{enumerate}
\end{thm}

\nd
[{\it Proof of Theorem \ref{mainthm2} (2)}]

For $s\in\mathbb{Z}_{\geq1}$,
let $\mathcal{M}'_{s,1,\iota}$ and $\mathcal{M}'_{s,2,\iota}$ be the
sets of monomials in Theorem \ref{prev} (1) and (2), respectively:
\[
\mathcal{M}'_{s,1,\iota}:=
\left\{
\frac{X^{P^{(1)}_{m,1}}_{s+m,1}}{X^{P^{(1)}_{m,2}}_{s+m,2}}\left| 1\leq m<\frac{1}{2}N \right.\right\}
\cup
\left\{
\frac{X^{P^{(1)}_{m+1,2}}_{s+m,2}}{X^{P^{(1)}_{m,1}}_{s+m+1,1}}\left| 1\leq m<\frac{1}{2}(N-1) \right.
\right\},
\]
\[
\mathcal{M}'_{s,2,\iota}:=
\left\{
\frac{X^{P^{(2)}_{m,2}}_{s+m,2}}{X^{P^{(2)}_{m,1}}_{s+m+1,1}}\left| 1\leq m<\frac{1}{2}N\right.
\right\}
\cup
\left\{
\frac{X^{P^{(2)}_{m+1,1}}_{s+m+1,1}}{X^{P^{(2)}_{m,2}}_{s+m+1,2}}\left| 1\leq m<\frac{1}{2}(N-1)\right.
\right\}.
\]
Then it is clear that the map $Trop$ induces a bijection between 
$\mathcal{M}'_{s,i,\iota}$ and $\Xi_{s,i,\iota}$ for $i=1,2$.
For $M\in \mathcal{M}_{s,i,\iota}\setminus \mathcal{M}'_{s,i,\iota}$ and $\mathbf{a}=(a_{t,j})_{t\in\mathbb{Z}_{\geq1},j\in I}\in\mathbb{Z}^{\infty}$,
if $Trop(M')(\mathbf{a})\geq0$ for all $M'\in \mathcal{M}'_{s,i,\iota}$ then
the inequality $Trop(M)(\mathbf{a})\geq0$ follows
by Theorem \ref{prev}.
Therefore, our claim is a consequence of Theorem \ref{polyhthm}. \qed

\section{Classical affine case}

Let $X=A^{(1)}_{n-1},B^{(1)}_{n-1}$, $C^{(1)}_{n-1}$, $D^{(1)}_{n-1}$, $A^{(2)\dagger}_{2n-2}$, $A^{(2)}_{2n-3}$ or
$D^{(2)}_{n}$ and $\iota$ be adapted.

\subsection{Young walls and truncated walls}

Let us review 
Young walls and truncated walls following \cite{Kang} and \cite{Ka24a}.
Each wall consists of $I$-colored blocks of three different shapes:
\begin{enumerate}
\item[(1)] a block with unit width, unit height and unit thickness:
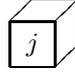
\[
\begin{xy}
(3,3) *{j}="0",
(0,0) *{}="1",
(6,0)*{}="2",
(6,6)*{}="3",
(0,6)*{}="4",
(3,9)*{}="5",
(9,9)*{}="6",
(9,3)*{}="7",
\ar@{-} "1";"2"^{}
\ar@{-} "1";"4"^{}
\ar@{-} "2";"3"^{}
\ar@{-} "3";"4"^{}
\ar@{-} "5";"4"^{}
\ar@{-} "5";"6"^{}
\ar@{-} "3";"6"^{}
\ar@{-} "2";"7"^{}
\ar@{-} "6";"7"^{}
\end{xy}
\]
\item[(2)] a block with unit width, unit height and half-unit thickness:
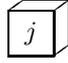
\[
\begin{xy}
(3,3) *{j}="0",
(0,0) *{}="1",
(6,0)*{}="2",
(6,6)*{}="3",
(0,6)*{}="4",
(2,7.5)*{}="5",
(8,7.5)*{}="6",
(8,1.5)*{}="7",
\ar@{-} "1";"2"^{}
\ar@{-} "1";"4"^{}
\ar@{-} "2";"3"^{}
\ar@{-} "3";"4"^{}
\ar@{-} "5";"4"^{}
\ar@{-} "5";"6"^{}
\ar@{-} "3";"6"^{}
\ar@{-} "2";"7"^{}
\ar@{-} "6";"7"^{}
\end{xy}
\]
\item[(3)] a block with unit width, half-unit height and unit thickness:
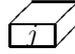
\[
\begin{xy}
(3,1.5) *{j}="0",
(0,0) *{}="1",
(6,0)*{}="2",
(6,3)*{}="3",
(0,3)*{}="4",
(3,6)*{}="5",
(9,6)*{}="6",
(9,3)*{}="7",
\ar@{-} "1";"2"^{}
\ar@{-} "1";"4"^{}
\ar@{-} "2";"3"^{}
\ar@{-} "3";"4"^{}
\ar@{-} "5";"4"^{}
\ar@{-} "5";"6"^{}
\ar@{-} "3";"6"^{}
\ar@{-} "2";"7"^{}
\ar@{-} "6";"7"^{}
\end{xy}
\]
\end{enumerate}
we simply express the block (1) with color $j\in I$
as
\begin{equation}\label{smpl1}
\begin{xy}
(3,3) *{j}="0",
(0,0) *{}="1",
(6,0)*{}="2",
(6,6)*{}="3",
(0,6)*{}="4",
\ar@{-} "1";"2"^{}
\ar@{-} "1";"4"^{}
\ar@{-} "2";"3"^{}
\ar@{-} "3";"4"^{}
\end{xy}
\end{equation}
block (2) with color $j\in I$ as
\begin{equation}\label{smpl12}
\begin{xy}
(1.5,4.5) *{j}="0",
(0,0) *{}="1",
(6,0)*{}="2",
(6,6)*{}="3",
(0,6)*{}="4",
\ar@{-} "1";"4"^{}
\ar@{-} "1";"3"^{}
\ar@{-} "3";"4"^{}
\end{xy}
\end{equation}
or
\begin{equation}\label{smpl11}
\begin{xy}
(4.5,1.5) *{j}="0",
(0,0) *{}="1",
(6,0)*{}="2",
(6,6)*{}="3",
(0,6)*{}="4",
\ar@{-} "1";"2"^{}
\ar@{-} "1";"3"^{}
\ar@{-} "2";"3"^{}
\end{xy}
\end{equation}
and (3) with color $j\in I$ as
\begin{equation}\label{smpl2}
\begin{xy}
(1.5,1.5) *{\ \ j}="0",
(0,0) *{}="1",
(6,0)*{}="2",
(6,3.5)*{}="3",
(0,3.5)*{}="4",
\ar@{-} "1";"2"^{}
\ar@{-} "1";"4"^{}
\ar@{-} "2";"3"^{}
\ar@{-} "3";"4"^{}
\end{xy}
\end{equation}
We call the blocks (\ref{smpl1}), (\ref{smpl12}), (\ref{smpl11}) and (\ref{smpl2}) $j$-blocks.

\begin{defn}\label{classdef}
For $k\in I$ and $X$ other than $C^{(1)}_{n-1}$, the index $k$ is said to be in class $1$ if the fundamental weight $\Lambda_k$ is level $1$.
The index $k\in I$ is said to be in class $2$ if $k$ is not in class $1$. When $X=C^{(1)}_{n-1}$, we understand all indices
$k\in I$ are in class $2$. 
The list of indices $k\in I$ in class $1$ is as follows:
\begin{table}[H]
  \begin{tabular}{|c|c|c|c|c|c|c|} \hline
  $X$ & $A^{(1)}_{n-1}$ & $B^{(1)}_{n-1}$ & $D^{(1)}_{n-1}$ & $A^{(2)\dagger}_{2n-2}$ & $A^{(2)}_{2n-3}$ & $D^{(2)}_{n}$ \\ \hline
  $k$ & $1,2,\cdots,n$ & $1,2,n$ & $1,2,n-1,n$ & $1$ & $1,2$ & $1,n$  \\ \hline
  \end{tabular}
\end{table}

\end{defn}

Let us recall {\it ground state walls} $Y_{\Lambda_k}$ of type $X$ for $k\in I$ in class $1$ \cite{Kang}:
\begin{itemize}
\item For $X=A^{(1)}_{n-1}$ and $\lambda=\Lambda_k$ ($k\in I$), let $Y_{\lambda}$ be the wall that has no blocks.
\item
For $X=A^{(2)\dagger}_{2n-2}$ and $\lambda=\Lambda_1$, we define
\[
Y_{\Lambda_1}:=
\begin{xy}
(-15.5,1.5) *{\ \cdots}="0000",
(-9.5,1.5) *{\ 1}="000",
(-3.5,1.5) *{\ 1}="00",
(1.5,1.5) *{\ \ 1}="0",
(0,0) *{}="1",
(6,0)*{}="2",
(6,3.5)*{}="3",
(0,3.5)*{}="4",
(-6,3.5)*{}="5",
(-6,0)*{}="6",
(-12,3.5)*{}="7",
(-12,0)*{}="8",
\ar@{-} "1";"2"^{}
\ar@{-} "1";"4"^{}
\ar@{-} "2";"3"^{}
\ar@{-} "3";"4"^{}
\ar@{-} "5";"6"^{}
\ar@{-} "5";"4"^{}
\ar@{-} "1";"6"^{}
\ar@{-} "7";"8"^{}
\ar@{-} "7";"5"^{}
\ar@{-} "6";"8"^{}
\end{xy}
\]
The wall $Y_{\Lambda_1}$ has infinitely many $1$-blocks with half-unit height
and
extends infinitely to the left.
\item
For $X=D^{(2)}_{n}$ and $\lambda=\Lambda_1$, $\Lambda_n$, we define
\[
Y_{\Lambda_1}:=
\begin{xy}
(-15.5,1.5) *{\ \cdots}="0000",
(-9.5,1.5) *{\ 1}="000",
(-3.5,1.5) *{\ 1}="00",
(1.5,1.5) *{\ \ 1}="0",
(0,0) *{}="1",
(6,0)*{}="2",
(6,3.5)*{}="3",
(0,3.5)*{}="4",
(-6,3.5)*{}="5",
(-6,0)*{}="6",
(-12,3.5)*{}="7",
(-12,0)*{}="8",
\ar@{-} "1";"2"^{}
\ar@{-} "1";"4"^{}
\ar@{-} "2";"3"^{}
\ar@{-} "3";"4"^{}
\ar@{-} "5";"6"^{}
\ar@{-} "5";"4"^{}
\ar@{-} "1";"6"^{}
\ar@{-} "7";"8"^{}
\ar@{-} "7";"5"^{}
\ar@{-} "6";"8"^{}
\end{xy}
\]
and
\[
Y_{\Lambda_n}:=
\begin{xy}
(-15.5,1.5) *{\ \cdots}="0000",
(-9.5,1.5) *{\ n}="000",
(-3.5,1.5) *{\ n}="00",
(1.5,1.5) *{\ \ n}="0",
(0,0) *{}="1",
(6,0)*{}="2",
(6,3.5)*{}="3",
(0,3.5)*{}="4",
(-6,3.5)*{}="5",
(-6,0)*{}="6",
(-12,3.5)*{}="7",
(-12,0)*{}="8",
\ar@{-} "1";"2"^{}
\ar@{-} "1";"4"^{}
\ar@{-} "2";"3"^{}
\ar@{-} "3";"4"^{}
\ar@{-} "5";"6"^{}
\ar@{-} "5";"4"^{}
\ar@{-} "1";"6"^{}
\ar@{-} "7";"8"^{}
\ar@{-} "7";"5"^{}
\ar@{-} "6";"8"^{}
\end{xy}
\]
\item For $X=A^{(2)}_{2n-3}$ or $B^{(1)}_{n-1}$ and $\lambda=\Lambda_1$, $\Lambda_2$,
\[
Y_{\Lambda_1}:=
\begin{xy}
(3,-1) *{\ \ 2}="2half",
(-2,-1) *{\ 1}="1half2",
(-9,-1) *{\ \ 2}="2hal3",
(-15,-1) *{\ \ 1}="1hal3",
(-22,-1) *{\ \ \cdots}="dot",
(0,-2.5) *{}="1",
(6,-2.5)*{}="2",
(6,3.5)*{}="3",
(0,3.5)*{}="4",
(-6,3.5)*{}="5",
(-6,-2.5)*{}="6",
(-12,-2.5)*{}="8",
(-12,3.5)*{}="8-1",
(-18,-2.5)*{}="9",
\ar@{-} "9";"8-1"^{}
\ar@{-} "8-1";"8"^{}
\ar@{-} "9";"8"^{}
\ar@{-} "5";"8"^{}
\ar@{-} "1";"3"^{}
\ar@{-} "4";"6"^{}
\ar@{-} "1";"2"^{}
\ar@{-} "1";"4"^{}
\ar@{-} "2";"3"^{}
\ar@{-} "5";"6"^{}
\ar@{-} "1";"6"^{}
\ar@{-} "6";"8"^{}
\end{xy}\qquad
Y_{\Lambda_2}:=
\begin{xy}
(3,-1) *{\ \ 1}="2half",
(-2,-1) *{\ 2}="1half2",
(-9,-1) *{\ \ 1}="2hal3",
(-15,-1) *{\ \ 2}="1hal3",
(-22,-1) *{\ \ \cdots}="dot",
(0,-2.5) *{}="1",
(6,-2.5)*{}="2",
(6,3.5)*{}="3",
(0,3.5)*{}="4",
(-6,3.5)*{}="5",
(-6,-2.5)*{}="6",
(-12,-2.5)*{}="8",
(-12,3.5)*{}="8-1",
(-18,-2.5)*{}="9",
\ar@{-} "9";"8-1"^{}
\ar@{-} "8-1";"8"^{}
\ar@{-} "9";"8"^{}
\ar@{-} "5";"8"^{}
\ar@{-} "1";"3"^{}
\ar@{-} "4";"6"^{}
\ar@{-} "1";"2"^{}
\ar@{-} "1";"4"^{}
\ar@{-} "2";"3"^{}
\ar@{-} "5";"6"^{}
\ar@{-} "1";"6"^{}
\ar@{-} "6";"8"^{}
\end{xy}
\]
The walls $Y_{\Lambda_1}$, $Y_{\Lambda_2}$ extend infinitely to the left. For
$X=B^{(1)}_{n-1}$ and $\lambda=\Lambda_n$,
we define
\[
Y_{\Lambda_n}:=
\begin{xy}
(-15.5,1.5) *{\ \cdots}="0000",
(-9.5,1.5) *{\ n}="000",
(-3.5,1.5) *{\ n}="00",
(1.5,1.5) *{\ \ n}="0",
(0,0) *{}="1",
(6,0)*{}="2",
(6,3.5)*{}="3",
(0,3.5)*{}="4",
(-6,3.5)*{}="5",
(-6,0)*{}="6",
(-12,3.5)*{}="7",
(-12,0)*{}="8",
\ar@{-} "1";"2"^{}
\ar@{-} "1";"4"^{}
\ar@{-} "2";"3"^{}
\ar@{-} "3";"4"^{}
\ar@{-} "5";"6"^{}
\ar@{-} "5";"4"^{}
\ar@{-} "1";"6"^{}
\ar@{-} "7";"8"^{}
\ar@{-} "7";"5"^{}
\ar@{-} "6";"8"^{}
\end{xy}
\]

\item For $X=D^{(1)}_{n-1}$ and $\lambda=\Lambda_1$, $\Lambda_2$, $\Lambda_{n-1}$ or $\Lambda_n$,
\[
Y_{\Lambda_1}:=
\begin{xy}
(3,-1) *{\ \ 2}="2half",
(-2,-1) *{\ 1}="1half2",
(-9,-1) *{\ \ 2}="2hal3",
(-15,-1) *{\ \ 1}="1hal3",
(-22,-1) *{\ \ \cdots}="dot",
(0,-2.5) *{}="1",
(6,-2.5)*{}="2",
(6,3.5)*{}="3",
(0,3.5)*{}="4",
(-6,3.5)*{}="5",
(-6,-2.5)*{}="6",
(-12,-2.5)*{}="8",
(-12,3.5)*{}="8-1",
(-18,-2.5)*{}="9",
\ar@{-} "9";"8-1"^{}
\ar@{-} "8-1";"8"^{}
\ar@{-} "9";"8"^{}
\ar@{-} "5";"8"^{}
\ar@{-} "1";"3"^{}
\ar@{-} "4";"6"^{}
\ar@{-} "1";"2"^{}
\ar@{-} "1";"4"^{}
\ar@{-} "2";"3"^{}
\ar@{-} "5";"6"^{}
\ar@{-} "1";"6"^{}
\ar@{-} "6";"8"^{}
\end{xy}\qquad
Y_{\Lambda_2}:=
\begin{xy}
(3,-1) *{\ \ 1}="2half",
(-2,-1) *{\ 2}="1half2",
(-9,-1) *{\ \ 1}="2hal3",
(-15,-1) *{\ \ 2}="1hal3",
(-22,-1) *{\ \ \cdots}="dot",
(0,-2.5) *{}="1",
(6,-2.5)*{}="2",
(6,3.5)*{}="3",
(0,3.5)*{}="4",
(-6,3.5)*{}="5",
(-6,-2.5)*{}="6",
(-12,-2.5)*{}="8",
(-12,3.5)*{}="8-1",
(-18,-2.5)*{}="9",
\ar@{-} "9";"8-1"^{}
\ar@{-} "8-1";"8"^{}
\ar@{-} "9";"8"^{}
\ar@{-} "5";"8"^{}
\ar@{-} "1";"3"^{}
\ar@{-} "4";"6"^{}
\ar@{-} "1";"2"^{}
\ar@{-} "1";"4"^{}
\ar@{-} "2";"3"^{}
\ar@{-} "5";"6"^{}
\ar@{-} "1";"6"^{}
\ar@{-} "6";"8"^{}
\end{xy}
\]
\[
Y_{\Lambda_{n-1}}:=
\begin{xy}
(3,-1) *{\ \ n}="2half",
(-2,-1) *{_{n-1}}="1half2",
(-9,-1) *{\ \ n}="2hal3",
(-15,-1) *{\ _{n-1}}="1hal3",
(-22,-1) *{\ \ \cdots}="dot",
(0,-2.5) *{}="1",
(6,-2.5)*{}="2",
(6,3.5)*{}="3",
(0,3.5)*{}="4",
(-6,3.5)*{}="5",
(-6,-2.5)*{}="6",
(-12,-2.5)*{}="8",
(-12,3.5)*{}="8-1",
(-18,-2.5)*{}="9",
\ar@{-} "9";"8-1"^{}
\ar@{-} "8-1";"8"^{}
\ar@{-} "9";"8"^{}
\ar@{-} "5";"8"^{}
\ar@{-} "1";"3"^{}
\ar@{-} "4";"6"^{}
\ar@{-} "1";"2"^{}
\ar@{-} "1";"4"^{}
\ar@{-} "2";"3"^{}
\ar@{-} "5";"6"^{}
\ar@{-} "1";"6"^{}
\ar@{-} "6";"8"^{}
\end{xy}\qquad
Y_{\Lambda_n}:=
\begin{xy}
(3,-1) *{\ \ _{n-1}}="2half",
(-2,-1) *{\ n}="1half2",
(-9,-1) *{\ \ _{n-1}}="2hal3",
(-15,-1) *{\ \ n}="1hal3",
(-22,-1) *{\ \ \cdots}="dot",
(0,-2.5) *{}="1",
(6,-2.5)*{}="2",
(6,3.5)*{}="3",
(0,3.5)*{}="4",
(-6,3.5)*{}="5",
(-6,-2.5)*{}="6",
(-12,-2.5)*{}="8",
(-12,3.5)*{}="8-1",
(-18,-2.5)*{}="9",
\ar@{-} "9";"8-1"^{}
\ar@{-} "8-1";"8"^{}
\ar@{-} "9";"8"^{}
\ar@{-} "5";"8"^{}
\ar@{-} "1";"3"^{}
\ar@{-} "4";"6"^{}
\ar@{-} "1";"2"^{}
\ar@{-} "1";"4"^{}
\ar@{-} "2";"3"^{}
\ar@{-} "5";"6"^{}
\ar@{-} "1";"6"^{}
\ar@{-} "6";"8"^{}
\end{xy}
\]
\end{itemize}

\begin{defn}\cite{Kang}\label{def-YW1}
Let $k\in I$ in class $1$. If the following holds then
a wall $Y$ is said to be {\it Young wall} of ground state $\Lambda_k$ of type $X$:
\begin{enumerate}
\item The wall $Y$ is obtained from $Y_{\Lambda_k}$
by stacking finitely many colored blocks on the ground-state wall $Y_{\Lambda_k}$.
\item There are no blocks on the top of a single block with half-unit thickness.
\item Except for the right-most column,
there are no free spaces to the
right of all blocks.
\item The blocks are stacked following the patterns we give below. %in Sect.5 of \cite{Kang}.
\end{enumerate}
\end{defn}
For each type $X$ and $\Lambda_k$ with $k\in I$ in class $1$, 
the patterns of (iv) are as follows: 
\[
\begin{xy}
(-24.5,50) *{A^{(1)}_{n-1},}="type",
(-24.5,45) *{\Lambda_k:}="type2",
(-15.5,0) *{\ _{k-3}}="000-1",
(-9.5,0) *{\ _{k-2}}="00-1",
(-3.5,0) *{\ _{k-1}}="0-1",
(1.5,0) *{\ \ k}="0-1",
%(-15.5,1.5) *{\ 1}="0000",
%(-9.5,1.5) *{\ 1}="000",
%(-3.5,1.5) *{\ 1}="00",
%(1.5,1.5) *{\ \ 1}="0",
(1.5,6.5) *{\ \ _{k+1}}="02",
(-4,6.5) *{\ \ k}="002",
(-10,6.5) *{\ \ _{k-1}}="0002",
(-16,6.5) *{\ \ _{k-2}}="00002",
(1.5,12.5) *{\ \ _{k+2}}="03",
(-4,12.5) *{\ \ _{k+1}}="003",
(-10,12.5) *{\ \ k}="0003",
(-16,12.5) *{\ \ _{k-1}}="00003",
(1.5,20.5) *{\ \ \vdots}="04",
(-4,20.5) *{\ \ \vdots}="004",
(-10,20.5) *{\ \ \vdots}="0004",
(-16,20.5) *{\ \ \vdots}="00004",
(1.5,27.5) *{\ \ n}="05",
(-4,27.5) *{\ \ _{n-1}}="005",
(-10,27.5) *{\ \ _{n-2}}="0005",
(-16,27.5) *{\ \ _{n-3}}="00005",
(1.5,33.5) *{\ \ 1}="06",
(-4,33.5) *{\ \ n}="006",
(-10,33.5) *{\ \ _{n-1}}="0006",
(-16,33.5) *{\ \ _{n-2}}="00006",
(1.5,39.5) *{\ \ 2}="07",
(-4,39.5) *{\ \ 1}="007",
(-10,39.5) *{\ \ n}="0007",
(-16,39.5) *{\ \ _{n-1}}="00007",
(1.5,48) *{\ \ \vdots}="08",
(-4,48) *{\ \ \vdots}="008",
(-10,48) *{\ \ \vdots}="0008",
(-16,48) *{\ \ \vdots}="00008",
(1.5,55) *{\ \ n}="09",
(-4,55) *{\ \ _{n-1}}="009",
(-10,55) *{\ \ _{n-2}}="0009",
(-16,55) *{\ \ _{n-3}}="00009",
%(1.5,60) *{\ \ 1}="010",
%(-4,60) *{\ \ 1}="0010",
%(-10,60) *{\ \ 1}="00010",
%(-16,60) *{\ \ 1}="000010",
(-18,-3.5) *{}="0.5-l",
(-20,-3.5) *{}="0-l",
(-20,0) *{}="1-l",
(-20,3.5) *{}="2-l",
(-20,9.5) *{}="3-l",
(-20,15.5) *{}="4-l",
(-20,24.5) *{}="5-l",
(-20,30.5) *{}="6-l",
(-20,36.5) *{}="7-l",
(-20,42.5) *{}="8-l",
(-20,52) *{}="9-l",
(-20,58) *{}="10-l",
(-20,62) *{}="11-l",
(-20,66) *{}="12-l",
(-20,60) *{}="13-l",
(-18,60) *{}="13.5-l",
(6,-3.5)*{}="r--1",
(6,0)*{}="r",
(6,3.5)*{}="r-0",
(0,-3.5) *{}="b1",
(-6,-3.5)*{}="b2",
(-12,-3.5)*{}="b3",
(0,60) *{}="t1",
(-6,60)*{}="t2",
(-12,60)*{}="t3",
(6,9.5)*{}="r-1",
(6,15.5)*{}="r-2",
(6,24.5)*{}="r-3",
(6,30.5)*{}="r-4",
(6,36.5)*{}="r-5",
(6,42.5)*{}="r-6",
(6,52)*{}="r-7",
(6,58)*{}="r-8",
(6,62)*{}="r-9",
%(6,66)*{}="r-10",
%(6,72)*{}="r-11",
(6,60)*{}="r-11.5",
\ar@{-} "b1";"t1"^{}
\ar@{-} "b2";"t2"^{}
\ar@{-} "b3";"t3"^{}
\ar@{-} "0.5-l";"13.5-l"^{}
\ar@{-} "r-11.5";"r--1"^{}
\ar@{-} "0-l";"r--1"^{}
%\ar@{-} "1-l";"r"^{}
\ar@{-} "2-l";"r-0"^{}
\ar@{-} "3-l";"r-1"^{}
\ar@{-} "4-l";"r-2"^{}
\ar@{-} "5-l";"r-3"^{}
\ar@{-} "6-l";"r-4"^{}
\ar@{-} "7-l";"r-5"^{}
\ar@{-} "8-l";"r-6"^{}
\ar@{-} "9-l";"r-7"^{}
\ar@{-} "10-l";"r-8"^{}
%\ar@{-} "11-l";"r-9"^{}
%\ar@{-} "12-l";"r-10"^{}
%\ar@{-} "13-l";"r-11"^{}
\end{xy}
\]
For $j_1,j_2$ such that $\{j_1,j_2\}=\{1,2\}$, 
\[
\begin{xy}
(-24.5,50) *{A^{(2)}_{2n-3},}="type",
(-24.5,45) *{\Lambda_{j_1}:}="type2",
(-15.5,-2) *{\ \ \ \ _{j_2}}="000-1",
(-9.5,-2) *{\ \ \ _{j_1}}="00-1",
(-3.5,-2) *{\ \ \ _{j_2}}="0-1",
(1.5,-2) *{\ \ \ \ _{j_1}}="0-1",
(-15.5,1.5) *{_{j_1}}="0000",
(-9.5,1.5) *{_{j_2}}="000",
(-3.5,1.5) *{ _{j_1}}="00",
(1.5,1.5) *{\ _{j_2}}="0",
(1.5,6.5) *{\ \ 3}="02",
(-4,6.5) *{\ \ 3}="002",
(-10,6.5) *{\ \ 3}="0002",
(-16,6.5) *{\ \ 3}="00002",
(1.5,12.5) *{\ \ 4}="03",
(-4,12.5) *{\ \ 4}="003",
(-10,12.5) *{\ \ 4}="0003",
(-16,12.5) *{\ \ 4}="00003",
(1.5,20.5) *{\ \ \vdots}="04",
(-4,20.5) *{\ \ \vdots}="004",
(-10,20.5) *{\ \ \vdots}="0004",
(-16,20.5) *{\ \ \vdots}="00004",
(1.5,27.5) *{\ \ _{n-1}}="05",
(-4,27.5) *{\ \ _{n-1}}="005",
(-10,27.5) *{\ \ _{n-1}}="0005",
(-16,27.5) *{\ \ _{n-1}}="00005",
(1.5,33.5) *{\ \ n}="06",
(-4,33.5) *{\ \ n}="006",
(-10,33.5) *{\ \ n}="0006",
(-16,33.5) *{\ \ n}="00006",
(1.5,39.5) *{\ \ _{n-1}}="07",
(-4,39.5) *{\ \ _{n-1}}="007",
(-10,39.5) *{\ \ _{n-1}}="0007",
(-16,39.5) *{\ \ _{n-1}}="00007",
(1.5,48) *{\ \ \vdots}="08",
(-4,48) *{\ \ \vdots}="008",
(-10,48) *{\ \ \vdots}="0008",
(-16,48) *{\ \ \vdots}="00008",
(1.5,55) *{\ \ 3}="09",
(-4,55) *{\ \ 3}="009",
(-10,55) *{\ \ 3}="0009",
(-16,55) *{\ \ 3}="00009",
(1.5,60) *{\ \ \ \ _{j_1}}="010",
(-4,60) *{\ \ \ \ _{j_2}}="0010",
(-10,60) *{\ \ \ \ _{j_1}}="00010",
(-16,60) *{\ \ \ \ _{j_2}}="000010",
(1.5,64) *{_{j_2}}="011",
(-4,64) *{_{j_1}}="0011",
(-10,64) *{_{j_2}}="00011",
(-16,64) *{_{j_1}}="000011",
(1.5,69) *{\ \ 3}="012",
(-4,69) *{\ \ 3}="0012",
(-10,69) *{\ \ 3}="00012",
(-16,69) *{\ \ 3}="000012",
(-18,-3.5) *{}="0.5-l",
(-20,-3.5) *{}="0-l",
(-20,0) *{}="1-l",
(-20,3.5) *{}="2-l",
(-20,9.5) *{}="3-l",
(-20,15.5) *{}="4-l",
(-20,24.5) *{}="5-l",
(-20,30.5) *{}="6-l",
(-20,36.5) *{}="7-l",
(-20,42.5) *{}="8-l",
(-20,52) *{}="9-l",
(-20,58) *{}="10-l",
(-20,62) *{}="11-l",
(-20,66) *{}="12-l",
(-20,72) *{}="13-l",
(-18,74) *{}="13.5-l",
(6,-3.5)*{}="r--1",
(6,0)*{}="r",
(6,3.5)*{}="r-0",
(0,-3.5) *{}="b1",
(-6,-3.5)*{}="b2",
(-12,-3.5)*{}="b3",
(0,74) *{}="t1",
(-6,74)*{}="t2",
(-12,74)*{}="t3",
(6,9.5)*{}="r-1",
(6,15.5)*{}="r-2",
(6,24.5)*{}="r-3",
(6,30.5)*{}="r-4",
(6,36.5)*{}="r-5",
(6,42.5)*{}="r-6",
(6,52)*{}="r-7",
(6,58)*{}="r-8",
(6,62)*{}="r-9",
(6,66)*{}="r-10",
(6,72)*{}="r-11",
(6,74)*{}="r-11.5",
(0,3.5)*{}="x1",
(-6,-3.5)*{}="x2",
(-6,3.5)*{}="x3",
(-12,-3.5)*{}="x4",
(-12,3.5)*{}="x5",
(-18,-3.5)*{}="x6",
(6,66)*{}="y-1",
(0,58)*{}="y0",
(0,66)*{}="y1",
(-6,58)*{}="y2",
(-6,66)*{}="y3",
(-12,58)*{}="y4",
(-12,66)*{}="y5",
(-18,58)*{}="y6",
\ar@{-} "b1";"t1"^{}
\ar@{-} "b2";"t2"^{}
\ar@{-} "b3";"t3"^{}
\ar@{-} "0.5-l";"13.5-l"^{}
\ar@{-} "r-11.5";"r--1"^{}
\ar@{-} "0-l";"r--1"^{}
%\ar@{-} "1-l";"r"^{}
\ar@{-} "2-l";"r-0"^{}
\ar@{-} "3-l";"r-1"^{}
\ar@{-} "4-l";"r-2"^{}
\ar@{-} "5-l";"r-3"^{}
\ar@{-} "6-l";"r-4"^{}
\ar@{-} "7-l";"r-5"^{}
\ar@{-} "8-l";"r-6"^{}
\ar@{-} "9-l";"r-7"^{}
\ar@{-} "10-l";"r-8"^{}
%\ar@{-} "11-l";"r-9"^{}
\ar@{-} "12-l";"r-10"^{}
\ar@{-} "13-l";"r-11"^{}
\ar@{-} "b1";"r-0"^{}
\ar@{-} "x1";"x2"^{}
\ar@{-} "x3";"x4"^{}
\ar@{-} "x5";"x6"^{}
\ar@{-} "y-1";"y0"^{}
\ar@{-} "y1";"y2"^{}
\ar@{-} "y3";"y4"^{}
\ar@{-} "y5";"y6"^{}
\end{xy}\quad
\begin{xy}
(-24.5,50) *{B^{(1)}_{n-1},}="type",
(-24.5,45) *{\Lambda_{j_1}:}="type2",
(-15.5,-2) *{\ \ \ \ _{j_2}}="000-1",
(-9.5,-2) *{\ \ \ _{j_1}}="00-1",
(-3.5,-2) *{\ \ \ _{j_2}}="0-1",
(1.5,-2) *{\ \ \ \ _{j_1}}="0-1",
(-15.5,1.5) *{_{j_1}}="0000",
(-9.5,1.5) *{_{j_2}}="000",
(-3.5,1.5) *{ _{j_1}}="00",
(1.5,1.5) *{\ _{j_2}}="0",
(1.5,6.5) *{\ \ 3}="02",
(-4,6.5) *{\ \ 3}="002",
(-10,6.5) *{\ \ 3}="0002",
(-16,6.5) *{\ \ 3}="00002",
(1.5,12.5) *{\ \ 4}="03",
(-4,12.5) *{\ \ 4}="003",
(-10,12.5) *{\ \ 4}="0003",
(-16,12.5) *{\ \ 4}="00003",
(1.5,20.5) *{\ \ \vdots}="04",
(-4,20.5) *{\ \ \vdots}="004",
(-10,20.5) *{\ \ \vdots}="0004",
(-16,20.5) *{\ \ \vdots}="00004",
(1.5,27.5) *{\ \ _{n-1}}="05",
(-4,27.5) *{\ \ _{n-1}}="005",
(-10,27.5) *{\ \ _{n-1}}="0005",
(-16,27.5) *{\ \ _{n-1}}="00005",
(1.5,32) *{\ \ n}="06",
(-4,32) *{\ \ n}="006",
(-10,32) *{\ \ n}="0006",
(-16,32) *{\ \ n}="00006",
(1.5,35) *{\ \ n}="06-2",
(-4,35) *{\ \ n}="006-2",
(-10,35) *{\ \ n}="0006-2",
(-16,35) *{\ \ n}="00006-2",
(1.5,39.5) *{\ \ _{n-1}}="07",
(-4,39.5) *{\ \ _{n-1}}="007",
(-10,39.5) *{\ \ _{n-1}}="0007",
(-16,39.5) *{\ \ _{n-1}}="00007",
(1.5,48) *{\ \ \vdots}="08",
(-4,48) *{\ \ \vdots}="008",
(-10,48) *{\ \ \vdots}="0008",
(-16,48) *{\ \ \vdots}="00008",
(1.5,55) *{\ \ 3}="09",
(-4,55) *{\ \ 3}="009",
(-10,55) *{\ \ 3}="0009",
(-16,55) *{\ \ 3}="00009",
(1.5,60) *{\ \ \ \ _{j_1}}="010",
(-4,60) *{\ \ \ \ _{j_2}}="0010",
(-10,60) *{\ \ \ \ _{j_1}}="00010",
(-16,60) *{\ \ \ \ _{j_2}}="000010",
(1.5,64) *{_{j_2}}="011",
(-4,64) *{_{j_1}}="0011",
(-10,64) *{_{j_2}}="00011",
(-16,64) *{_{j_1}}="000011",
(1.5,69) *{\ \ 3}="012",
(-4,69) *{\ \ 3}="0012",
(-10,69) *{\ \ 3}="00012",
(-16,69) *{\ \ 3}="000012",
(-18,-3.5) *{}="0.5-l",
(-20,-3.5) *{}="0-l",
(-20,0) *{}="1-l",
(-20,3.5) *{}="2-l",
(-20,9.5) *{}="3-l",
(-20,15.5) *{}="4-l",
(-20,24.5) *{}="5-l",
(-20,30.5) *{}="6-l",
(-20,33.5) *{}="6-7-l",
(-20,36.5) *{}="7-l",
(-20,42.5) *{}="8-l",
(-20,52) *{}="9-l",
(-20,58) *{}="10-l",
(-20,62) *{}="11-l",
(-20,66) *{}="12-l",
(-20,72) *{}="13-l",
(-18,74) *{}="13.5-l",
(6,-3.5)*{}="r--1",
(6,0)*{}="r",
(6,3.5)*{}="r-0",
(0,-3.5) *{}="b1",
(-6,-3.5)*{}="b2",
(-12,-3.5)*{}="b3",
(0,74) *{}="t1",
(-6,74)*{}="t2",
(-12,74)*{}="t3",
(6,9.5)*{}="r-1",
(6,15.5)*{}="r-2",
(6,24.5)*{}="r-3",
(6,30.5)*{}="r-4",
(6,33.5)*{}="6-7-r",
(6,36.5)*{}="r-5",
(6,42.5)*{}="r-6",
(6,52)*{}="r-7",
(6,58)*{}="r-8",
(6,62)*{}="r-9",
(6,66)*{}="r-10",
(6,72)*{}="r-11",
(6,74)*{}="r-11.5",
(0,3.5)*{}="x1",
(-6,-3.5)*{}="x2",
(-6,3.5)*{}="x3",
(-12,-3.5)*{}="x4",
(-12,3.5)*{}="x5",
(-18,-3.5)*{}="x6",
(6,66)*{}="y-1",
(0,58)*{}="y0",
(0,66)*{}="y1",
(-6,58)*{}="y2",
(-6,66)*{}="y3",
(-12,58)*{}="y4",
(-12,66)*{}="y5",
(-18,58)*{}="y6",
\ar@{-} "6-7-l";"6-7-r"^{}
\ar@{-} "b1";"t1"^{}
\ar@{-} "b2";"t2"^{}
\ar@{-} "b3";"t3"^{}
\ar@{-} "0.5-l";"13.5-l"^{}
\ar@{-} "r-11.5";"r--1"^{}
\ar@{-} "0-l";"r--1"^{}
%\ar@{-} "1-l";"r"^{}
\ar@{-} "2-l";"r-0"^{}
\ar@{-} "3-l";"r-1"^{}
\ar@{-} "4-l";"r-2"^{}
\ar@{-} "5-l";"r-3"^{}
\ar@{-} "6-l";"r-4"^{}
\ar@{-} "7-l";"r-5"^{}
\ar@{-} "8-l";"r-6"^{}
\ar@{-} "9-l";"r-7"^{}
\ar@{-} "10-l";"r-8"^{}
%\ar@{-} "11-l";"r-9"^{}
\ar@{-} "12-l";"r-10"^{}
\ar@{-} "13-l";"r-11"^{}
\ar@{-} "b1";"r-0"^{}
\ar@{-} "x1";"x2"^{}
\ar@{-} "x3";"x4"^{}
\ar@{-} "x5";"x6"^{}
\ar@{-} "y-1";"y0"^{}
\ar@{-} "y1";"y2"^{}
\ar@{-} "y3";"y4"^{}
\ar@{-} "y5";"y6"^{}
\end{xy}\quad
\begin{xy}
(-24.5,50) *{B^{(1)}_{n-1},}="type",
(-24.5,45) *{\Lambda_n:}="type2",
(-15.5,-2) *{\ n}="000-1",
(-9.5,-2) *{\ n}="00-1",
(-3.5,-2) *{\ n}="0-1",
(1.5,-2) *{\ \ n}="0-1",
(-15.5,1.5) *{\ n}="0000",
(-9.5,1.5) *{\ n}="000",
(-3.5,1.5) *{\ n}="00",
(1.5,1.5) *{\ \ n}="0",
(1.5,6.5) *{\ \ _{n-1}}="02",
(-4,6.5) *{\ \ _{n-1}}="002",
(-10,6.5) *{\ \ _{n-1}}="0002",
(-16,6.5) *{\ \ _{n-1}}="00002",
(1.5,12.5) *{\ \ _{n-2}}="03",
(-4,12.5) *{\ \ _{n-2}}="003",
(-10,12.5) *{\ \ _{n-2}}="0003",
(-16,12.5) *{\ \ _{n-2}}="00003",
(1.5,20.5) *{\ \ \vdots}="04",
(-4,20.5) *{\ \ \vdots}="004",
(-10,20.5) *{\ \ \vdots}="0004",
(-16,20.5) *{\ \ \vdots}="00004",
(1.5,27.5) *{\ \ 3}="05",
(-4,27.5) *{\ \ 3}="005",
(-10,27.5) *{\ \ 3}="0005",
(-16,27.5) *{\ \ 3}="00005",
(1.5,32) *{\ \ \ \ _{j_1}}="06",
(-4,32) *{\ \ \ \ _{j_2}}="006",
(-10,32) *{\ \ \ \ _{j_1}}="0006",
(-16,32) *{\ \ \ \ _{j_2}}="00006",
(1.5,35) *{_{j_2}}="06-2",
(-4,35) *{_{j_1}}="006-2",
(-10,35) *{_{j_2}}="0006-2",
(-16,35) *{_{j_1}}="00006-2",
(1.5,39.5) *{\ \ 3}="07",
(-4,39.5) *{\ \ 3}="007",
(-10,39.5) *{\ \ 3}="0007",
(-16,39.5) *{\ \ 3}="00007",
(1.5,48) *{\ \ \vdots}="08",
(-4,48) *{\ \ \vdots}="008",
(-10,48) *{\ \ \vdots}="0008",
(-16,48) *{\ \ \vdots}="00008",
(1.5,55) *{\ \ _{n-1}}="09",
(-4,55) *{\ \ _{n-1}}="009",
(-10,55) *{\ \ _{n-1}}="0009",
(-16,55) *{\ \ _{n-1}}="00009",
(1.5,60) *{\ \ n}="010",
(-4,60) *{\ \ n}="0010",
(-10,60) *{\ \ n}="00010",
(-16,60) *{\ \ n}="000010",
(1.5,64) *{\ \ n}="011",
(-4,64) *{\ \ n}="0011",
(-10,64) *{\ \ n}="00011",
(-16,64) *{\ \ n}="000011",
(1.5,69) *{\ \ _{n-1}}="012",
(-4,69) *{\ \ _{n-1}}="0012",
(-10,69) *{\ \ _{n-1}}="00012",
(-16,69) *{\ \ _{n-1}}="000012",
(-18,-3.5) *{}="0.5-l",
(-20,-3.5) *{}="0-l",
(-20,0) *{}="1-l",
(-20,3.5) *{}="2-l",
(-20,9.5) *{}="3-l",
(-20,15.5) *{}="4-l",
(-20,24.5) *{}="5-l",
(-20,30.5) *{}="6-l",
(-20,33.5) *{}="6-7-l",
(-20,36.5) *{}="7-l",
(-20,42.5) *{}="8-l",
(-20,52) *{}="9-l",
(-20,58) *{}="10-l",
(-20,62) *{}="11-l",
(-20,66) *{}="12-l",
(-20,72) *{}="13-l",
(-18,74) *{}="13.5-l",
(6,-3.5)*{}="r--1",
(6,0)*{}="r",
(6,3.5)*{}="r-0",
(0,-3.5) *{}="b1",
(-6,-3.5)*{}="b2",
(-12,-3.5)*{}="b3",
(0,74) *{}="t1",
(-6,74)*{}="t2",
(-12,74)*{}="t3",
(6,9.5)*{}="r-1",
(6,15.5)*{}="r-2",
(6,24.5)*{}="r-3",
(6,30.5)*{}="r-4",
(6,33.5)*{}="6-7-r",
(6,36.5)*{}="r-5",
(6,42.5)*{}="r-6",
(6,52)*{}="r-7",
(6,58)*{}="r-8",
(6,62)*{}="r-9",
(6,66)*{}="r-10",
(6,72)*{}="r-11",
(6,74)*{}="r-11.5",
(6,36.5)*{}="x1",
(0,30.5)*{}="x2",
(0,36.5)*{}="x3",
(-6,30.5)*{}="x4",
(-6,36.5)*{}="x5",
(-12,30.5)*{}="x6",
(-12,36.5)*{}="x7",
(-18,30.5)*{}="x8",
%\ar@{-} "6-7-l";"6-7-r"^{}
\ar@{-} "b1";"t1"^{}
\ar@{-} "b2";"t2"^{}
\ar@{-} "b3";"t3"^{}
\ar@{-} "0.5-l";"13.5-l"^{}
\ar@{-} "r-11.5";"r--1"^{}
\ar@{-} "0-l";"r--1"^{}
\ar@{-} "1-l";"r"^{}
\ar@{-} "2-l";"r-0"^{}
\ar@{-} "3-l";"r-1"^{}
\ar@{-} "4-l";"r-2"^{}
\ar@{-} "5-l";"r-3"^{}
\ar@{-} "6-l";"r-4"^{}
\ar@{-} "7-l";"r-5"^{}
\ar@{-} "8-l";"r-6"^{}
\ar@{-} "9-l";"r-7"^{}
\ar@{-} "10-l";"r-8"^{}
\ar@{-} "11-l";"r-9"^{}
\ar@{-} "12-l";"r-10"^{}
\ar@{-} "13-l";"r-11"^{}
\ar@{-} "x1";"x2"^{}
\ar@{-} "x3";"x4"^{}
\ar@{-} "x5";"x6"^{}
\ar@{-} "x7";"x8"^{}
\end{xy}
\]

Let $\{j_1,j_2\}=\{1,2\}$ or $\{n-1,n\}$, $\{j_3,j_4\}=\{n-1,n\}$ or $\{1,2\}$,
\[
\begin{xy}
(-24.5,50) *{D^{(1)}_{n-1},}="type",
(-24.5,45) *{\Lambda_{j_1}:}="type2",
(-15.5,-2) *{\ \ \ \ _{j_2}}="000-1",
(-9.5,-2) *{\ \ \ _{j_1}}="00-1",
(-3.5,-2) *{\ \ \ _{j_2}}="0-1",
(1.5,-2) *{\ \ \ \ _{j_1}}="0-1",
(-15.5,1.5) *{_{j_1}}="0000",
(-9.5,1.5) *{_{j_2}}="000",
(-3.5,1.5) *{ _{j_1}}="00",
(1.5,1.5) *{\ _{j_2}}="0",
(1.5,6.5) *{\ \ 3}="02",
(-4,6.5) *{\ \ 3}="002",
(-10,6.5) *{\ \ 3}="0002",
(-16,6.5) *{\ \ 3}="00002",
(1.5,12.5) *{\ \ 4}="03",
(-4,12.5) *{\ \ 4}="003",
(-10,12.5) *{\ \ 4}="0003",
(-16,12.5) *{\ \ 4}="00003",
(1.5,20.5) *{\ \ \vdots}="04",
(-4,20.5) *{\ \ \vdots}="004",
(-10,20.5) *{\ \ \vdots}="0004",
(-16,20.5) *{\ \ \vdots}="00004",
(1.5,27.5) *{\ \ _{n-2}}="05",
(-4,27.5) *{\ \ _{n-2}}="005",
(-10,27.5) *{\ \ _{n-2}}="0005",
(-16,27.5) *{\ \ _{n-2}}="00005",
(1.5,32) *{\ \ \ \ _{j_3}}="06",
(-4,32) *{\ \ \ \ _{j_4}}="006",
(-10,32) *{\ \ \ \ _{j_3}}="0006",
(-16,32) *{\ \ \ \ _{j_4}}="00006",
(1.5,35) *{_{j_4}}="06-2",
(-4,35) *{_{j_3}}="006-2",
(-10,35) *{_{j_4}}="0006-2",
(-16,35) *{_{j_3}}="00006-2",
(1.5,39.5) *{\ \ _{n-2}}="07",
(-4,39.5) *{\ \ _{n-2}}="007",
(-10,39.5) *{\ \ _{n-2}}="0007",
(-16,39.5) *{\ \ _{n-2}}="00007",
(1.5,48) *{\ \ \vdots}="08",
(-4,48) *{\ \ \vdots}="008",
(-10,48) *{\ \ \vdots}="0008",
(-16,48) *{\ \ \vdots}="00008",
(1.5,55) *{\ \ 3}="09",
(-4,55) *{\ \ 3}="009",
(-10,55) *{\ \ 3}="0009",
(-16,55) *{\ \ 3}="00009",
(1.5,60) *{\ \ \ \ _{j_1}}="010",
(-4,60) *{\ \ \ \ _{j_2}}="0010",
(-10,60) *{\ \ \ \ _{j_1}}="00010",
(-16,60) *{\ \ \ \ _{j_2}}="000010",
(1.5,64) *{_{j_2}}="011",
(-4,64) *{_{j_1}}="0011",
(-10,64) *{_{j_2}}="00011",
(-16,64) *{_{j_1}}="000011",
(1.5,69) *{\ \ 3}="012",
(-4,69) *{\ \ 3}="0012",
(-10,69) *{\ \ 3}="00012",
(-16,69) *{\ \ 3}="000012",
(-18,-3.5) *{}="0.5-l",
(-20,-3.5) *{}="0-l",
(-20,0) *{}="1-l",
(-20,3.5) *{}="2-l",
(-20,9.5) *{}="3-l",
(-20,15.5) *{}="4-l",
(-20,24.5) *{}="5-l",
(-20,30.5) *{}="6-l",
(-20,36.5) *{}="7-l",
(-20,42.5) *{}="8-l",
(-20,52) *{}="9-l",
(-20,58) *{}="10-l",
(-20,62) *{}="11-l",
(-20,66) *{}="12-l",
(-20,72) *{}="13-l",
(-18,74) *{}="13.5-l",
(6,-3.5)*{}="r--1",
(6,0)*{}="r",
(6,3.5)*{}="r-0",
(0,-3.5) *{}="b1",
(-6,-3.5)*{}="b2",
(-12,-3.5)*{}="b3",
(0,74) *{}="t1",
(-6,74)*{}="t2",
(-12,74)*{}="t3",
(6,9.5)*{}="r-1",
(6,15.5)*{}="r-2",
(6,24.5)*{}="r-3",
(6,30.5)*{}="r-4",
(6,36.5)*{}="r-5",
(6,42.5)*{}="r-6",
(6,52)*{}="r-7",
(6,58)*{}="r-8",
(6,62)*{}="r-9",
(6,66)*{}="r-10",
(6,72)*{}="r-11",
(6,74)*{}="r-11.5",
(0,3.5)*{}="x1",
(-6,-3.5)*{}="x2",
(-6,3.5)*{}="x3",
(-12,-3.5)*{}="x4",
(-12,3.5)*{}="x5",
(-18,-3.5)*{}="x6",
(6,66)*{}="y-1",
(0,58)*{}="y0",
(0,66)*{}="y1",
(-6,58)*{}="y2",
(-6,66)*{}="y3",
(-12,58)*{}="y4",
(-12,66)*{}="y5",
(-18,58)*{}="y6",
(6,36.5)*{}="xx1",
(0,30.5)*{}="xx2",
(0,36.5)*{}="xx3",
(-6,30.5)*{}="xx4",
(-6,36.5)*{}="xx5",
(-12,30.5)*{}="xx6",
(-12,36.5)*{}="xx7",
(-18,30.5)*{}="xx8",
\ar@{-} "b1";"t1"^{}
\ar@{-} "b2";"t2"^{}
\ar@{-} "b3";"t3"^{}
\ar@{-} "0.5-l";"13.5-l"^{}
\ar@{-} "r-11.5";"r--1"^{}
\ar@{-} "0-l";"r--1"^{}
%\ar@{-} "1-l";"r"^{}
\ar@{-} "2-l";"r-0"^{}
\ar@{-} "3-l";"r-1"^{}
\ar@{-} "4-l";"r-2"^{}
\ar@{-} "5-l";"r-3"^{}
\ar@{-} "6-l";"r-4"^{}
\ar@{-} "7-l";"r-5"^{}
\ar@{-} "8-l";"r-6"^{}
\ar@{-} "9-l";"r-7"^{}
\ar@{-} "10-l";"r-8"^{}
%\ar@{-} "11-l";"r-9"^{}
\ar@{-} "12-l";"r-10"^{}
\ar@{-} "13-l";"r-11"^{}
\ar@{-} "b1";"r-0"^{}
\ar@{-} "x1";"x2"^{}
\ar@{-} "x3";"x4"^{}
\ar@{-} "x5";"x6"^{}
\ar@{-} "y-1";"y0"^{}
\ar@{-} "y1";"y2"^{}
\ar@{-} "y3";"y4"^{}
\ar@{-} "y5";"y6"^{}
\ar@{-} "xx1";"xx2"^{}
\ar@{-} "xx3";"xx4"^{}
\ar@{-} "xx5";"xx6"^{}
\ar@{-} "xx7";"xx8"^{}
\end{xy}
\]

\[
\begin{xy}
(-24.5,50) *{A^{(2)\dagger}_{2n-2},}="type",
(-24.5,45) *{\Lambda_1:}="type2",
(-15.5,-2) *{\ 1}="000-1",
(-9.5,-2) *{\ 1}="00-1",
(-3.5,-2) *{\ 1}="0-1",
(1.5,-2) *{\ \ 1}="0-1",
(-15.5,1.5) *{\ 1}="0000",
(-9.5,1.5) *{\ 1}="000",
(-3.5,1.5) *{\ 1}="00",
(1.5,1.5) *{\ \ 1}="0",
(1.5,6.5) *{\ \ 2}="02",
(-4,6.5) *{\ \ 2}="002",
(-10,6.5) *{\ \ 2}="0002",
(-16,6.5) *{\ \ 2}="00002",
(1.5,12.5) *{\ \ 3}="03",
(-4,12.5) *{\ \ 3}="003",
(-10,12.5) *{\ \ 3}="0003",
(-16,12.5) *{\ \ 3}="00003",
(1.5,20.5) *{\ \ \vdots}="04",
(-4,20.5) *{\ \ \vdots}="004",
(-10,20.5) *{\ \ \vdots}="0004",
(-16,20.5) *{\ \ \vdots}="00004",
(1.5,27.5) *{\ \ _{n-1}}="05",
(-4,27.5) *{\ \ _{n-1}}="005",
(-10,27.5) *{\ \ _{n-1}}="0005",
(-16,27.5) *{\ \ _{n-1}}="00005",
(1.5,33.5) *{\ \ n}="06",
(-4,33.5) *{\ \ n}="006",
(-10,33.5) *{\ \ n}="0006",
(-16,33.5) *{\ \ n}="00006",
(1.5,39.5) *{\ \ _{n-1}}="07",
(-4,39.5) *{\ \ _{n-1}}="007",
(-10,39.5) *{\ \ _{n-1}}="0007",
(-16,39.5) *{\ \ _{n-1}}="00007",
(1.5,48) *{\ \ \vdots}="08",
(-4,48) *{\ \ \vdots}="008",
(-10,48) *{\ \ \vdots}="0008",
(-16,48) *{\ \ \vdots}="00008",
(1.5,55) *{\ \ 2}="09",
(-4,55) *{\ \ 2}="009",
(-10,55) *{\ \ 2}="0009",
(-16,55) *{\ \ 2}="00009",
(1.5,60) *{\ \ 1}="010",
(-4,60) *{\ \ 1}="0010",
(-10,60) *{\ \ 1}="00010",
(-16,60) *{\ \ 1}="000010",
(1.5,64) *{\ \ 1}="011",
(-4,64) *{\ \ 1}="0011",
(-10,64) *{\ \ 1}="00011",
(-16,64) *{\ \ 1}="000011",
(1.5,69) *{\ \ 2}="012",
(-4,69) *{\ \ 2}="0012",
(-10,69) *{\ \ 2}="00012",
(-16,69) *{\ \ 2}="000012",
(-18,-3.5) *{}="0.5-l",
(-20,-3.5) *{}="0-l",
(-20,0) *{}="1-l",
(-20,3.5) *{}="2-l",
(-20,9.5) *{}="3-l",
(-20,15.5) *{}="4-l",
(-20,24.5) *{}="5-l",
(-20,30.5) *{}="6-l",
(-20,36.5) *{}="7-l",
(-20,42.5) *{}="8-l",
(-20,52) *{}="9-l",
(-20,58) *{}="10-l",
(-20,62) *{}="11-l",
(-20,66) *{}="12-l",
(-20,72) *{}="13-l",
(-18,74) *{}="13.5-l",
(6,-3.5)*{}="r--1",
(6,0)*{}="r",
(6,3.5)*{}="r-0",
(0,-3.5) *{}="b1",
(-6,-3.5)*{}="b2",
(-12,-3.5)*{}="b3",
(0,74) *{}="t1",
(-6,74)*{}="t2",
(-12,74)*{}="t3",
(6,9.5)*{}="r-1",
(6,15.5)*{}="r-2",
(6,24.5)*{}="r-3",
(6,30.5)*{}="r-4",
(6,36.5)*{}="r-5",
(6,42.5)*{}="r-6",
(6,52)*{}="r-7",
(6,58)*{}="r-8",
(6,62)*{}="r-9",
(6,66)*{}="r-10",
(6,72)*{}="r-11",
(6,74)*{}="r-11.5",
\ar@{-} "b1";"t1"^{}
\ar@{-} "b2";"t2"^{}
\ar@{-} "b3";"t3"^{}
\ar@{-} "0.5-l";"13.5-l"^{}
\ar@{-} "r-11.5";"r--1"^{}
\ar@{-} "0-l";"r--1"^{}
\ar@{-} "1-l";"r"^{}
\ar@{-} "2-l";"r-0"^{}
\ar@{-} "3-l";"r-1"^{}
\ar@{-} "4-l";"r-2"^{}
\ar@{-} "5-l";"r-3"^{}
\ar@{-} "6-l";"r-4"^{}
\ar@{-} "7-l";"r-5"^{}
\ar@{-} "8-l";"r-6"^{}
\ar@{-} "9-l";"r-7"^{}
\ar@{-} "10-l";"r-8"^{}
\ar@{-} "11-l";"r-9"^{}
\ar@{-} "12-l";"r-10"^{}
\ar@{-} "13-l";"r-11"^{}
\end{xy}\quad
\begin{xy}
(-24.5,50) *{D^{(2)}_{n},}="type",
(-24.5,45) *{\Lambda_1:}="type2",
(-15.5,-2) *{\ 1}="000-1",
(-9.5,-2) *{\ 1}="00-1",
(-3.5,-2) *{\ 1}="0-1",
(1.5,-2) *{\ \ 1}="0-1",
(-15.5,1.5) *{\ 1}="0000",
(-9.5,1.5) *{\ 1}="000",
(-3.5,1.5) *{\ 1}="00",
(1.5,1.5) *{\ \ 1}="0",
(1.5,6.5) *{\ \ 2}="02",
(-4,6.5) *{\ \ 2}="002",
(-10,6.5) *{\ \ 2}="0002",
(-16,6.5) *{\ \ 2}="00002",
(1.5,12.5) *{\ \ 3}="03",
(-4,12.5) *{\ \ 3}="003",
(-10,12.5) *{\ \ 3}="0003",
(-16,12.5) *{\ \ 3}="00003",
(1.5,20.5) *{\ \ \vdots}="04",
(-4,20.5) *{\ \ \vdots}="004",
(-10,20.5) *{\ \ \vdots}="0004",
(-16,20.5) *{\ \ \vdots}="00004",
(1.5,27.5) *{\ \ _{n-1}}="05",
(-4,27.5) *{\ \ _{n-1}}="005",
(-10,27.5) *{\ \ _{n-1}}="0005",
(-16,27.5) *{\ \ _{n-1}}="00005",
(1.5,32) *{\ \ n}="06",
(-4,32) *{\ \ n}="006",
(-10,32) *{\ \ n}="0006",
(-16,32) *{\ \ n}="00006",
(1.5,35) *{\ \ n}="06-2",
(-4,35) *{\ \ n}="006-2",
(-10,35) *{\ \ n}="0006-2",
(-16,35) *{\ \ n}="00006-2",
(1.5,39.5) *{\ \ _{n-1}}="07",
(-4,39.5) *{\ \ _{n-1}}="007",
(-10,39.5) *{\ \ _{n-1}}="0007",
(-16,39.5) *{\ \ _{n-1}}="00007",
(1.5,48) *{\ \ \vdots}="08",
(-4,48) *{\ \ \vdots}="008",
(-10,48) *{\ \ \vdots}="0008",
(-16,48) *{\ \ \vdots}="00008",
(1.5,55) *{\ \ 2}="09",
(-4,55) *{\ \ 2}="009",
(-10,55) *{\ \ 2}="0009",
(-16,55) *{\ \ 2}="00009",
(1.5,60) *{\ \ 1}="010",
(-4,60) *{\ \ 1}="0010",
(-10,60) *{\ \ 1}="00010",
(-16,60) *{\ \ 1}="000010",
(1.5,64) *{\ \ 1}="011",
(-4,64) *{\ \ 1}="0011",
(-10,64) *{\ \ 1}="00011",
(-16,64) *{\ \ 1}="000011",
(1.5,69) *{\ \ 2}="012",
(-4,69) *{\ \ 2}="0012",
(-10,69) *{\ \ 2}="00012",
(-16,69) *{\ \ 2}="000012",
(-18,-3.5) *{}="0.5-l",
(-20,-3.5) *{}="0-l",
(-20,0) *{}="1-l",
(-20,3.5) *{}="2-l",
(-20,9.5) *{}="3-l",
(-20,15.5) *{}="4-l",
(-20,24.5) *{}="5-l",
(-20,30.5) *{}="6-l",
(-20,33.5) *{}="6-7-l",
(-20,36.5) *{}="7-l",
(-20,42.5) *{}="8-l",
(-20,52) *{}="9-l",
(-20,58) *{}="10-l",
(-20,62) *{}="11-l",
(-20,66) *{}="12-l",
(-20,72) *{}="13-l",
(-18,74) *{}="13.5-l",
(6,-3.5)*{}="r--1",
(6,0)*{}="r",
(6,3.5)*{}="r-0",
(0,-3.5) *{}="b1",
(-6,-3.5)*{}="b2",
(-12,-3.5)*{}="b3",
(0,74) *{}="t1",
(-6,74)*{}="t2",
(-12,74)*{}="t3",
(6,9.5)*{}="r-1",
(6,15.5)*{}="r-2",
(6,24.5)*{}="r-3",
(6,30.5)*{}="r-4",
(6,33.5)*{}="6-7-r",
(6,36.5)*{}="r-5",
(6,42.5)*{}="r-6",
(6,52)*{}="r-7",
(6,58)*{}="r-8",
(6,62)*{}="r-9",
(6,66)*{}="r-10",
(6,72)*{}="r-11",
(6,74)*{}="r-11.5",
\ar@{-} "6-7-l";"6-7-r"^{}
\ar@{-} "b1";"t1"^{}
\ar@{-} "b2";"t2"^{}
\ar@{-} "b3";"t3"^{}
\ar@{-} "0.5-l";"13.5-l"^{}
\ar@{-} "r-11.5";"r--1"^{}
\ar@{-} "0-l";"r--1"^{}
\ar@{-} "1-l";"r"^{}
\ar@{-} "2-l";"r-0"^{}
\ar@{-} "3-l";"r-1"^{}
\ar@{-} "4-l";"r-2"^{}
\ar@{-} "5-l";"r-3"^{}
\ar@{-} "6-l";"r-4"^{}
\ar@{-} "7-l";"r-5"^{}
\ar@{-} "8-l";"r-6"^{}
\ar@{-} "9-l";"r-7"^{}
\ar@{-} "10-l";"r-8"^{}
\ar@{-} "11-l";"r-9"^{}
\ar@{-} "12-l";"r-10"^{}
\ar@{-} "13-l";"r-11"^{}
\end{xy}\quad
\begin{xy}
(-24.5,50) *{D^{(2)}_{n},}="type",
(-24.5,45) *{\Lambda_n:}="type2",
(-15.5,-2) *{\ n}="000-1",
(-9.5,-2) *{\ n}="00-1",
(-3.5,-2) *{\ n}="0-1",
(1.5,-2) *{\ \ n}="0-1",
(-15.5,1.5) *{\ n}="0000",
(-9.5,1.5) *{\ n}="000",
(-3.5,1.5) *{\ n}="00",
(1.5,1.5) *{\ \ n}="0",
(1.5,6.5) *{\ \ _{n-1}}="02",
(-4,6.5) *{\ \ _{n-1}}="002",
(-10,6.5) *{\ \ _{n-1}}="0002",
(-16,6.5) *{\ \ _{n-1}}="00002",
(1.5,12.5) *{\ \ _{n-2}}="03",
(-4,12.5) *{\ \ _{n-2}}="003",
(-10,12.5) *{\ \ _{n-2}}="0003",
(-16,12.5) *{\ \ _{n-2}}="00003",
(1.5,20.5) *{\ \ \vdots}="04",
(-4,20.5) *{\ \ \vdots}="004",
(-10,20.5) *{\ \ \vdots}="0004",
(-16,20.5) *{\ \ \vdots}="00004",
(1.5,27.5) *{\ \ 2}="05",
(-4,27.5) *{\ \ 2}="005",
(-10,27.5) *{\ \ 2}="0005",
(-16,27.5) *{\ \ 2}="00005",
(1.5,32) *{\ \ 1}="06",
(-4,32) *{\ \ 1}="006",
(-10,32) *{\ \ 1}="0006",
(-16,32) *{\ \ 1}="00006",
(1.5,35) *{\ \ 1}="06-2",
(-4,35) *{\ \ 1}="006-2",
(-10,35) *{\ \ 1}="0006-2",
(-16,35) *{\ \ 1}="00006-2",
(1.5,39.5) *{\ \ 2}="07",
(-4,39.5) *{\ \ 2}="007",
(-10,39.5) *{\ \ 2}="0007",
(-16,39.5) *{\ \ 2}="00007",
(1.5,48) *{\ \ \vdots}="08",
(-4,48) *{\ \ \vdots}="008",
(-10,48) *{\ \ \vdots}="0008",
(-16,48) *{\ \ \vdots}="00008",
(1.5,55) *{\ \ _{n-1}}="09",
(-4,55) *{\ \ _{n-1}}="009",
(-10,55) *{\ \ _{n-1}}="0009",
(-16,55) *{\ \ _{n-1}}="00009",
(1.5,60) *{\ \ n}="010",
(-4,60) *{\ \ n}="0010",
(-10,60) *{\ \ n}="00010",
(-16,60) *{\ \ n}="000010",
(1.5,64) *{\ \ n}="011",
(-4,64) *{\ \ n}="0011",
(-10,64) *{\ \ n}="00011",
(-16,64) *{\ \ n}="000011",
(1.5,69) *{\ \ _{n-1}}="012",
(-4,69) *{\ \ _{n-1}}="0012",
(-10,69) *{\ \ _{n-1}}="00012",
(-16,69) *{\ \ _{n-1}}="000012",
(-18,-3.5) *{}="0.5-l",
(-20,-3.5) *{}="0-l",
(-20,0) *{}="1-l",
(-20,3.5) *{}="2-l",
(-20,9.5) *{}="3-l",
(-20,15.5) *{}="4-l",
(-20,24.5) *{}="5-l",
(-20,30.5) *{}="6-l",
(-20,33.5) *{}="6-7-l",
(-20,36.5) *{}="7-l",
(-20,42.5) *{}="8-l",
(-20,52) *{}="9-l",
(-20,58) *{}="10-l",
(-20,62) *{}="11-l",
(-20,66) *{}="12-l",
(-20,72) *{}="13-l",
(-18,74) *{}="13.5-l",
(6,-3.5)*{}="r--1",
(6,0)*{}="r",
(6,3.5)*{}="r-0",
(0,-3.5) *{}="b1",
(-6,-3.5)*{}="b2",
(-12,-3.5)*{}="b3",
(0,74) *{}="t1",
(-6,74)*{}="t2",
(-12,74)*{}="t3",
(6,9.5)*{}="r-1",
(6,15.5)*{}="r-2",
(6,24.5)*{}="r-3",
(6,30.5)*{}="r-4",
(6,33.5)*{}="6-7-r",
(6,36.5)*{}="r-5",
(6,42.5)*{}="r-6",
(6,52)*{}="r-7",
(6,58)*{}="r-8",
(6,62)*{}="r-9",
(6,66)*{}="r-10",
(6,72)*{}="r-11",
(6,74)*{}="r-11.5",
\ar@{-} "6-7-l";"6-7-r"^{}
\ar@{-} "b1";"t1"^{}
\ar@{-} "b2";"t2"^{}
\ar@{-} "b3";"t3"^{}
\ar@{-} "0.5-l";"13.5-l"^{}
\ar@{-} "r-11.5";"r--1"^{}
\ar@{-} "0-l";"r--1"^{}
\ar@{-} "1-l";"r"^{}
\ar@{-} "2-l";"r-0"^{}
\ar@{-} "3-l";"r-1"^{}
\ar@{-} "4-l";"r-2"^{}
\ar@{-} "5-l";"r-3"^{}
\ar@{-} "6-l";"r-4"^{}
\ar@{-} "7-l";"r-5"^{}
\ar@{-} "8-l";"r-6"^{}
\ar@{-} "9-l";"r-7"^{}
\ar@{-} "10-l";"r-8"^{}
\ar@{-} "11-l";"r-9"^{}
\ar@{-} "12-l";"r-10"^{}
\ar@{-} "13-l";"r-11"^{}
\end{xy}
\]
The following is an example of Young wall of ground state $\Lambda_1$ of type ${A}^{(2)\dagger}_{4}$. 
\begin{equation*}
\begin{xy}
(-15.5,-2) *{\ 1}="000-1",
(-9.5,-2) *{\ 1}="00-1",
(-3.5,-2) *{\ 1}="0-1",
(1.5,-2) *{\ \ 1}="0-1",
(-21.5,-2) *{\dots}="00000",
(-9.5,1.5) *{\ 1}="000",
(-3.5,1.5) *{\ 1}="00",
(1.5,1.5) *{\ \ 1}="0",
(1.5,6.5) *{\ \ 2}="02",
(-4,6.5) *{\ \ 2}="002",
(1.5,12.5) *{\ \ 3}="03",
(1.5,18.5) *{\ \ 2}="04",
(0,0) *{}="1",
(0,-3.5) *{}="1-u",
(6,0)*{}="2",
(6,-3.5)*{}="2-u",
(6,3.5)*{}="3",
(0,3.5)*{}="4",
(-6,3.5)*{}="5",
(-6,0)*{}="6",
(-6,-3.5)*{}="6-u",
(-12,3.5)*{}="7",
(-12,0)*{}="8",
(-12,-3.6)*{}="8-u",
(-18,3.5)*{}="9",
(-18,0)*{}="10-a",
(-18,0)*{}="10",
(-18,-3.5)*{}="10-u",
(6,9.5)*{}="3-1",
(6,15.5)*{}="3-2",
(6,21.5)*{}="3-3",
(0,9.5)*{}="4-1",
(0,15.5)*{}="4-2",
(0,21.5)*{}="4-3",
(-6,9.5)*{}="5-1",
\ar@{-} "8";"10-a"^{}
\ar@{-} "10-u";"10-a"^{}
\ar@{-} "10-u";"2-u"^{}
\ar@{-} "8";"8-u"^{}
\ar@{-} "6";"6-u"^{}
\ar@{-} "2";"2-u"^{}
\ar@{-} "1";"1-u"^{}
\ar@{-} "1";"2"^{}
\ar@{-} "1";"4"^{}
\ar@{-} "2";"3"^{}
\ar@{-} "3";"4"^{}
\ar@{-} "5";"6"^{}
\ar@{-} "5";"4"^{}
\ar@{-} "1";"6"^{}
\ar@{-} "7";"8"^{}
\ar@{-} "7";"5"^{}
\ar@{-} "6";"8"^{}
\ar@{-} "3";"3-1"^{}
\ar@{-} "4-1";"3-1"^{}
\ar@{-} "4-1";"4"^{}
\ar@{-} "5-1";"5"^{}
\ar@{-} "4-1";"5-1"^{}
\ar@{-} "4-1";"4-2"^{}
\ar@{-} "3-2";"3-1"^{}
\ar@{-} "3-2";"4-2"^{}
\ar@{-} "3-2";"3-3"^{}
\ar@{-} "4-3";"4-2"^{}
\ar@{-} "4-3";"3-3"^{}
\end{xy}
\end{equation*}

Although, when $X$ is of type $C^{(1)}_{n-1}$, the index $1\in I$ is not in class $1$, we define a pattern as follows:
\[
\begin{xy}
(-24.5,50) *{C^{(1)}_{n-1},}="type",
(-24.5,45) *{\Lambda_1:}="type2",
(-15.5,0) *{\ 1}="000-1",
(-9.5,0) *{\ 1}="00-1",
(-3.5,0) *{\ 1}="0-1",
(1.5,0) *{\ \ 1}="0-1",
%(-15.5,1.5) *{\ 1}="0000",
%(-9.5,1.5) *{\ 1}="000",
%(-3.5,1.5) *{\ 1}="00",
%(1.5,1.5) *{\ \ 1}="0",
(1.5,6.5) *{\ \ 2}="02",
(-4,6.5) *{\ \ 2}="002",
(-10,6.5) *{\ \ 2}="0002",
(-16,6.5) *{\ \ 2}="00002",
(1.5,12.5) *{\ \ 3}="03",
(-4,12.5) *{\ \ 3}="003",
(-10,12.5) *{\ \ 3}="0003",
(-16,12.5) *{\ \ 3}="00003",
(1.5,20.5) *{\ \ \vdots}="04",
(-4,20.5) *{\ \ \vdots}="004",
(-10,20.5) *{\ \ \vdots}="0004",
(-16,20.5) *{\ \ \vdots}="00004",
(1.5,27.5) *{\ \ _{n-1}}="05",
(-4,27.5) *{\ \ _{n-1}}="005",
(-10,27.5) *{\ \ _{n-1}}="0005",
(-16,27.5) *{\ \ _{n-1}}="00005",
(1.5,33.5) *{\ \ n}="06",
(-4,33.5) *{\ \ n}="006",
(-10,33.5) *{\ \ n}="0006",
(-16,33.5) *{\ \ n}="00006",
(1.5,39.5) *{\ \ _{n-1}}="07",
(-4,39.5) *{\ \ _{n-1}}="007",
(-10,39.5) *{\ \ _{n-1}}="0007",
(-16,39.5) *{\ \ _{n-1}}="00007",
(1.5,48) *{\ \ \vdots}="08",
(-4,48) *{\ \ \vdots}="008",
(-10,48) *{\ \ \vdots}="0008",
(-16,48) *{\ \ \vdots}="00008",
(1.5,55) *{\ \ 2}="09",
(-4,55) *{\ \ 2}="009",
(-10,55) *{\ \ 2}="0009",
(-16,55) *{\ \ 2}="00009",
%(1.5,60) *{\ \ 1}="010",
%(-4,60) *{\ \ 1}="0010",
%(-10,60) *{\ \ 1}="00010",
%(-16,60) *{\ \ 1}="000010",
(1.5,62) *{\ \ 1}="011",
(-4,62) *{\ \ 1}="0011",
(-10,62) *{\ \ 1}="00011",
(-16,62) *{\ \ 1}="000011",
(1.5,69) *{\ \ 2}="012",
(-4,69) *{\ \ 2}="0012",
(-10,69) *{\ \ 2}="00012",
(-16,69) *{\ \ 2}="000012",
(-18,-3.5) *{}="0.5-l",
(-20,-3.5) *{}="0-l",
(-20,0) *{}="1-l",
(-20,3.5) *{}="2-l",
(-20,9.5) *{}="3-l",
(-20,15.5) *{}="4-l",
(-20,24.5) *{}="5-l",
(-20,30.5) *{}="6-l",
(-20,36.5) *{}="7-l",
(-20,42.5) *{}="8-l",
(-20,52) *{}="9-l",
(-20,58) *{}="10-l",
(-20,62) *{}="11-l",
(-20,66) *{}="12-l",
(-20,72) *{}="13-l",
(-18,74) *{}="13.5-l",
(6,-3.5)*{}="r--1",
(6,0)*{}="r",
(6,3.5)*{}="r-0",
(0,-3.5) *{}="b1",
(-6,-3.5)*{}="b2",
(-12,-3.5)*{}="b3",
(0,74) *{}="t1",
(-6,74)*{}="t2",
(-12,74)*{}="t3",
(6,9.5)*{}="r-1",
(6,15.5)*{}="r-2",
(6,24.5)*{}="r-3",
(6,30.5)*{}="r-4",
(6,36.5)*{}="r-5",
(6,42.5)*{}="r-6",
(6,52)*{}="r-7",
(6,58)*{}="r-8",
(6,62)*{}="r-9",
(6,66)*{}="r-10",
(6,72)*{}="r-11",
(6,74)*{}="r-11.5",
\ar@{-} "b1";"t1"^{}
\ar@{-} "b2";"t2"^{}
\ar@{-} "b3";"t3"^{}
\ar@{-} "0.5-l";"13.5-l"^{}
\ar@{-} "r-11.5";"r--1"^{}
\ar@{-} "0-l";"r--1"^{}
%\ar@{-} "1-l";"r"^{}
\ar@{-} "2-l";"r-0"^{}
\ar@{-} "3-l";"r-1"^{}
\ar@{-} "4-l";"r-2"^{}
\ar@{-} "5-l";"r-3"^{}
\ar@{-} "6-l";"r-4"^{}
\ar@{-} "7-l";"r-5"^{}
\ar@{-} "8-l";"r-6"^{}
\ar@{-} "9-l";"r-7"^{}
\ar@{-} "10-l";"r-8"^{}
%\ar@{-} "11-l";"r-9"^{}
\ar@{-} "12-l";"r-10"^{}
\ar@{-} "13-l";"r-11"^{}
\end{xy}
\]

Let $k\in I$ be in class $2$.
The pattern obtained by the following procedure is called a supporting pattern (resp. covering pattern) for $\Lambda_k$ of type $X$:
\begin{enumerate}
\item[(1)]
Truncating all boxes of a Young wall pattern for $\Lambda_1$ of type $X$ below first (resp. second) $k$-boxes from the bottom. 
\item[(2)]
Spliting the bottom $k$-boxes into two half-unit height boxes.
\end{enumerate}

The following is an example of procedure to obtain the supporting pattern
for $\Lambda_k$ of type $A^{(2)\dagger}_{2n-2}$:
\[
\begin{xy}
(-24.5,50) *{A^{(2)\dagger}_{2n-2},}="type",
(-24.5,45) *{\Lambda_1:}="type2",
(-15.5,-2) *{\ 1}="000-1",
(-9.5,-2) *{\ 1}="00-1",
(-3.5,-2) *{\ 1}="0-1",
(1.5,-2) *{\ \ 1}="0-1",
(-15.5,1.5) *{\ 1}="0000",
(-9.5,1.5) *{\ 1}="000",
(-3.5,1.5) *{\ 1}="00",
(1.5,1.5) *{\ \ 1}="0",
(1.5,6.5) *{\ \ \vdots}="02",
(-4,6.5) *{\ \ \vdots}="002",
(-10,6.5) *{\ \ \vdots}="0002",
(-16,6.5) *{\ \ \vdots}="00002",
(1.5,12.5) *{\ \ k}="03",
(-4,12.5) *{\ \ k}="003",
(-10,12.5) *{\ \ k}="0003",
(-16,12.5) *{\ \ k}="00003",
(1.5,20.5) *{\ \ \vdots}="04",
(-4,20.5) *{\ \ \vdots}="004",
(-10,20.5) *{\ \ \vdots}="0004",
(-16,20.5) *{\ \ \vdots}="00004",
(1.5,27.5) *{\ \ _{n-1}}="05",
(-4,27.5) *{\ \ _{n-1}}="005",
(-10,27.5) *{\ \ _{n-1}}="0005",
(-16,27.5) *{\ \ _{n-1}}="00005",
(1.5,33.5) *{\ \ n}="06",
(-4,33.5) *{\ \ n}="006",
(-10,33.5) *{\ \ n}="0006",
(-16,33.5) *{\ \ n}="00006",
(1.5,39.5) *{\ \ _{n-1}}="07",
(-4,39.5) *{\ \ _{n-1}}="007",
(-10,39.5) *{\ \ _{n-1}}="0007",
(-16,39.5) *{\ \ _{n-1}}="00007",
(1.5,48) *{\ \ \vdots}="08",
(-4,48) *{\ \ \vdots}="008",
(-10,48) *{\ \ \vdots}="0008",
(-16,48) *{\ \ \vdots}="00008",
(1.5,55) *{\ \ 2}="09",
(-4,55) *{\ \ 2}="009",
(-10,55) *{\ \ 2}="0009",
(-16,55) *{\ \ 2}="00009",
(1.5,60) *{\ \ 1}="010",
(-4,60) *{\ \ 1}="0010",
(-10,60) *{\ \ 1}="00010",
(-16,60) *{\ \ 1}="000010",
(1.5,64) *{\ \ 1}="011",
(-4,64) *{\ \ 1}="0011",
(-10,64) *{\ \ 1}="00011",
(-16,64) *{\ \ 1}="000011",
(1.5,69) *{\ \ 2}="012",
(-4,69) *{\ \ 2}="0012",
(-10,69) *{\ \ 2}="00012",
(-16,69) *{\ \ 2}="000012",
(-18,-3.5) *{}="0.5-l",
(-20,-3.5) *{}="0-l",
(-20,0) *{}="1-l",
(-20,3.5) *{}="2-l",
(-20,9.5) *{}="3-l",
(-20,15.5) *{}="4-l",
(-20,24.5) *{}="5-l",
(-20,30.5) *{}="6-l",
(-20,36.5) *{}="7-l",
(-20,42.5) *{}="8-l",
(-20,52) *{}="9-l",
(-20,58) *{}="10-l",
(-20,62) *{}="11-l",
(-20,66) *{}="12-l",
(-20,72) *{}="13-l",
(-18,74) *{}="13.5-l",
(6,-3.5)*{}="r--1",
(6,0)*{}="r",
(6,3.5)*{}="r-0",
(0,-3.5) *{}="b1",
(-6,-3.5)*{}="b2",
(-12,-3.5)*{}="b3",
(0,74) *{}="t1",
(-6,74)*{}="t2",
(-12,74)*{}="t3",
(6,9.5)*{}="r-1",
(6,15.5)*{}="r-2",
(6,24.5)*{}="r-3",
(6,30.5)*{}="r-4",
(6,36.5)*{}="r-5",
(6,42.5)*{}="r-6",
(6,52)*{}="r-7",
(6,58)*{}="r-8",
(6,62)*{}="r-9",
(6,66)*{}="r-10",
(6,72)*{}="r-11",
(6,74)*{}="r-11.5",
\ar@{-} "b1";"t1"^{}
\ar@{-} "b2";"t2"^{}
\ar@{-} "b3";"t3"^{}
\ar@{-} "0.5-l";"13.5-l"^{}
\ar@{-} "r-11.5";"r--1"^{}
\ar@{-} "0-l";"r--1"^{}
\ar@{-} "1-l";"r"^{}
\ar@{-} "2-l";"r-0"^{}
\ar@{-} "3-l";"r-1"^{}
\ar@{-} "4-l";"r-2"^{}
\ar@{-} "5-l";"r-3"^{}
\ar@{-} "6-l";"r-4"^{}
\ar@{-} "7-l";"r-5"^{}
\ar@{-} "8-l";"r-6"^{}
\ar@{-} "9-l";"r-7"^{}
\ar@{-} "10-l";"r-8"^{}
\ar@{-} "11-l";"r-9"^{}
\ar@{-} "12-l";"r-10"^{}
\ar@{-} "13-l";"r-11"^{}
\end{xy}\quad
\begin{xy}
(-24.5,50) *{\overset{(1)}{\rightarrow}}="type",
(12,50) *{\overset{(2)}{\rightarrow}}="typetype",
(-15.5,0) *{\ k}="000-1",
(-9.5,0) *{\ k}="00-1",
(-3.5,0) *{\ k}="0-1",
(1.5,0) *{\ \ k}="0-1",
%(-15.5,1.5) *{\ k}="0000",
%(-9.5,1.5) *{\ k}="000",
%(-3.5,1.5) *{\ k}="00",
%(1.5,1.5) *{\ \ k}="0",
(1.5,6.5) *{\ \ _{k+1}}="02",
(-4,6.5) *{\ \ _{k+1}}="002",
(-10,6.5) *{\ \ _{k+1}}="0002",
(-16,6.5) *{\ \ _{k+1}}="00002",
(1.5,12.5) *{\ \ _{k+2}}="03",
(-4,12.5) *{\ \ _{k+2}}="003",
(-10,12.5) *{\ \ _{k+2}}="0003",
(-16,12.5) *{\ \ _{k+2}}="00003",
(1.5,20.5) *{\ \ \vdots}="04",
(-4,20.5) *{\ \ \vdots}="004",
(-10,20.5) *{\ \ \vdots}="0004",
(-16,20.5) *{\ \ \vdots}="00004",
(1.5,27.5) *{\ \ _{n-1}}="05",
(-4,27.5) *{\ \ _{n-1}}="005",
(-10,27.5) *{\ \ _{n-1}}="0005",
(-16,27.5) *{\ \ _{n-1}}="00005",
(1.5,33.5) *{\ \ n}="06",
(-4,33.5) *{\ \ n}="006",
(-10,33.5) *{\ \ n}="0006",
(-16,33.5) *{\ \ n}="00006",
(1.5,39.5) *{\ \ _{n-1}}="07",
(-4,39.5) *{\ \ _{n-1}}="007",
(-10,39.5) *{\ \ _{n-1}}="0007",
(-16,39.5) *{\ \ _{n-1}}="00007",
(1.5,48) *{\ \ \vdots}="08",
(-4,48) *{\ \ \vdots}="008",
(-10,48) *{\ \ \vdots}="0008",
(-16,48) *{\ \ \vdots}="00008",
(1.5,55) *{\ \ 2}="09",
(-4,55) *{\ \ 2}="009",
(-10,55) *{\ \ 2}="0009",
(-16,55) *{\ \ 2}="00009",
(1.5,60) *{\ \ 1}="010",
(-4,60) *{\ \ 1}="0010",
(-10,60) *{\ \ 1}="00010",
(-16,60) *{\ \ 1}="000010",
(1.5,64) *{\ \ 1}="011",
(-4,64) *{\ \ 1}="0011",
(-10,64) *{\ \ 1}="00011",
(-16,64) *{\ \ 1}="000011",
(1.5,69) *{\ \ 2}="012",
(-4,69) *{\ \ 2}="0012",
(-10,69) *{\ \ 2}="00012",
(-16,69) *{\ \ 2}="000012",
(-18,-3.5) *{}="0.5-l",
(-20,-3.5) *{}="0-l",
(-20,0) *{}="1-l",
(-20,3.5) *{}="2-l",
(-20,9.5) *{}="3-l",
(-20,15.5) *{}="4-l",
(-20,24.5) *{}="5-l",
(-20,30.5) *{}="6-l",
(-20,36.5) *{}="7-l",
(-20,42.5) *{}="8-l",
(-20,52) *{}="9-l",
(-20,58) *{}="10-l",
(-20,62) *{}="11-l",
(-20,66) *{}="12-l",
(-20,72) *{}="13-l",
(-18,74) *{}="13.5-l",
(6,-3.5)*{}="r--1",
(6,0)*{}="r",
(6,3.5)*{}="r-0",
(0,-3.5) *{}="b1",
(-6,-3.5)*{}="b2",
(-12,-3.5)*{}="b3",
(0,74) *{}="t1",
(-6,74)*{}="t2",
(-12,74)*{}="t3",
(6,9.5)*{}="r-1",
(6,15.5)*{}="r-2",
(6,24.5)*{}="r-3",
(6,30.5)*{}="r-4",
(6,36.5)*{}="r-5",
(6,42.5)*{}="r-6",
(6,52)*{}="r-7",
(6,58)*{}="r-8",
(6,62)*{}="r-9",
(6,66)*{}="r-10",
(6,72)*{}="r-11",
(6,74)*{}="r-11.5",
\ar@{-} "b1";"t1"^{}
\ar@{-} "b2";"t2"^{}
\ar@{-} "b3";"t3"^{}
\ar@{-} "0.5-l";"13.5-l"^{}
\ar@{-} "r-11.5";"r--1"^{}
\ar@{-} "0-l";"r--1"^{}
%\ar@{-} "1-l";"r"^{}
\ar@{-} "2-l";"r-0"^{}
\ar@{-} "3-l";"r-1"^{}
\ar@{-} "4-l";"r-2"^{}
\ar@{-} "5-l";"r-3"^{}
\ar@{-} "6-l";"r-4"^{}
\ar@{-} "7-l";"r-5"^{}
\ar@{-} "8-l";"r-6"^{}
\ar@{-} "9-l";"r-7"^{}
\ar@{-} "10-l";"r-8"^{}
\ar@{-} "11-l";"r-9"^{}
\ar@{-} "12-l";"r-10"^{}
\ar@{-} "13-l";"r-11"^{}
\end{xy}
\begin{xy}
(-24.5,50) *{A^{(2)\dagger}_{2n-2},}="type",
(-24.5,45) *{\Lambda_k:}="type2",
(-15.5,-2) *{\ k}="000-1",
(-9.5,-2) *{\ k}="00-1",
(-3.5,-2) *{\ k}="0-1",
(1.5,-2) *{\ \ k}="0-1",
(-15.5,1.5) *{\ k}="0000",
(-9.5,1.5) *{\ k}="000",
(-3.5,1.5) *{\ k}="00",
(1.5,1.5) *{\ \ k}="0",
(1.5,6.5) *{\ \ _{k+1}}="02",
(-4,6.5) *{\ \ _{k+1}}="002",
(-10,6.5) *{\ \ _{k+1}}="0002",
(-16,6.5) *{\ \ _{k+1}}="00002",
(1.5,12.5) *{\ \ _{k+2}}="03",
(-4,12.5) *{\ \ _{k+2}}="003",
(-10,12.5) *{\ \ _{k+2}}="0003",
(-16,12.5) *{\ \ _{k+2}}="00003",
(1.5,20.5) *{\ \ \vdots}="04",
(-4,20.5) *{\ \ \vdots}="004",
(-10,20.5) *{\ \ \vdots}="0004",
(-16,20.5) *{\ \ \vdots}="00004",
(1.5,27.5) *{\ \ _{n-1}}="05",
(-4,27.5) *{\ \ _{n-1}}="005",
(-10,27.5) *{\ \ _{n-1}}="0005",
(-16,27.5) *{\ \ _{n-1}}="00005",
(1.5,33.5) *{\ \ n}="06",
(-4,33.5) *{\ \ n}="006",
(-10,33.5) *{\ \ n}="0006",
(-16,33.5) *{\ \ n}="00006",
(1.5,39.5) *{\ \ _{n-1}}="07",
(-4,39.5) *{\ \ _{n-1}}="007",
(-10,39.5) *{\ \ _{n-1}}="0007",
(-16,39.5) *{\ \ _{n-1}}="00007",
(1.5,48) *{\ \ \vdots}="08",
(-4,48) *{\ \ \vdots}="008",
(-10,48) *{\ \ \vdots}="0008",
(-16,48) *{\ \ \vdots}="00008",
(1.5,55) *{\ \ 2}="09",
(-4,55) *{\ \ 2}="009",
(-10,55) *{\ \ 2}="0009",
(-16,55) *{\ \ 2}="00009",
(1.5,60) *{\ \ 1}="010",
(-4,60) *{\ \ 1}="0010",
(-10,60) *{\ \ 1}="00010",
(-16,60) *{\ \ 1}="000010",
(1.5,64) *{\ \ 1}="011",
(-4,64) *{\ \ 1}="0011",
(-10,64) *{\ \ 1}="00011",
(-16,64) *{\ \ 1}="000011",
(1.5,69) *{\ \ 2}="012",
(-4,69) *{\ \ 2}="0012",
(-10,69) *{\ \ 2}="00012",
(-16,69) *{\ \ 2}="000012",
(-18,-3.5) *{}="0.5-l",
(-20,-3.5) *{}="0-l",
(-20,0) *{}="1-l",
(-20,3.5) *{}="2-l",
(-20,9.5) *{}="3-l",
(-20,15.5) *{}="4-l",
(-20,24.5) *{}="5-l",
(-20,30.5) *{}="6-l",
(-20,36.5) *{}="7-l",
(-20,42.5) *{}="8-l",
(-20,52) *{}="9-l",
(-20,58) *{}="10-l",
(-20,62) *{}="11-l",
(-20,66) *{}="12-l",
(-20,72) *{}="13-l",
(-18,74) *{}="13.5-l",
(6,-3.5)*{}="r--1",
(6,0)*{}="r",
(6,3.5)*{}="r-0",
(0,-3.5) *{}="b1",
(-6,-3.5)*{}="b2",
(-12,-3.5)*{}="b3",
(0,74) *{}="t1",
(-6,74)*{}="t2",
(-12,74)*{}="t3",
(6,9.5)*{}="r-1",
(6,15.5)*{}="r-2",
(6,24.5)*{}="r-3",
(6,30.5)*{}="r-4",
(6,36.5)*{}="r-5",
(6,42.5)*{}="r-6",
(6,52)*{}="r-7",
(6,58)*{}="r-8",
(6,62)*{}="r-9",
(6,66)*{}="r-10",
(6,72)*{}="r-11",
(6,74)*{}="r-11.5",
\ar@{-} "b1";"t1"^{}
\ar@{-} "b2";"t2"^{}
\ar@{-} "b3";"t3"^{}
\ar@{-} "0.5-l";"13.5-l"^{}
\ar@{-} "r-11.5";"r--1"^{}
\ar@{-} "0-l";"r--1"^{}
\ar@{-} "1-l";"r"^{}
\ar@{-} "2-l";"r-0"^{}
\ar@{-} "3-l";"r-1"^{}
\ar@{-} "4-l";"r-2"^{}
\ar@{-} "5-l";"r-3"^{}
\ar@{-} "6-l";"r-4"^{}
\ar@{-} "7-l";"r-5"^{}
\ar@{-} "8-l";"r-6"^{}
\ar@{-} "9-l";"r-7"^{}
\ar@{-} "10-l";"r-8"^{}
\ar@{-} "11-l";"r-9"^{}
\ar@{-} "12-l";"r-10"^{}
\ar@{-} "13-l";"r-11"^{}
\end{xy}
\]

\begin{defn}\cite{Ka24a}
Let $k\in I$ be in class $2$.
A wall $Y$ is called a {\it truncated wall} of supporting (resp. covering) ground state $\Lambda_k$ of type $X$ if it satisfies the following:
\begin{enumerate}
\item The wall $Y$ is obtained from $Y_{\Lambda_k}$
by stacking finitely many colored blocks on the ground-state wall $Y_{\Lambda_k}$.
\item There are no blocks on the top of a single block with half-unit thickness.
\item Except for the right-most column,
there are no free spaces to the
right of all blocks.
\item The colored blocks are stacked following the supporting (resp. covering) pattern for $\Lambda_k$ of type $X$.
\end{enumerate}
\end{defn}

Here, the ground state wall $Y_{\Lambda_k}$ is as follows:
\[
\begin{xy}
(-15.5,1.5) *{\ \cdots}="0000",
(-9.5,1.5) *{\ k}="000",
(-3.5,1.5) *{\ k}="00",
(1.5,1.5) *{\ \ k}="0",
(0,0) *{}="1",
(6,0)*{}="2",
(6,3.5)*{}="3",
(0,3.5)*{}="4",
(-6,3.5)*{}="5",
(-6,0)*{}="6",
(-12,3.5)*{}="7",
(-12,0)*{}="8",
\ar@{-} "1";"2"^{}
\ar@{-} "1";"4"^{}
\ar@{-} "2";"3"^{}
\ar@{-} "3";"4"^{}
\ar@{-} "5";"6"^{}
\ar@{-} "5";"4"^{}
\ar@{-} "1";"6"^{}
\ar@{-} "7";"8"^{}
\ar@{-} "7";"5"^{}
\ar@{-} "6";"8"^{}
\end{xy}
\]
The conditions (i), (ii) and (iii) are same as in the Definition \ref{def-YW1}.

\begin{defn}\cite{Kang}\label{def-YW2}
Let $Y$ be a Young wall or truncated wall of ground state $\lambda$ of type $X$.
\begin{enumerate}
\item A column of $Y$ is called a {\it full column} if its top has unit thickness
and its height is a multiple of the unit length.
\item For $X=A^{(1)}_{n-1}$, we understand all Young walls are proper.
\item For $X=A^{(2)}_{2n-3}$, $A^{(2)\dagger}_{2n-2}$, $B^{(1)}_{n-1}$, $D^{(1)}_{n-1}$ and $D^{(2)}_n$,
the Young wall $Y$ is said to be {\it proper} if none of two full columns of $Y$ have the same height.
\end{enumerate}
\end{defn}

\begin{defn}
Let $k\in I$ be in class $2$ and $Y$ be a proper truncated wall of ground state $\Lambda_k$.
A $k$-block with half-unit height put on the ground state wall is called a {\it second half} $k$-{\it block}. 
\end{defn}

\begin{defn}\cite{Kang}\label{def-YW2a}
Let $Y$ be a proper Young wall or proper truncated wall of ground state $\Lambda_k$ and $i\in I$.
\begin{enumerate}
\item An $i$-block in $Y$ other than second half $k$-block
is said to be a {\it removable} $i$-{\it block} 
if the wall obtained from $Y$ by removing this $i$-block 
remains proper.
\item
A place is said to be an $i$-{\it admissible} {\it slot}
if the wall obtained from $Y$ by adding
an $i$-block other than second half $k$-block to this place
remains proper.
\end{enumerate}
\end{defn}

\begin{defn}\cite{Ka23a}\label{def-YW2b}
Let $Y$ be a proper Young wall or proper truncated wall and $t=1$ or $t=n$. 
\begin{enumerate}
\item[(i)]
Let $Y'$ be a wall
obtained from $Y$
by putting two $t$-blocks of shape (\ref{smpl2}) on the top of a column of $Y$:
\[
Y=
\begin{xy}
(14,10) *{\leftarrow A}="A",
(-3,3) *{\cdots}="dot1",
(12,3) *{\cdots}="dot2",
(0,-4) *{}="-1",
(8,-4)*{}="-2",
(4,-6)*{\vdots}="dot3",
(0,0) *{}="1",
(8,0)*{}="2",
(8,8)*{}="3",
(0,8)*{}="4",
(8,12)*{}="5",
(0,12)*{}="6",
(8,16)*{}="7",
(0,16)*{}="8",
\ar@{-} "1";"-1"^{}
\ar@{-} "-2";"2"^{}
\ar@{-} "1";"2"^{}
\ar@{-} "1";"4"^{}
\ar@{-} "2";"3"^{}
\ar@{-} "3";"4"^{}
\ar@{--} "3";"5"^{}
\ar@{--} "4";"6"^{}
\ar@{--} "6";"5"^{}
\ar@{--} "7";"5"^{}
\ar@{--} "6";"8"^{}
\ar@{--} "7";"8"^{}
\end{xy}\qquad
Y'=
\begin{xy}
(14,14) *{}="B",
(4,10)*{t}="t1",
(4,14)*{t}="t2",
(-3,3) *{\cdots}="dot1",
(12,3) *{\cdots}="dot2",
(0,-4) *{}="-1",
(8,-4)*{}="-2",
(4,-6)*{\vdots}="dot3",
(0,0) *{}="1",
(8,0)*{}="2",
(8,8)*{}="3",
(0,8)*{}="4",
(8,12)*{}="5",
(0,12)*{}="6",
(8,16)*{}="7",
(0,16)*{}="8",
\ar@{-} "1";"-1"^{}
\ar@{-} "-2";"2"^{}
\ar@{-} "1";"2"^{}
\ar@{-} "1";"4"^{}
\ar@{-} "2";"3"^{}
\ar@{-} "3";"4"^{}
\ar@{-} "3";"5"^{}
\ar@{-} "4";"6"^{}
\ar@{-} "6";"5"^{}
\ar@{-} "7";"5"^{}
\ar@{-} "6";"8"^{}
\ar@{-} "7";"8"^{}
\end{xy}
\]
In $Y$, we give a name to a slot as above.
The slot $A$ in $Y$ is said to be {\it double} $t$-{\it admissible}
if $Y'$ is proper.
\item[(ii)]
We assume that we obtain a wall $Y''$ by removing two $t$-blocks of shape (\ref{smpl2}) from the top of a column of $Y$:
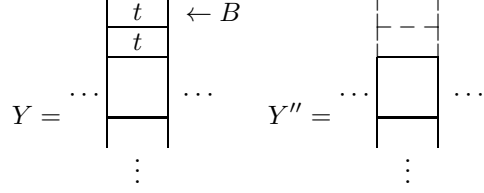
\[
Y=
\begin{xy}
(14,14) *{\leftarrow B}="B",
(4,10)*{t}="t1",
(4,14)*{t}="t2",
(-3,3) *{\cdots}="dot1",
(12,3) *{\cdots}="dot2",
(0,-4) *{}="-1",
(8,-4)*{}="-2",
(4,-6)*{\vdots}="dot3",
(0,0) *{}="1",
(8,0)*{}="2",
(8,8)*{}="3",
(0,8)*{}="4",
(8,12)*{}="5",
(0,12)*{}="6",
(8,16)*{}="7",
(0,16)*{}="8",
\ar@{-} "1";"-1"^{}
\ar@{-} "-2";"2"^{}
\ar@{-} "1";"2"^{}
\ar@{-} "1";"4"^{}
\ar@{-} "2";"3"^{}
\ar@{-} "3";"4"^{}
\ar@{-} "3";"5"^{}
\ar@{-} "4";"6"^{}
\ar@{-} "6";"5"^{}
\ar@{-} "7";"5"^{}
\ar@{-} "6";"8"^{}
\ar@{-} "7";"8"^{}
\end{xy}\quad
Y''=
\begin{xy}
(14,10) *{}="A",
(-3,3) *{\cdots}="dot1",
(12,3) *{\cdots}="dot2",
(0,-4) *{}="-1",
(8,-4)*{}="-2",
(4,-6)*{\vdots}="dot3",
(0,0) *{}="1",
(8,0)*{}="2",
(8,8)*{}="3",
(0,8)*{}="4",
(8,12)*{}="5",
(0,12)*{}="6",
(8,16)*{}="7",
(0,16)*{}="8",
\ar@{-} "1";"-1"^{}
\ar@{-} "-2";"2"^{}
\ar@{-} "1";"2"^{}
\ar@{-} "1";"4"^{}
\ar@{-} "2";"3"^{}
\ar@{-} "3";"4"^{}
\ar@{--} "3";"5"^{}
\ar@{--} "4";"6"^{}
\ar@{--} "6";"5"^{}
\ar@{--} "7";"5"^{}
\ar@{--} "6";"8"^{}
\ar@{--} "7";"8"^{}
\end{xy}
\]
A block in $Y$ is named $B$ as above.
The block $B$ in $Y$ is said to be {\it double} $t$-{\it removable}
if $Y''$ is proper.
\item[(iii)] For $i\in I$, other $i$-admissible slots (resp. removable $i$-blocks)
are said to be single admissible (resp. single removable).
\end{enumerate}
\end{defn}

\subsection{Assignment of homomorphisms to walls}

Following \cite{Ka24a}, we assign some 
homomorphisms 
to Young walls and truncated walls.
As for the things we will not use in this paper, we just cited from the previous paper.

For $k\in I$, we define $T_k^X$, $\overline{T}_k^X$ and $\overline{\overline{T}}_k^X\in\mathbb{Z}_{\geq1}$ just as in subsection 3.1 of \cite{Ka24a}. From now on, we will
draw proper Young walls of ground state $\Lambda_k$ of type $X$ in $\mathbb{R}_{\leq0}\times \mathbb{R}_{\geq T^X_k}$
for each $k\in I$ in class $1$.
We draw truncated walls of supporting (resp. covering) ground state $\Lambda_k$ of type $X$ in $\mathbb{R}_{\leq0}\times \mathbb{R}_{\geq \overline{T}^X_k}$
(resp. in $\mathbb{R}_{\leq0}\times \mathbb{R}_{\geq \overline{\overline{T}}^X_k}$) for each $k\in I$ in class $2$. 

Let $S$ be one of the following slots or blocks in a proper Young wall or truncated wall 
with $i\in\mathbb{Z}_{\geq0}$, $l\in\mathbb{Z}_{\geq1}$: 
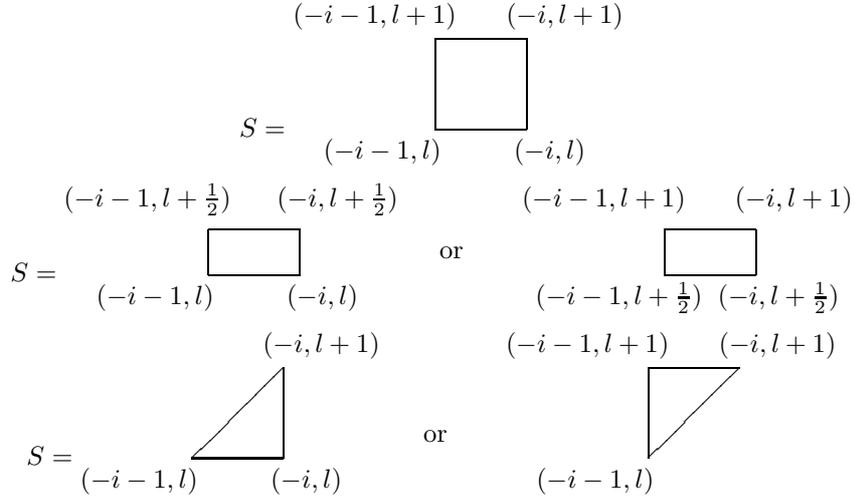
\[
S=
\begin{xy}
(-8,15) *{(-i-1,l+1)}="0000",
(17,15) *{(-i,l+1)}="000",
(15,-3) *{(-i,l)}="00",
(-7,-3) *{(-i-1,l)}="0",
(0,0) *{}="1",
(12,0)*{}="2",
(12,12)*{}="3",
(0,12)*{}="4",
\ar@{-} "1";"2"^{}
\ar@{-} "1";"4"^{}
\ar@{-} "2";"3"^{}
\ar@{-} "3";"4"^{}
\end{xy}
\]
\[
S=
\begin{xy}
(-8,10) *{(-i-1,l+\frac{1}{2})}="0000",
(17,10) *{(-i,l+\frac{1}{2})}="000",
(15,-3) *{(-i,l)}="00",
(-7,-3) *{(-i-1,l)}="0",
(0,0) *{}="1",
(12,0)*{}="2",
(12,6)*{}="3",
(0,6)*{}="4",
(32,3) *{{\rm or}}="or",
(52,10) *{(-i-1,l+1)}="a0000",
(77,10) *{(-i,l+1)}="a000",
(75,-3) *{(-i,l+\frac{1}{2})}="a00",
(54,-3) *{(-i-1,l+\frac{1}{2})}="a0",
(60,0) *{}="a1",
(72,0)*{}="a2",
(72,6)*{}="a3",
(60,6)*{}="a4",
\ar@{-} "1";"2"^{}
\ar@{-} "1";"4"^{}
\ar@{-} "2";"3"^{}
\ar@{-} "3";"4"^{}
\ar@{-} "a1";"a2"^{}
\ar@{-} "a1";"a4"^{}
\ar@{-} "a2";"a3"^{}
\ar@{-} "a3";"a4"^{}
\end{xy}
\]
\[
S=
\begin{xy}
(17,15) *{(-i,l+1)}="000",
(15,-3) *{(-i,l)}="00",
(-7,-3) *{(-i-1,l)}="0",
(0,0) *{}="1",
(12,0)*{}="2",
(12,12)*{}="3",
(32,3) *{{\rm or}}="or",
(52,15) *{(-i-1,l+1)}="0000a",
(77,15) *{(-i,l+1)}="000a",
(53,-3) *{(-i-1,l)}="0a",
(60,0) *{}="1a",
(72,12)*{}="3a",
(60,12)*{}="4a",
\ar@{-} "1";"2"^{}
\ar@{-} "1";"3"^{}
\ar@{-} "2";"3"^{}
\ar@{-} "1a";"3a"^{}
\ar@{-} "1a";"4a"^{}
\ar@{-} "3a";"4a"^{}
\end{xy}
\]
Then we say the position of $S$ is $(-i,l)$.
If $S$ is colored by $t\in I$ then
we assign a homomorphism
\begin{equation}\label{l1def1}
L^X_{s,k,{\rm ad}}(S):=
\begin{cases}
x_{r,t} & \text{if }S\text{ is admissible,}
\\
0 & \text{otherwise,}
\end{cases}\quad
L^X_{s,k,{\rm re}}(S):=
\begin{cases}
x_{r+1,t} & \text{if }S\text{ is removable,}
\\
0 & \text{otherwise.}
\end{cases}
\end{equation}
Here, $r\in\mathbb{Z}_{\geq s}$ is determined from $i,l,s$ and $\iota$. The precise definition of $r$ is
in subsection 4.3 of \cite{Ka24a}.

Let $Y$ be a proper Young wall or proper truncated wall of ground state
$\Lambda_k$ of type $X=A^{(1)}_{n-1}$, $B^{(1)}_{n-1}$, $C^{(1)}_{n-1}$, $D^{(1)}_{n-1}$, $A^{(2)\dagger}_{2n-2}$, $A^{(2)}_{2n-3}$ or
$D^{(2)}_{n}$. For $s\in \mathbb{Z}$, we define
\begin{eqnarray}
L^X_{s,k,\iota}(Y)&:=&
\sum_{t\in I} \left(
\sum_{P:\text{single }t\text{-admissible slot}} L^X_{s,k,{\rm ad}}(P)
-\sum_{P:\text{single removable }t\text{-block}} L^X_{s,k,{\rm re}}(P)\right)\nonumber\\
& & + 
\sum_{P:\text{double }t\text{-admissible slot}} 2L^X_{s,k,{\rm ad}}(P)
-\sum_{P:\text{double removable }t\text{-block}} 2L^X_{s,k,{\rm re}}(P). \label{L11-def}
\end{eqnarray}

\subsection{Set of walls}

\begin{defn}\cite{Ka24a}\label{YWk-def}
\begin{enumerate}
\item
For $k\in I$ in class $1$,
we define ${\rm YW}^{X}_{k}$ as
the set of all proper Young walls of ground state $\Lambda_k$ of type $X$.
\item
For $k\in I$ in class $2$,
one defines $\overline{{\rm YW}^{X}_k}$ 
(resp. $\overline{\overline{{\rm YW}^{X}_k}}$)
as
the set of all proper truncated walls of supporting (resp. covering) ground state $\Lambda_k$ of type $X$.
\item
Let $k\in I$ be in class $2$. We define
\[
{\rm YW}^{X}_{k}:=\left\{(\overline{Y},\overline{\overline{Y}})\in \overline{{\rm YW}^{X}_k}\times \overline{\overline{{\rm YW}^{X}_k}} |
h_j=\frac{1}{2} \text{ if and only if }h_j'=\frac{1}{2}\ \text{for }j\in\mathbb{Z}_{\geq1}
\right\},
\]
where $h_j$ (resp. $h_j'$) is the height of $j$-th column of $\overline{Y}$ (resp. $\overline{\overline{Y}}$) from the right.
We also define
\[
Y_{\Lambda_k}:=(
\begin{xy}
(-15.5,1.5) *{\ \cdots}="0000",
(-9.5,1.5) *{\ k}="000",
(-3.5,1.5) *{\ k}="00",
(1.5,1.5) *{\ \ k}="0",
(0,0) *{}="1",
(6,0)*{}="2",
(6,3.5)*{}="3",
(0,3.5)*{}="4",
(-6,3.5)*{}="5",
(-6,0)*{}="6",
(-12,3.5)*{}="7",
(-12,0)*{}="8",
\ar@{-} "1";"2"^{}
\ar@{-} "1";"4"^{}
\ar@{-} "2";"3"^{}
\ar@{-} "3";"4"^{}
\ar@{-} "5";"6"^{}
\ar@{-} "5";"4"^{}
\ar@{-} "1";"6"^{}
\ar@{-} "7";"8"^{}
\ar@{-} "7";"5"^{}
\ar@{-} "6";"8"^{}
\end{xy},
\begin{xy}
(-15.5,1.5) *{\ \cdots}="0000",
(-9.5,1.5) *{\ k}="000",
(-3.5,1.5) *{\ k}="00",
(1.5,1.5) *{\ \ k}="0",
(0,0) *{}="1",
(6,0)*{}="2",
(6,3.5)*{}="3",
(0,3.5)*{}="4",
(-6,3.5)*{}="5",
(-6,0)*{}="6",
(-12,3.5)*{}="7",
(-12,0)*{}="8",
\ar@{-} "1";"2"^{}
\ar@{-} "1";"4"^{}
\ar@{-} "2";"3"^{}
\ar@{-} "3";"4"^{}
\ar@{-} "5";"6"^{}
\ar@{-} "5";"4"^{}
\ar@{-} "1";"6"^{}
\ar@{-} "7";"8"^{}
\ar@{-} "7";"5"^{}
\ar@{-} "6";"8"^{}
\end{xy})\in {\rm YW}^X_k.
\]
\end{enumerate}
\end{defn}

\begin{defn}\cite{Ka24a}
Let $k\in I$ be in class $2$.
\begin{enumerate}
\item
For a pair $Y=(\overline{Y},\overline{\overline{Y}})\in {\rm YW}^{X}_{k}$,
we assume that
$Y'=(\overline{Y}',\overline{\overline{Y}}')\in {\rm YW}^{X}_{k}$ is obtained from $Y$
by removing one second half $k$-block
\[
\begin{xy}
(1.5,1.5) *{\ \ k}="0",
(0,0) *{}="1",
(6,0)*{}="2",
(6,3.5)*{}="3",
(0,3.5)*{}="4",
\ar@{-} "1";"2"^{}
\ar@{-} "1";"4"^{}
\ar@{-} "2";"3"^{}
\ar@{-} "3";"4"^{}
\end{xy}
\]
from the top of $j$-th column of $\overline{Y}$ and one from
the top of $j$-th column of $\overline{\overline{Y}}$ with some $j\in \mathbb{Z}_{\geq1}$.
Then the pair of these two $k$-blocks are said to be a {\it removable} $k$-{\it pair}.
\item
For $Y=(\overline{Y},\overline{\overline{Y}})\in {\rm YW}^{X}_{k}$,
we assume that
$Y''=(\overline{Y}'',\overline{\overline{Y}}'')\in {\rm YW}^{X}_{k}$ is obtained from $Y$
by adding one second half $k$-block
\[
\begin{xy}
(1.5,1.5) *{\ \ k}="0",
(0,0) *{}="1",
(6,0)*{}="2",
(6,3.5)*{}="3",
(0,3.5)*{}="4",
\ar@{-} "1";"2"^{}
\ar@{-} "1";"4"^{}
\ar@{-} "2";"3"^{}
\ar@{-} "3";"4"^{}
\end{xy}
\]
to the top of $j$-th column of $\overline{Y}$ and one to 
the top of $j$-th column of
$\overline{\overline{Y}}$ with some $j\in \mathbb{Z}_{\geq1}$. 
Then the pair of these two $k$-slots are said to be a $k$-{\it admissible pair}.
\end{enumerate}
\end{defn}

For $k\in I$ in class $2$,
let $P$ be a pair of $k$-blocks in some $Y\in {\rm YW}^{X}_{k}$ with the following coordinate:
\[
P=
\begin{xy}
(-8,10) *{(-i-1,\overline{T}_k+1)}="0000",
(17,10) *{(-i,\overline{T}_k+1)}="000",
(15,-3) *{(-i,\overline{T}_k+\frac{1}{2})}="00",
(-7,-3) *{(-i-1,\overline{T}_k+\frac{1}{2})}="0",
(6,3) *{k}="0",
(0,0) *{}="1",
(12,0)*{}="2",
(12,6)*{}="3",
(0,6)*{}="4",
(32,3) *{}="or",
(52,10) *{(-i-1,\overline{\overline{T}}_k+1)}="a0000",
(80,10) *{(-i,\overline{\overline{T}}_k+1)}="a000",
(80,-3) *{(-i,\overline{\overline{T}}_k+\frac{1}{2})}="a00",
(54,-3) *{(-i-1,\overline{\overline{T}}_k+\frac{1}{2})}="a0",
(60,0) *{}="a1",
(72,0)*{}="a2",
(72,6)*{}="a3",
(60,6)*{}="a4",
(66,3) *{k}="0",
\ar@{-} "1";"2"^{}
\ar@{-} "1";"4"^{}
\ar@{-} "2";"3"^{}
\ar@{-} "3";"4"^{}
\ar@{-} "a1";"a2"^{}
\ar@{-} "a1";"a4"^{}
\ar@{-} "a2";"a3"^{}
\ar@{-} "a3";"a4"^{}
\end{xy}
\]
Then we define
\begin{equation}\label{l1def5}
L^X_{s,k,{\rm ad}}(P):=
\begin{cases}
x_{s+i,k} & \text{if }P\text{ is a }k\text{-admissible pair}, \\
0 & \text{otherwise,}
\end{cases}
\quad L^X_{s,k,{\rm re}}(P):=
\begin{cases}
x_{s+i+1,k} & \text{if }P\text{ is a removable }k\text{-pair}, \\
0 & \text{otherwise.}
\end{cases}
\end{equation}

For above pair  $P$,
we say the {\it position} of $P$ is $(-i,\overline{T}_k)$.
For $k\in I$ in class $2$ and
$Y=(\overline{Y},\overline{\overline{Y}})\in {\rm YW}^{X}_{k}$, 
using (\ref{L11-def}) and (\ref{l1def5}),
one
defines
\begin{eqnarray}
L^X_{s,k,\iota}(Y)&:=& L^X_{s,k,\iota}(\overline{Y})+L^X_{s,k,\iota}(\overline{\overline{Y}})\nonumber\\
& &+
\sum_{P: k\text{-admissible pair in }Y} L^X_{s,k,{\rm ad}}(P)
-\sum_{P: \text{removable }k\text{-pair in }Y} L^X_{s,k,{\rm re}}(P). \label{L-pair-def}
\end{eqnarray}
For each $k\in I$ and $s\in\mathbb{Z}_{\geq1}$, we set
${\rm COMB}^{X}_{s,k,\iota}[\infty]
:=
\{
L^X_{s,k,\iota}(Y) | Y\in {\rm YW}^X_k
\}$.

\begin{thm}\cite{Ka24a}\label{pre-thm1}
Let $\mathfrak{g}$ be of type $X^L$ ($=A^{(1)}_{n-1}$, $B^{(1)}_{n-1}$, $C^{(1)}_{n-1}$, $D^{(1)}_{n-1}$, $A^{(2)}_{2n-2}$, $A^{(2)}_{2n-3}$ or
$D^{(2)}_{n}$), where $X^L$ is defined in (\ref{Langtype}).
Let $\iota$ be an adapted sequence. Then $\iota$ satisfies the $\Xi'$-positivity condition and we have
\[
\Xi_{s,k,\iota}'={\rm COMB}^{X}_{s,k,\iota}[\infty].
\]
\end{thm}

\begin{prop}\cite{Ka24a}\label{pre-prop1}
Let $\mathfrak{g}$ be of type $X^L$ ($=A^{(1)}_{n-1}$, $B^{(1)}_{n-1}$, $C^{(1)}_{n-1}$, $D^{(1)}_{n-1}$, $A^{(2)}_{2n-2}$, $A^{(2)}_{2n-3}$ or
$D^{(2)}_{n}$) and $s\in\mathbb{Z}_{\geq1}$, $k\in I$.
\begin{enumerate}
\item[(1)] For $Y_{\Lambda_k}\in {\rm YW}^X_k$,
we have $L^X_{s,k,\iota}(Y_{\Lambda_k})=x_{s,k}$.
\item[(2)] Let $Y,Y'\in {\rm YW}^X_k$ and we assume that $Y'$ is obtained from $Y$ by adding a $t$-block or pair of $t$-blocks
to a $t$-admissible slot or $t$-admissible pair $P$ for $t\in I$. Then we have
\[
L^{X}_{s,k,\iota}(Y')=L^{X}_{s,k,\iota}(Y)-\beta_{r,t},
\]
where $r\in\mathbb{Z}_{\geq s}$ and
$L^X_{s,k,{\rm ad}}(P)=x_{r,t}$.
\end{enumerate}
\end{prop}

\subsection{Proof of Theorem \ref{mainthm2} (3)}

We assume $\mathfrak{g}$ is of type $X^L$.
Let us show Theorem \ref{mainthm2} (3).
By Theorem \ref{mainthm1} and \ref{pre-thm1}, we need to show ${\rm COMB}^{X}_{s,k,\iota}[\infty]\subset \Xi_{s,k,\iota}'^+$.
We take any $Y\in {\rm YW}^X_k$ and let us prove $L^X_{s,k,\iota}(Y)\in \Xi_{s,k,\iota}'^+$ by the induction on the number of blocks. 
Here, 
we say the number of blocks in $Y$ is $m$ if we get the wall $Y$ by adding $m$ blocks to $Y_{\Lambda_k}$.

If the number of blocks is $0$ then we have $Y=Y_{\Lambda_k}$ and $L^X_{s,k,\iota}(Y)=x_{s,k}\in \Xi_{s,k,\iota}'^+$ (Proposition \ref{pre-prop1} (1))
so that we may assume $Y$ has at least one block.
It follows from Proposition \ref{pre-prop1} (2) that
$L^X_{s,k,\iota}(Y)$ can be expressed as
\[
L^X_{s,k,\iota}(Y)=x_{s,k}-\sum_{(r,j)\in\mathbb{Z}_{\geq s}\times I} c_{r,j}\beta_{r,j}
\]
with non-negative integers $\{c_{r,j}\}$. It holds $c_{r,j}=0$ except for finitely many $(r,j)$.
By $Y\neq Y_{\Lambda_k}$, there is a pair $(m',t')\in\mathbb{Z}_{\geq s}\times I$ such as
$(m',t')
={\rm max}\{(r,j)\in\mathbb{Z}_{\geq s}\times I |c_{r,j}>0 \}$. Recall that we defined the order on $\mathbb{Z}_{\geq 1}\times I$ in the subsection \ref{s-d}.

Taking (\ref{betask}) into the account, it is easy to see
the coefficient of $x_{m'+1,t'}$ is negative in $L^X_{s,k,\iota}(Y)$.
Thus,
one can also take $(m'',t'')$ as
\[
(m'',t'')={\rm min}\{(r,j)\in\mathbb{Z}_{\geq 1}\times I | \text{the coefficient of }x_{r,j}\text{ in }L^X_{s,k,\iota}(Y)\text{ is negative}\}.
\]
The definitions (\ref{L11-def}), (\ref{L-pair-def}) of $L^X_{s,k,\iota}(Y)$ yield that
there is a removable $t''$-block or removable $t''$-pair $B$ in $Y$ such that $L^X_{s,k,{\rm re}}(B)=x_{m'',t''}$ and $m''\geq2$. 
Let $Y''\in {\rm YW}^X_k$
be the wall obtained from $Y$ by removing $B$.
Using the induction assumption, it follows $L^X_{s,k,\iota}(Y'')\in \Xi_{s,k,\iota}'^+$.
Proposition \ref{pre-prop1} (2) yields
\[
L^{X}_{s,k,\iota}(Y)=L^{X}_{s,k,\iota}(Y'')-\beta_{m''-1,t''}%=S'_{m''-1,t''}L^{X}_{s,k,\iota}(Y'').
\]
and by the minimality of $(m'',t'')$, the coefficient of $x_{m''-1,t''}$ in $L^{X}_{s,k,\iota}(Y'')$ is positive. Hence,
$L^{X}_{s,k,\iota}(Y)=S'_{m''-1,t''}L^{X}_{s,k,\iota}(Y'')$.
Combining with $L^X_{s,k,\iota}(Y'')\in \Xi_{s,k,\iota}'^+$, one obtains $L^X_{s,k,\iota}(Y)\in \Xi_{s,k,\iota}'^+$. \qed

\end{document}